\documentclass{amsart}

\input xy
\xyoption{all}

\usepackage{geometry}
\usepackage{amsmath}

\usepackage{amssymb}
\usepackage{graphicx}
\usepackage[noadjust]{cite}
\usepackage{hyperref}
\usepackage{setspace}
\usepackage{amsthm}  % see geometry.pdf on how to lay out the page. There's lots.
\usepackage{quiver}

\newtheorem{same}{This should never appear}[section]
\newtheorem{defin}[same]{Definition}
\newtheorem{definition}[same]{Definition}
\newtheorem{claim}[same]{Claim}

\newtheorem{remark}[same]{Remark}
\newtheorem{theorem}[same]{Theorem}
\newtheorem{example}[same]{Example}
\newtheorem{lemma}[same]{Lemma}
\newtheorem{fact}[same]{Fact}
\newtheorem{question}[same]{Question}
\newtheorem{cor}[same]{Corollary}
\newtheorem{corollary}[same]{Corollary}
\newtheorem{prop}[same]{Proposition}

\newtheorem{hyp}[same]{Hypothesis}

\numberwithin{equation}{section}

\newbox\noforkbox \newdimen\forklinewidth
\setbox0\hbox{$\textstyle\smile$}
\setbox1\hbox to \wd0{\hfil\vrule width \forklinewidth depth-2pt
 height 10pt \hfil}
\setbox\noforkbox\hbox{\lower 2pt\box1\lower 2pt\box0\relax}
\def\unionstick{\mathop{\copy\noforkbox}\limits}

\def\nonfork_#1{\unionstick_{\textstyle #1}}

\setbox0\hbox{$\textstyle\smile$}
\setbox1\hbox to \wd0{\hfil{\sl /\/}\hfil}
\setbox2\hbox to \wd0{\hfil\vrule height 10pt depth -2pt width
              \forklinewidth\hfil}
\newbox\doesforkbox
\setbox\doesforkbox\hbox{\lower 2pt\box1 \lower 2pt\box2\lower2pt\box0\relax}
\def\nunionstick{\mathop{\copy\doesforkbox}\limits}

\def\fork_#1{\nunionstick_{\textstyle #1}}

\newcommand{\LS}{\text{LS}}

\newcommand{\ba}{\bold{a}}
\newcommand{\bb}{\bold{b}}

\newcommand{\bm}{\bold{m}}

\newcommand{\bj}{\bold{j}}
\newcommand{\K}{\mathcal{K}}
\newcommand{\bx}{\bold{x}}

\newcommand{\cf}{\text{cf }}
\newcommand{\rest}{\upharpoonright}
\newcommand{\id}{\textrm{id}}

\newcommand{\ds}{\text{ds}}

\newcommand{\Mod}{\te{Mod }}

\newcommand{\tp}{tp}

\newcommand{\seq}[1]{\langle #1 \rangle}
\newcommand{\te}[1]{\textrm{#1}}

\newcommand{\bg}{\bold{g}}

\newcommand{\bi}{\bold{i}}

\newcommand{\bK}{\mathbb{K}}
\newcommand{\bL}{\mathbb{L}}
\newcommand{\bQ}{\mathbb{Q}}

\newcommand{\bZ}{\mathbb{Z}}

\newcommand{\cK}{\mathcal{K}}
\newcommand{\cL}{\mathcal{L}}

\newcommand{\cP}{\mathcal{P}}

\newcommand{\gS}{\text{gS}}

\newcommand{\cl}{\text{cl}}
\newcommand{\cof}{\text{cof }}

\newcommand{\Card}{Card}
\renewcommand{\bg}{\textbf{big}}
\renewcommand{\big}{\text{{\bf big} }}

\renewcommand{\S}{\text{S}}

\renewcommand{\subset}{\subseteq}

\newcommand{\footnotei}[1]{} \usepackage[final]{showkeys}
\newcommand{\END}{\bibliographystyle{amsalpha}\bibliography{bib}\end{document}}  \newcommand{\comment}[1]{}

\title{Erd\H{o}s-Rado Classes}
\author{Will Boney}
\email{wb1011@txstate.edu}
\address{Department of Mathematics, Texas State University, San Marcos, TX, USA}
\date{\today\\Keywords: nonelemenary classes, Ehrenfeucht-Mostowski models, Erd\H{o}s-Rado theorem\\AMS 2010 Subject Classification: 03C45, 03E02(primary) 03C55, 3C75 (secondary)\\The author was supported by the National Science Foundation under grant DMS-2339018.}

\begin{document}

\begin{abstract}
We amalgamate two generalizations of Ramsey's Theorem--Ramsey classes and the Erd\H{o}s-Rado Theorem--into the notion of a combinatorial Erd\H{o}s-Rado class.  These classes are closely related to Erd\H{o}s-Rado classes, which are those from which we can build generalized indiscernibles and blueprints in nonelementary classes, especially Abstract Elementary Classes.  We give several examples and some applications.
\end{abstract}

\maketitle

\footnotei{For future editing, remember that
\begin{itemize}
	\item \bg is without the space: `\bg'
	\item \big is with the space: `\big'
\end{itemize}
Remove comments!
}

\section{Introduction}

The motivation for this paper is to amalgamate two distinct generalizations of the classic Ramsey's Theorem.  Ramsey's Theorem \cite[Theorem A]{r-ramsey} says that, fixing finite $n$ and $c$ in advance, one can find large, finite homogeneous subsets for colorings of $n$-sized sets with $c$ colors, as long as the set originally colored was big enough.  In the well-known arrow notation\footnote{The notation
$$\alpha \xrightarrow{} (\beta)^r_\gamma$$
means that for any coloring $c:[\alpha]^r \to \gamma$, there is $X \subset \alpha$ of type $\beta$ such that $c"[X]^r$ is a single element; such an $X$ is called \emph{homogeneous}.  Hajnal and Larson \cite[p. 130]{handbook-st-hajnal-larson} point out ``[t]here are cases in mathematical history when a
well-chosen notation can enormously enhance the development of a branch
of mathematics and a case in point is the ordinary partition symbol." }, this can be stated as follows.

\begin{fact}[Ramsey]
For any finite $k, n, c$, there is finite $R$ such that
$$R\xrightarrow{}(k)^n_c$$
\end{fact}

There are two ways for this to be generalized.  The first generalization is to coloring other classes of structures.  An important observation is that coloring subsets of a given set is the same as coloring increasing tuples of that length according to some fixed linear order, so Ramsey's Theorem can be seen as a result about coloring linear orders and finding homogeneous copies of linear orders within it.  A Ramsey class $\cK_0$ is a class of finite structures where a variant of Ramsey's Theorem holds: given finite $k < \omega$ and $A, B \in \cK_0$, there is some $C \in \cK_0$ such that any coloring of the copies of $B$ appearing in $C$ by $k$ colors gives rise to a copy of $A$ in $C$ that is homogeneous for this coloring.  This is written as 
$$C \xrightarrow{} (A)^B_k$$
Independently, Ne\v{s}et\v{r}il and R\"{o}dl \cite{nr-ramsey} and Abramson and Harrington \cite{ah-ramsey} showed that the class of finite, linearly ordered $\tau$-structures is a Ramsey class when $\tau$ is a finite relational language.  Since then the theory of Ramsey classes has become a productive area connecting combinatorics, dynamics, and model theory (the connection to model theory is partially explained below; a nice survey on Ramsey classes is Bodirsky \cite{b-ramsey-survey}).

The second generalization is to remove the restriction `finite' in the statement of Ramsey's Theorem.  Allowing the arity of the coloring (the upper exponent in the arrow relation) to be infinite would make positive results contradict the axiom of choice (see \cite[Theorem 12.1]{ehmr-cst}), so we focus on finite arity colorings.  Ramsey's Theorem can be easily generalized to $\omega \xrightarrow{} (\omega)^n_c$ for all finite $n, c$.  Moving to infinitely many colors and uncountable homogeneous sets, Erd\H{o}s and Rado \cite[Theorem 39]{er-erdosrado} proved the following (and, unlike most bounds in finite Ramsey theory, the left-hand cardinal is known to be optimal).

\begin{fact}[Erd\H{o}s-Rado]
For any finite $n$ and infinite $\kappa$,
$$\beth_{n-1}(\kappa)^+ \xrightarrow{} (\kappa^+)^n_\kappa$$
\end{fact}

This has been generalized in many directions, including unbalanced and polarized partition relations.  Excellent surveys can be found in Erd\H{o}s, Hajnal, M\'{a}t\'{e}, and Rado \cite{ehmr-cst} and Hajnal and Larson \cite{handbook-st-hajnal-larson}.

We give a general framework for generalizations of the Erd\H{o}s-Rado Theorem along the lines of Ramsey classes, appropriately called combinatorial Erd\H{o}s-Rado classes (Definition \ref{cer-def}, see later in this introduction for a discussion of Erd\H{o}s-Rado classes).  Roughly, a class $\cK$ is a combinatorial Erd\H{o}s-Rado class if it satisfies enough instances of $\lambda \xrightarrow{\cK} (\kappa)^n_\kappa$, where this means any coloring of $n$-tuples from any $\lambda$-big structure in $\cK$ with $\kappa$-many colors has a homogeneous substructure that is $\kappa$-big (Section \ref{spr-sec} makes these notions of `big' and `homogeneous' precise).  Note that we require that all $n$-tuples are colored, rather than coloring copies of a single structure as in Ramsey classes.  Many partition relations of this sort (positive and negative) already exist in the literature, and we collect the most relevant and place them in this framework in Section \ref{spr-example-ssec}.

Our main interest in these results comes from model theory, specifically building generalized indiscernibles in nonelementary classes.  An (order) indiscernible sequence indexed by a linear order $I$ is a sequence  $\{\ba_i : i \in I\}$ in a structure $M$ where the information (specifically, the type) about the elements $\ba_{i_1}, \dots, \ba_{i_n}$ computed in the structure $M$ only depends on the ordering of the indices $i_1, \dots, i_n$.   Generalized indiscernibles replace the linear orders with some other index class: trees, functions spaces, etc. Generalized indiscernibles (and the related notion of generalized blueprints) appear in Shelah \cite{sh:c}, and we recount the definitions in Section \ref{prelim-sec}.

In elementary classes (those axiomatizable in first-order logic), indiscernibles exist because of Ramsey's Theorem and compactness.  Moving to more complicated index classes $\cK$, the combinatorics necessary to build generalized indiscernibles from $\cK$ are exactly the same as requiring that $\cK$ be the directed colimits of a Ramsey class $\cK_0$ (see \cite[Theorem 4.31]{s-nip-collapse}).  In both of these constructions, restriction to finite structures is sufficient to build indiscernibles because the compactness theorem reduces satisfiability to satisfiability of finite sets.

The study of nonelementary classes typically focuses on those axiomatizable in nice logics beyond first-order and, slightly more broadly, on Abstract Elementary Classes.  Abstract Elementary Classes (introduced by Shelah \cite{sh88}) give an axiomatic framework for a class of structures $\bK$ and a strong substructure notion $\prec_\bK$ meant to encompass a wide variety of nonelementary classes.  A key feature of nonelementary classes is that they lack the structure that the compactness theorem endows on elementary classes.  Indeed, Lindstr\"{o}m's Theorem \cite{l-lindstrom} says that no logic stronger than first-order can satisfy the classical (countable) compactness theorem and the downward L\"{o}wenheim-Skolem property.  In practice, stronger logics tend to fail compactness (the cofinality quantifier logics $\bL(Q^{\cof}_\alpha)$ are a notable exception \cite{sh43, b-cofquant}).  Thus, different methods are necessary to build indiscernibles in Abstract Elementary Classes.

For order indiscernibles, this method comes by way of Morley's Omitting Types Theorem \cite{m-mott} using the Erd\H{o}s-Rado Theorem mentioned above (an exposition appears in \cite[Chapter 4 and Appendix A]{baldwinbook}).  For generalized indiscernibles, the generalization of the Erd\H{o}s-Rado Theorem to combinatorial Erd\H{o}s-Rado classes described above gives the desired tools.  We call $\cK$ an Erd\H{o}s-Rado class if we can build $\cK$-indiscernibles in any Abstract Elementary Class (Definition \ref{er-def} and Theorem \ref{gmott-ap-thm}).  Generalized indiscernibles have occasionally seen use in nonelementary classses (for instance, \cite{gs-unsuperstab}, \cite[Chapter V.F]{sh:h}, \cite{bs-atomic-stability}).  

The use of structural partition relations allows us to present a unified framework for generating generalized indiscernibles in nonelementary classes.  This allows us to generalize Morley's result as Generalized Morley's Omitting Types Theorem \ref{gmott-thm}.  There is also some work in this direction in Shelah \cite{she59}, and we compare them in Remark \ref{er-ram-remark}.

We also make explicit category theoretic formulations of (generalized) Ehrenfeucht-Mostowski models and indiscernible collapse.  This is motivated by a statement of Morley's Omitting Types Theorem by Makkai and Par\'{e} in their work on accessible categories \cite[Theorem 3.4.1]{mp-accessible}.  Essentially, generalized blueprints correspond to nice functors, and we prove a converse to this as well (Theorem \ref{der-bp-thm}).

We would like to thank Sebastian Vasey, Lynn Scow, and a very thorough referee for comments on initial drafts.  We thank Saharon Shelah for pointing out the argument for Proposition \ref{clo-part-prop}, although he suggests it is well-known.

\subsection{Outline}

Section \ref{prelim-sec} gives the necessary preliminaries on abstract classes of structures, types, and generalized indiscernibles and blueprints.  We also include a description in Section \ref{prelim-example-ssec} of the examples we will consider in this paper.  Section \ref{spr-sec} gives the definition of combinatorial Erd\H{o}s-Rado classes and of the structural partition relation that defines them.  Section \ref{spr-example-ssec} gives several known (and a few new) examples and counterexamples of these classses.  Section \ref{er-sec} defines Erd\H{o}s-Rado classes and proves the main link between the two notions, Generalized Morley's Omitting Types Theorem \ref{gmott-thm}.  Section \ref{ext-sec} describes several extensions and partial converses to this result, including the category theoretic perspective on blueprints.  Section \ref{app-sec} gives three applications of this technology: stability spectra of tame AECs, indiscernible collapse in nonelementary classes, and the interpretability order.

Note that the definition of Erd\H{o}s-Rado classes (Definition \ref{er-def}) does not actually depend on that of combinatorial Erd\H{o}s-Rado classes (Definition \ref{cer-def}) or any of Section \ref{spr-sec}.  However, we give the combinatorial definitions first, as they provide the largest class of examples of Erd\H{o}s-Rado classes.

\subsection{Conventions} Throughout the paper, we deal with different classes of structures, normally referred to by $K$ in some font with some decoration.  To aid the reader, we observe the following convention:
\begin{itemize}
	\item the script or calligraphic $K$--typeset as $\cK$--will be used as the domain or index class that we wish to build generalized indiscernibles \emph{from}.  They typically have few assumptions of model-theoretic structure on them.  Erd\H{o}s-Rado classes will be of this type, and the class of linear orders form the prototypical example.
	\item the bold $K$--typeset as $\bK$--will be used as the target class that we wish to build generalized indiscernibles \emph{in}.  They will typically be well-structured in some model-theoretic sense.  Elementary classes and Abstract Elementary Classes form the prototypical examples.
\end{itemize}
We also observe two important conventions with respect to types that might be missed by the model-theoretically inclined reader that skips the Preliminaries Section (see Definition \ref{stone-def}):
\begin{enumerate}
	\item Since we never deal with types over some parameter set\footnote{Indiscernibles over a set of parameters can be recovered by adding those parameters to the language. The exception to this rule is Section \ref{unsuper-ssec}, which deals with a specific application to Abstract Elementary Classes, and uses types over parameter sets and other techniques.}, we omit the domain of types throughout.  For example, we write $\tp_{\cK}(\ba; I)$ for the $\cK$-type of $\ba$ over the empty set computed in $I$, rather than $\tp_{\cK}(\ba/\emptyset; I)$; and
	\item $\bK^\tau$ is the class of all $\tau$-structures with $\tau$-substructure $\subset_\tau$ as the strong substructure relation.  In particular, $\tp_\tau$ is the type in this class, which turns out to be quantifier-free type (Proposition \ref{qf-type-prop}).
\end{enumerate}

Note that Sections \ref{gsott-ssec} and \ref{unsuper-ssec} require more knowledge about Abstract Elementary Classes.  This can be found in, e.g., Baldwin \cite{baldwinbook}.

\section{Preliminaries}\label{prelim-sec}

\subsection{Classes of structures and types}

We want to have a very general framework for classes of structures in a common language along with a distinguished substructure relation.  Although much more general than we need, we can use the notion of an abstract class (this formalization is originally due to Grossberg).  Additionally, we expect our Erd\H{o}s-Rado classes to have orderings (similar to, e.g., \cite[Proposition 2.2]{b-ramsey-survey} for Ramsey classes), so we introduce the notion of an ordered abstract class.  An alternative would be to consider equivalence classes of types in the Stone space after modding out by permutation of the indices, but requiring an ordering seems simpler.

Note that the examples presented tend to be universal classes (and mostly relational), in which case the type is determined by the quantifier-free type in $\bL_{\omega, \omega}$.  However, we offer a more general framework because it adds little technical difficulty and offers the possibility to wider applicability.  For instance, well-founded trees (Example \ref{wftr-prelim-ex}, \ref{wftr-spr-ex}) are {\bf not} a universal class.

\begin{defin}\
\begin{enumerate}
	\item $(\cK, \leq_\cK)$ is an \emph{abstract class} iff there a language $\tau = \tau(\cK)$ such that each $M \in \cK$ is a $\tau$-structure, $\leq_\cK$ is a partial order contained in $\subset_\tau$, and membership in $\cK$ and $\leq_\cK$ both respect isomorphism.  We often refer to the class simply as $\cK$.
    \begin{enumerate}
        \item In an abstract class $(\cK, \leq_\cK)$, a \emph{$\cK$-embedding} is $\tau$-embedding (an injection that preserves and refelects atomic sentences) $f:M \to N$ between elements of $\cK$ such that $f(M) \leq_\cK N$.
    \end{enumerate}
	\item $(\cK, \leq_\cK)$ is an \emph{ordered abstract class} iff it is an abstract class with a distinguished binary relation $<$ in $\tau(\cK)$ such that $<^I$ is a total order of $I$ for every $I \in \cK$.
\end{enumerate}
\end{defin}

We will also use the types of elements.  Most of the classes we consider will not be elementary (either in axiomatization of $\cK$ or ordering $\leq_\cK$), so syntactic types give way to semantic notions.  Specifically, we use the notion of Galois types (also called semantic or orbital types) used in the study of Abstract Elementary Classes (and originated in \cite{sh300}).  However, in most cases, this will be the same as quantifier-free types.  Note that we typically drop any adjective and use `type'  or sometimes `$\cK$-type' to refer to the following semantic definition, although we will decorate the symbol with the ambient class.  Although we use the `index class' notation $\cK$ throughout Definition \ref{stone-def}, we will use these ideas for both index classes and target classes.

\begin{defin} \label{stone-def}
Let $\cK$ be an abstract class.
\begin{enumerate}
	\item Given $I_1, I_2 \in \cK$ and $\ba_1 \in I_1$, $\ba_2 \in I_2$, we say that $\ba_1$ and $\ba_2$ \emph{have the same $\cK$-type} iff there are $J_1, \dots, J_n; I^*_1, \dots, I^*_{n+1} \in \cK$, $\bb_\ell \in I^*_\ell$, and $\cK$-embeddings $f_\ell:I^*_{\ell+1}\to J_\ell$ such that
	\begin{enumerate}
		\item $I^*_1=I_1$, $I^*_{n+1}=I_2$, $\ba_1 = \bb_1$, and $\ba_2=\bb_{n+1}$;
		\item $I^*_\ell \leq_\cK J_\ell$; and
		\item $\bb_\ell = f_\ell(\bb_{\ell+1})$.
	\end{enumerate}
	$$
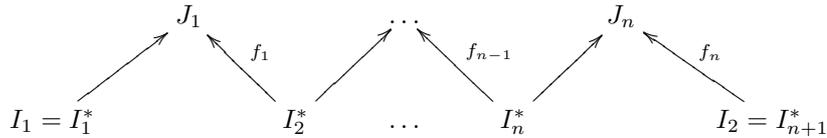
\xymatrix{ & J_1 &  & \dots &  & J_n & \\
	I_1=I_1^*\ar[ur] &  & I_2^* \ar[ul]_{f_1} \ar[ur] & \dots & I_n^*\ar[ul]_{f_{n-1}} \ar[ur] &  & I_2 = I_{n+1}^*\ar[ul]_{f_n}}$$
	
	We write $tp_\cK(\ba; I)$ to be the equivalence class\footnote{This is a proper class, but we can use Scott's trick (see \cite[p. 65]{jech}) or some other method to only deal with sets.} of all tuples in all structures that have the same type as $\ba$.  Thus, `$tp_\cK(\ba_1; I_1) = tp_\cK(\ba_2; I_2)$' has the same meaning as `$\ba_1$ and $\ba_2$ have the same $\cK$-type.'
	\item $\S_\cK:=\left\{ tp_\cK(\ba; I) \mid \ba \in I\in \cK\right\}$ is the \emph{Stone space} or space of types.
	\item If $\cK$ is an ordered abstract class, then $\S_\cK^{inc}$ is the subset of $\S_\cK$ whose realizations are in increasing order, namely,
	$$\S^{inc}_\cK:=\left\{ tp_\cK(\ba; I) \mid \ba \in I\in \cK \text{ and }a_1 < \dots <a_n\right\}$$
	\item Adding a superscript $n < \omega$ to either $\S_\cK$ or $\S^{inc}_\cK$ restricts to looking at types of $n$-tuples.
	\item \label{type-rest-item} Let $p \in \S_\cK^n$ be $\tp_\cK(i_1, \dots, i_n; I)$ and $s \subset n$ be the set of $k_1 < \dots < k_m$ for $m = |s|$.  Then $p^s := \tp_\cK(i_{k_1}, \dots, i_{k_m}; I) \in \S_\cK^m$.
	\item If we have an ordered abstract class decorated with a superscript $\cK^x$, then we often use this superscript in place of the whole class in this notation, e.g., the Stone space of $\cK^{\chi-or}$ is denoted $\S_{\chi-or}$ rather than $\S_{\cK^{\chi-or}}$.  
    \item For a language $\tau$, we use $\bK^\tau$ to be the abstract class of all $\tau$-structures with $\tau$-substructure.
\end{enumerate}
\end{defin}

Understanding this notation is key to the rest of the paper, so we unravel these notions in some examples.

First, consider $\bK^\tau$.  $\bK^\tau$-embeddings are injections that preserve and reflect the $\tau$-structure. Then $\S_\tau$ refers to all $\bK^\tau$-types in this class, and these types turn out to be exactly quantifier-free types in the language.

\begin{prop}\label{qf-type-prop}
Fix a language $\tau$.  Then
\begin{enumerate}
    \item if $M_0, M_1 \in \bK^\tau$ are structures and $\ba_\ell \in M_\ell$ such that
    $$\tp_{qf}(\ba_0; M_0) = \tp_{qf}(\ba_1; M_1)$$
    then the $\bK^\tau$-types are also equal, in particular witnessed by $\bK^\tau$-embeddings $f_\ell: M_\ell \to N$ such that $f_0(\ba_0)=f_1(\ba_1)$; and
    \item $\bK^\tau$ has amalgamation.
\end{enumerate}
\end{prop}

Note that although we say `$\bK^\tau$-types are quantifier-free types', the statement
$$\tp_\tau\left(\ba; M\right) = \tp_{qf}\left(\ba;M\right)$$
is technically false: the left-hand side is a (rank-initial subset of) equivalence class of pairs, while the right hand side is a $|\tau|+\aleph_0$-sized collection of quantifier-free formulas.

{\bf Proof:}  The second item is well-known and can be proved in many ways.  One way is to take $M_0\subset M_1, M_2$ and freely generate the $\tau$-structure on $M_1 \cup M_2$ by closing under functions and adding no relations; this even proves a disjoint version of amalgamation.  The second item is proved similarly after first identifying $\ba_0$ and $\ba_1$ in the union; this is possible exactly because they share the same quantifier-free type.\comment{For the first item, write down the first order (and quantifier free) theory that is the union of the atomic diagram of $M_0$ and $M_1$ using the same constants for $\ba_0$ and $\ba_1$.  By compactness, if this is unsatisfiable, then it witnessed by $\phi_\ell$ satisfied in $M_\ell$ such that 
$$\vDash \phi_0 \to \neg \phi_1$$
By Craig's Interpolation Theorem, there is $\psi$ in the common language of $\phi_0$ and $\phi_1$--which is $\tau$ and the constant symbols for $\ba_\ell$--that interpolates this.  But $\psi$ is part of each atomic diagram, so we cannot have $\psi \to \neg \phi_1$.  Thus, this union of atomic diagrams is satisfiable, and that model witnesses the type equality as described.

For the second item, nothing in the above proof use that $\ba_\ell$ was finite.  So given $M_0$ that is a $\tau$-substrucutre of $M_1$ and $M_2$, this means
$$\tp_{qf}(M_0; M_1) = \tp_{qf}(M_0; M_2)$$
The same argument above gives an amalgam of these models. }\hfill \dag\\

Note that $\bK^\tau$ might fail to have the joint embedding property.  In fact, this property is equivalent to $\bK^\tau$ having a unique $0$-type, and is satisfied whenever $\bK^\tau$ has no constants (or $0$-ary functions).

As a second example, consider the class $\cK^{2-or}$ (which will be used in Example \ref{gem-ex}).  This is defined more generally in Example \ref{chior-prelim-ex}, but this case consists two disjoint linear orders $(I_0, I_1)$ where everything in $I_0$ is less than everything in $I_1$.  Then $\S_{2-or}$ is the collection of all $\bK^{2-or}$-types.  In particular, an increasing tuple $\ba$ from a model $(I_0, I_1)$ in $\bK^{2-or}$ can be written as $\ba_0, \ba_1$ where each $\ba_\ell$ is an increasing tuple (possibly empty) from $I_\ell$.  A proof similar to (and simpler than) Proposition \ref{qf-type-prop} shows that $\bK^{2-or}$-types are also quantifier free types.  Thus, an element $p \in \S^{inc, n}_{2-or}$ is determined by some number $k\leq n$ that indicates that the tuple is an increasing $k$-tuple from the first linear order, followed by an increasing $(n-k)$-tuple from the second part.

A similar statement is true for $\cK^{\chi-or}$, although a type $p \in \S^{inc, n}_{\chi-or}$ is determined by a partition of $n$ into $\chi$-many pieces, most of which are empty when $\chi$ is infinite.

\subsection{Generalized indiscernibles and blueprints}

The following generalizes the normal theory of blueprints and Ehrenfeucht-Mostowski models begun in \cite{em-emmodels}.  These generalized notions appear in \cite[Section VII.2]{sh:c}.

\begin{definition} \label{genbp-def} Let $\cK$ be an ordered abstract class.
\begin{enumerate}
    \item A \emph{blueprint $\Phi$ proper for $\cK$} is a function $\Phi:\S_\K^{inc} \to \S_\tau$ for some $\tau = \tau(\Phi)$ that satisfies the following coherence conditions:
    \begin{enumerate}
        \item the free variables of $\Phi(p)$ are the free variables of $p$; and
        \item given variables $s \subset n$ and $p \in \S_\cK^{inc, n}$, we have that 
        $$\Phi\left(p^s\right) = \Phi(p)^s$$
    \end{enumerate}
    
    $\Upsilon^{\K}$ is the collection of all blueprints proper for $\K$.
    
    $\Upsilon^{\K}_\kappa$ is the collection of all blueprints proper for $\K$ such that $|\tau(\Phi)|\leq \kappa$.
    \item \label{genem-item} Let $I \in \cK$ and $\Phi$ be a blueprint proper for $\cK$.  Then, we can build a $\tau(\Phi)$-structure  $EM(I, \Phi)$ such that, for all $i_1< \dots < i_n \in I$, we have that 
    $$\tp_\tau\left(i_1, \dots, i_n; EM(I, \Phi)\right) = \Phi\left(\tp_\cK(i_1, \dots, i_n; I)\right)$$
     and that every element of $EM(I, \Phi)$ is a $\tau(\Phi)$-term of a sequence from $I$.
    
    If $\tau\subset \tau(\Phi)$, then $EM_\tau(I, \Phi) := EM(I, \Phi)\rest \tau$.
    \item Given an class $K$ of $\tau$-structures and a blueprint $\Phi$  with $\tau \subset \tau(\Phi)$, we say that \emph{$\Phi$ is proper for $(\cK, K)$} iff it is proper for $\cK$ and, for any $I \in \cK$, $EM_{\tau}(I, \Phi) \in K$.
    
    $\Upsilon^{\cK}[K]$ is the collection of all blueprints proper for $(\cK, K)$.
    
    $\Upsilon^{\cK}_\kappa[K]$ is the collection of all blueprints proper for $(\cK, K)$ such that $|\tau(\Phi)|\leq \kappa$.
    
    \item Given an abstract class $\bK=(K, \prec_\bK)$ and a blueprint $\Phi$ with $\tau(\bK) \subset \tau(\Phi)$, we say that \emph{$\Phi$ is proper for $(\cK, \bK)$} iff it is proper for $(\cK, K)$ and, for any $I \leq_\cK J \in \cK$, $EM_{\tau(\bK)}(I, \Phi) \leq_\bK EM_{\tau(\bK)}(J, \Phi)$.
    
    $\Upsilon^{\cK}[\bK]$ is the collection of all blueprints proper for $(\cK, \bK)$.
    
    $\Upsilon^{\cK}_\kappa[\bK]$ is the collection of all blueprints proper for $(\cK, \bK)$ such that $|\tau(\Phi)|\leq \kappa$.
\end{enumerate}
\end{definition}

Being proper for $\cK$ is the same as being proper for $(\cK, \bK^{\tau(\Phi)})$.  Note that if $\Phi \in \Upsilon^{\cK}[\bK]$, then the blueprint $\Phi$ actually maps $\S^{inc}_{\cK} \to S_\bK$.  An observant reader might complain that the description in Definition \ref{genbp-def}.(\ref{genem-item}) uniquely describes a model, but is short on proving it's existence.  However, the existence of such a model follows from standard arguments about EM models, see, e.g., \cite[Section 5.2]{m-modeltheory}.  Our formalism has $I$ be the generating set for $EM(I, \Phi)$ (and later indiscernibles), rather than passing to a skeleton.

From a category-theoretic perspective, a blueprint $\Phi \in \Upsilon^{\cK}[\bK]$ induces a functor $\Phi:\cK \to \bK$ that is faithful, preserves direceted colimits, and induces a natural transformation between the `underlying set' functor of each concrete category.  We return to this perspective in Section \ref{ct-ind-ssec} and derive a converse Theorem \ref{gmott-con-thm} of the Generalized Morley's Omitting Types Theorem \ref{gmott-thm}.

\begin{example} \label{gem-ex}\
\begin{enumerate}
    \item These definitions generalize the standard notions of blueprints and Ehrenfeucht-Mostowski models when $\K$ is the class of linear orders.
    \item Consider a bidimensional theory like the theory $T$ of a predicate $P$ that is infinite and coinfinite\footnote{A more mathematically complex example of a bidimensional theory is the theory $Th( \oplus \bZ(p^\infty))$ of the direct sum of countably many copies of the Pr\"{u}fer p-group.  The same analysis applies there.}.  Each model is determined by two infinite cardinals, the size of $P$ and the size of its complement. Using standard Ehrenfeucht-Mostowski models, one could only get blueprints that either vary one dimension and not the other {\bf or} make the dimensions the same.
       
    However, there is a generalized blueprint\footnote{$\Upsilon^{2-or}$ is $\Upsilon^{\cK^{2-or}}$, and this class is consists of two disjoint linear orders; see Example \ref{chior-prelim-ex}.} $\Phi \in \Upsilon^{2-or}_{\aleph_0}[\text{Mod}(T)]$ for the class of two disjoint linear orders that takes $(I, J)$ to the model $M_{(I, J)}$ with universe 
    $$(I+\omega)\cup (J+\omega)$$
    and predicate $P^{M_{(I, J)}}:= I+\omega$.
    Thus, every model of $T$ is isomorphic to $EM_\tau\left((I, J), \Phi\right)$ for some $I$ and $J$.
\end{enumerate}
\end{example}

Using generalized blueprints, we can build models with generalized indiscernibles (see Theorem \ref{gmott-ap-thm} for this in action).  

\begin{defin}\label{gen-ind-def}
Let $\cK$ an ordered abstract class and $\bK$ be an abstract class.  Then, given $I \in \cK$ and $M \in \bK$, a collection $\{\ba_i \in {}^{<\omega}M \mid i \in I\}$ is a \emph{$\cK$-indiscernible sequence} iff for every $i_1< \dots< i_n; j_1< \dots< j_n \in I$, if 
$$\tp_\cK(i_1, \dots, i_n; I) = \tp_\cK(j_1, \dots, j_n; I)$$
then
$$\tp_\bK(\ba_{i_1}, \dots, \ba_{i_n}; M) = \tp_\bK(\ba_{j_1}, \dots, \ba_{j_n}; M)$$
\end{defin}

An important fact to keep in mind is that, in nonelementary classes, not every collection of indiscernibles can be turned into a blueprint; \cite[Example 18.9]{baldwinbook} provides such an example.  This is in contrast to first-order, where every infinite set of indiscernibles can be stretched (see \cite[Lemma 5.1.3]{tz-model-theory}).

\subsection{Our examples}\label{prelim-example-ssec}
There will be several examples that we will develop here and in Section \ref{spr-example-ssec}.  Here, we define the relevant classes and note the syntactic characterization of their types (normally quantifier-free).  Section \ref{spr-example-ssec} explains how these classes fit within the framework of Erd\H{o}s-Rado classes.  In each case, the strong substructure relation is just substructure for the appropriate language unless noted otherwise.

\begin{example}[Linear orders] \label{or-prelim-ex}
$\cK^{or}$ is the class of linear orders in the language with a single binary relation $<$.  This is an ordered abstract class and is universal, so $\cK^{or}$-type is simply quantifier-free type.  This is our prototypical Erd\H{o}s-Rado class.
\end{example}

\begin{example}[$\chi$ disjoint linear orders] \label{chior-prelim-ex}
$\cK^{\chi-or}$ is the class of $\chi$ disjoint linear orders.  In order to make this an ordered abstract class, we say $\bar{I} \in \cK^{\chi-or}$ consists of disjoint sets $\{I_i\}_{i<\chi}$ and a total ordering $<$ such that $i < j < \chi$ implies that $I_i << I_j$ ($X<<Y$ means that every element of $X$ is below every element of $Y$).   Note that if $\chi$ is infinite, then this is not an elementary class.
\end{example}

\begin{example}[$\chi$-colored linear orders] \label{clo-prelim-ex}
We set $\cK^{\chi-color}$ to be a particular class of colored linear orders.  $(I, <, P_\beta)_{\beta<\chi} \in \cK^{\chi-color}$ consists of a well-ordering $(I, <)$ such that $P_\beta = \{i \in I : i\text{ is the }(\chi\cdot\gamma+\beta)\text{th element of $I$ for some }\gamma\}$.
\end{example}

\begin{example}[Trees of height $n<\omega$]\label{ntr-prelim-ex}
Fix the language $\tau_{n-tr} = \left( P_k,<, \prec, \wedge\right)_{k<n}$.  Then $\cK^{n-tr}$ consists of all $\tau_{n-tr}$-structures $I$ such that
\begin{itemize}
	\item $(I, \prec)$ is a tree of height $n$;
	\item $P_k$ are all vertices on level $k$;
	\item $<$ is a total order of $I$ coming from a lexicographic ordering of the tree; and
	\item $\wedge$ is the meet operation on this tree.
\end{itemize}
Then $\cK^{n-tr}$-type is just quantifier-free type in this language.
\end{example}

\begin{example}[Trees of height $\omega$]\label{wtr-prelim-ex}
$\cK^{\omega-tr}$ are the trees of height $\omega$ formalized in the language $\tau_{\omega-tr}=\cup_{n<\omega}\tau_{n-tr}$.
\end{example}

Of course, these tree examples can be continued on past height $\omega$, but we know of no results (positive or negative) on these classes in terms of the Erd\H{o}s-Rado notions.

\begin{example}[Well-founded trees]\label{wftr-prelim-ex}
$\cK^{wf-tr}$ are the well-founded trees formalized in the language $\tau_{\omega-tr}$; recall a tree is well-founded iff it contains no infinite branch.  Given any ordinal $\alpha$, we can build a well-founded tree whose nodes are decreasing sequences of ordinals starting with $\alpha$.  Going the other way, if $T$ is a well-founded tree, then we can relabel the nodes with ordinals such that each path is a decreasing sequence.  Note that this relabeling can be done in many different ways and it probably disagrees with the lexicographical ordering on successors in $\tau_{\omega-tr}$.%.Well-founded trees can be identified with decreasing sequences of ordinals.  \footnote{WB: Explain this example more, look at GS reference. Referee \# 31}%Further, we can assign a rank to the nodes of well-founded trees by stipulating that the rank of a node $\eta$ is $\geq \alpha+1$ iff it has a successor of node $\geq \alpha$.  When we allow for arbitrary splitting (as we do in $\cK^{wf-tr}$), any ordinal rank is attainable; if we let $\ds(\alpha)$ be the well-founded tree consisting of all descending sequences from $\alpha$, then the rank of the root $\seq{}$ in this tree is $\alpha$.
\end{example}

\begin{example}[Convexly-ordered equivalence relations]\label{ceq-prelim-ex}
A convexly ordered equivalence relation is $(I, <, E)$, where $E$ is an equivalence relation on $I$, $<$ is a total order, and 
$$\forall x, y, z \in I \left( x E z \wedge x < y < z \to x E y\right)$$
$\cK^{ceq}$ is the collection of all such structures.  These are similar to the class $\cK^{\chi-or}$ except the $\chi$ is allowed to vary.  However, the type of, e.g., singletons in different equivalence classes is the same.  This will make finding type homogeneous sets for colorings more difficult.
\end{example}

\begin{example}[$n$-multi-linear orders]\label{nmlo-prelim-ex}
A $n$-multi-linear order is $(I, <_1, \dots, <_n)$ where each $<_i$ is a linear order of $I$.  $\cK^{n-mlo}$ is the class of these.  We take $<_1$ as the distinguished linear order to view this as an ordered abstract class.
\end{example}

\begin{example}[Ordered graphs]\label{og-prelim-ex}
$\cK^{og}$ consists of the class of all ordered graphs.
\end{example}

\begin{example}[Colored hypergraphs]\label{chg-prelim-ex}
Fix $k \leq \omega$ and a cardinal $\sigma$.  $\cK^{(k, \sigma)-hg}$ consists of all $(I, \sigma; <, F, \alpha)_{\alpha < \sigma}$ where $<$ is a linear ordering and $F:[I]^{<k} \to \sigma$ is a function. If $\sigma=2$, then one can think of $\cK^{(k, 2)-hg}$ as the collection of all hypergraphs with all edge arities $<k$.
\end{example}

\section{Structural Partition Relations and Combinatorial Erd\H{o}s-Rado Classes} \label{spr-sec}

We will formulate a version of the normal partition relation for classes other than linear orders in Definition \ref{spr-def}.  This will encapsulate the idea that any coloring of $n$-tuples from a large structure will have a large substructure that behaves the same way with respect to this coloring.  First, we consider an example that indicates some of the difficulties and the need for new concepts, namely bigness notions (Definition \ref{big-def}) and type-homogeneity (Definition \ref{type-hom-def}).

\begin{example}\label{bad-color-ex}
Let $\bar{I}=(I_0, I_1) \in \cK^{2-or}$ (recall Example \ref{chior-prelim-ex}).  Define a coloring $c:[\bar{I}]^2 \to 2$ based on wether the two elements are in the same partition: given $i, j \in (I_0, I_1)$, set
$$c(\{i, j\}) = \begin{cases}
0 & i \in I_0 \iff j \in I_1\\
1 & \text{otherwise}
\end{cases}$$
Then any $\bar{I^*} \subset \bar{I}$ that contains at least one element from one partition and two from the other will not be homogeneous for this coloring no matter what $\bar{I}$ is.  Alternatively, any $\bar{I^*} \subset \bar{I}$ that contains elements from just one partition will always be homogeneous.
\end{example}

This example exposes two issues.  
\begin{itemize}
	\item First, we could take $\bar{I^*}$ to be $(\emptyset, I_1)$, which is homogeneous for this coloring.  However, taking one of the partitions to be empty goes against the point of working in $\cK^{2-or}$.  So we will attach to these classes a notion of size (or bigness) that takes the structure of the class into account.
	\item Second, we colored the pairs using information about their type.  This meant that we could place restrictions on the structure of any homogeneous subset.  To allow for big homogeneous sets we will allow for the `single color' to depend on the type of tuple.
\end{itemize}

For the first issue, we define abstractly what it means to be a bigness notion.  The only requirements are a monotonicty condition and some weak degree of universality.  For each class from Subsection \ref{prelim-example-ssec}, we make its associated bigness notion explicit in Subsection \ref{spr-example-ssec}.  Many natural bigness notions correspond to a degrees of universality after removing the order, but other cases are more complex (see Examples \ref{clo-part-prop} or \ref{wftr-spr-ex}).

\begin{defin}\label{big-def}
Let $\cK$ be an abstract class.  A \emph{bigness notion} \big for $\cK$ is a class $\{\cK_\mu^\bg \subset \cK \mid \mu \in \Card\}$ such that
\begin{enumerate}
    \item each $\cK_\mu^\big$ is nonempty;
	\item if $\mu_1 \leq \mu_2$ and $M \leq_{\cK} N$, then $M \in \cK^\bg_{\mu_2}$ implies that $N \in \cK^\bg_{\mu_1}$; and
	\item if $M \in \cK_{\aleph_0}^\bg$, then every type in $\S_\cK$ is realized in $M$.
\end{enumerate}
We write `$M \in \cK$ is $\mu$-\bg' for `$M \in \cK_\mu^\bg$.'  Also, we will typically only have one have one bigness notion for a given class, so we will omit it.
\end{defin} 

Note that the omission of \big will lead to some nonstandard notation, e.g., $\cK^{\chi-or}_\mu$ are the at least $\mu$-big elements of $\cK^{\chi-or}$ according to the bigness notion given in Example \ref{chior-spr-ex}, rather than all elements of $\cK^{\chi-or}$ whose universe has cardinality $\mu$.

\begin{remark} \label{0-type-remark}
    Note that the existence of a bigness notion for $\cK$ (or just a model satisfying the conclusion of Definition \ref{big-def}.(3)) implies that there is a unique $0$-type realized by any model in $\cK$.  This is because each model realizes a single $0$-type, and so every model must realize the $0$-type realized in any $\aleph_0$-big model.

    In most examples below, there are no constants, so the $0$-type does not contain information about any `prime structure.' However, the semantic definition of types used here (recall Definition \ref{stone-def}) means that the uniqueness of the $0$-type is equivalent to a natural notion of connectedness of the category $\cK$; this connectedness is a reflection of the completeness of a first-order theory in elementary classes.
\end{remark}

Turning to homogeneity, the key observation from Example \ref{bad-color-ex} was that the types of tuples are extra information that can be used to define a coloring.  In the class of linear orders, there is only one increasing type of an $n$-tuple, so this issue doesn't arise.  In the general case, we can always use the type as information to color a tuple, so we want homogeneity to mean that the type is the \emph{only} information that can be used to determine the color of a tuple.

\begin{defin}\label{type-hom-def}
Let $\cK$ be an ordered abstract class, $I \in \cK$, and $c:[I]^n\to \kappa$.  We say that $I_0 \leq_\cK I$ is \emph{type-homogeneous for $c$} iff the color of a tuple from $I_0$ is determined by the $\cK$-type of that tuple listed in increasing order; that is, there is a function $c^*:\S^{inc,n}_\cK \to \kappa$ such that, for any $i_1 <\dots < i_n \in I_0$, we have that
$$c\left(\{i_1, \dots, i_n\}\right) = c^*\left( \tp_\cK\left(i_1, \dots, i_n; I\right)\right)$$
\end{defin} 

In Example \ref{bad-color-ex}, the entire set $\bar{I}$ is type-homogeneous for the given coloring.

With these new concepts in hand, we can define the structural partition relation.

\begin{defin}\label{spr-def}
Let $\cK$ be an ordered abstract class with a bigness notion \bg.  Given cardinal $\mu, \lambda, \alpha, \kappa$, we write
$$(\lambda) \xrightarrow[\big]{\cK}(\mu)^\alpha_\kappa$$
to mean that given any $\lambda$-\big $I \in \cK$ and coloring $c:[I]^\alpha \to \kappa$, there is a $\mu$-\big $I_0 \leq_\cK I$ from $\cK$ that is type-homogeneous for $c$.

If $\cK$ is one of our examples with an associated bigness notion and is denoted $\cK^x$, then we simply write
$$(\lambda)\xrightarrow{x} (\mu)^\alpha_\kappa$$
for $(\lambda) \xrightarrow[\big]{\cK^x}(\mu)^\alpha_\kappa$.
\end{defin} 

Since the associated bigness notion for $\cK^{or}$ is simply cardinality, $(\lambda) \xrightarrow{or} (\mu)^\alpha_\kappa$ is the normal partition relation.  In particular, positive instances of the structural partition relation are guaranteed by the Erd\H{o}s-Rado Theorem, which states that $\beth_{n-1}(\kappa)^+ \xrightarrow{or} (\kappa^+)^n_\kappa$ for every cardinal $\kappa$ and every $n < \omega$.  We only consider this relation with $\alpha$ finite.\footnotei{WB: Could we reprove the restriction in any ERC?}

Polarized partition relations (see \cite[Section III.8.7]{ehmr-cst}) are similar to $\xrightarrow{\chi-or}$, but typically specify (in our language) the type of the tuple to be considered (and so are more like the Ramsey class-style partition relations, although we use them in Remark \ref{ehr-rem}).

We have focused on the difference between the Ramsey-style structural partition relation (that fixes the type of the tuples colored) and the Erd\H{o}s-Rado-style structural partition relation (that colors all tuples of a fixed length).  These differ on a case-by-case basis, but the referee points out that if there are finitely many $n$-types in the class, then enough Ramsey-style structural partition results can be strung together to get a type-homogeneous coloring.  That is, if $\cK$ is a Ramsey class with finitely many $n$-types for each $n<\omega$, then for each $n, k<\omega$ and $A \in \cK$, there is a $C^* \in \cK$ such that any coloring of $n$-tuples of $C^*$ with $k$ colors has a type-homogeneous copy of $A$.  We also use a version of this argument in Remark \ref{ehr-rem}.

We will list several further positive instances of structural partition relations (new and old) in Subsection \ref{spr-example-ssec}.

From the structural partition relation, we can define combinatorial Erd\H{o}s-Rado classes as those that satisfy structural partition relations for all inputs on the right side.

\begin{defin}\label{cer-def}
Let $\cK$ be an ordered abstract class with a bigness notion \bg.  We say that $\cK$ is a \emph{combinatorial Erd\H{o}s-Rado class} iff there is some function $F:\Card\times \omega \to \Card$ such that, for every $\kappa < \mu$ and $n<\omega$, we have that
$$\left(F(\mu, n)\right) \xrightarrow[\big]{\cK} \left(\mu\right)^n_\kappa$$
We refer to the function $F$ as a witness.
\end{defin}

The primary results and discussions deal with colorings that fix the arity $n$ of the coloring.  However, it is sometimes useful to consider colorings that color all tuples of arity $\leq n$ at once (the `sometimes' in this paper is Theorem \ref{endapp-er-thm}).  It turns out that this stronger relation holds in any combinatorial Erd\H{o}s-Rado class and is witnessed by an iterate of the normal witnessing functions.

\begin{defin}
Let $\cK$ be an ordered abstract class with a bigness notion \bg.  Given cardinal $\mu, \lambda, \alpha, \kappa$, we write
$$(\lambda) \xrightarrow[\big]{\cK}(\mu)^{\leq\alpha}_\kappa$$
to mean that given any $\lambda$-\big $I \in \cK$ and coloring $c:[I]^{\leq\alpha} \to \kappa$, there is a $\mu$-\big $I_0 \leq_\cK I$ from $\cK$ that is type-homogeneous for $c$.

We use the same simplifying notations as in Definition \ref{spr-def}.
\end{defin}

\begin{prop}\label{leqn-prop}
Suppose that $\cK$ is a combinatorial Erd\H{o}s-Rado class witnessed by $F$.  Define $F^*:\Card\times\omega\to \Card$ by induction:
\begin{eqnarray*}
    F^*(\mu,1)&=&F(\mu, 1)\\
    F^*(\mu, n+1) &=&F\left(F^*(\mu, n), n+1\right)
\end{eqnarray*}
Then for each $\kappa <\mu$ and $n<\omega$,
$$\left(F^*(\mu, n)\right) \xrightarrow{\cK} \left(\mu\right)^{\leq n}_{\kappa}$$
\end{prop}
{\bf Proof:} We work by induction on $n<\omega$. For $n=1$, there is nothing to prove.  Suppose this holds for some $n<\omega$.

Let $X \in \cK$ by $F^*(\mu, n+1)$-big and we have a coloring
$$c:[X]^{\leq n+1}\to \kappa$$ 
From definition of $F^*$, we have
$$\left(F^*(\mu, n+1)\right) \xrightarrow{\cK} \left(F^*(\mu, n)\right)^{n+1}_{\kappa}$$
so there is $X_0 \subset X$ that is $F^*(\mu, n)$-big and type-homogeneous for $c$ restricted to the $(n+1)$-tuples.  By induction,
$$\left(F^*(\mu, n)\right) \xrightarrow{\cK} \left(\mu\right)^{\leq n}_{\kappa}$$
There is $X_1 \subset X_0$ that is $\mu$-big and type-homogeneous for $c$ restricted to the $(\leq n)$-tuples.  Then $X_1 \subset X$ is $\mu$-big and type-homogeneous for all of $c$, proving the proposition. \hfill\dag\\

\begin{remark}
\begin{enumerate}    
    \item We could swap the order the definition of the function in Proposition \ref{leqn-prop} by setting 
    \begin{eqnarray*}
    F^{**}(\mu,1)&=&F(\mu, 1)\\
    F^{**}(\mu, n+1) &=&F^{**}\left(F(\mu, n+1), n\right)
\end{eqnarray*}
and the result still holds.  This would be useful if there is an example where $F^{**}$ was less than $F^*$, but it gives the same function in the known examples.
    \item Note that for a cardinal $\lambda$, the following are equivalent:
    \begin{enumerate}
        \item for all $\mu < \lambda$ and $n<\omega$,
        $$F(\mu, n)<\lambda$$
        \item for all $\mu < \lambda$ and $n<\omega$,
        $$F^*(\mu, n) < \lambda$$
    \end{enumerate}
    This observation will be helpful in Theorem \ref{endapp-er-thm} because it means that the `accumulation cardinals' of the two functions are the same.
\end{enumerate}
\end{remark}

\subsection{Examples and some counter-examples} \label{spr-example-ssec}

We show that the examples introduced in Section \ref{prelim-example-ssec} are combinatorial Erd\H{o}s-Rado classes or mention results indicating they are not.  In most cases, no claim of the optimality of the witnessing functions is made.  While interesting from a combinatorial perspective, any reduction of the bounds on the order of `finitely many power set operations' will not affect the witnesses for these classes being Erd\H{o}s-Rado via an application of the Generalized Morley's Omitting Types Theorem \ref{gmott-thm}.  Note that none of these results were originally stated in the notation of Definition \ref{spr-def} (especially since that notation was originated for this paper); however, we have translated those results into this language to illustrate our notions.

\begin{example}[Linear orders]
In $\cK^{or}$, the canonical bigness notion is just cardinality, so $I \in \cK^{or}_\mu$ iff $|I|\geq \mu$.  The classic Erd\H{o}s-Rado theorem \cite[Theorem 39]{er-erdosrado} states that, for all $n < \omega$ and $\kappa$
$$\beth_{n-1}(\kappa)^+ \xrightarrow{or} (\kappa^+)^n_\kappa$$
Thus, $\cK^{or}$ is a combinatorial Erd\H{o}s-Rado class witnessed by $(\kappa^+, n) \mapsto \beth_{n-1}(\kappa)^+$ and $(\delta, n) \mapsto \beth_{n-1}(\delta)^+$ for limit $\delta$.
\end{example}

Note that the classic results on the Sierpinski coloring \cite{s-relations}\footnote{For a more accessible reference, see \cite[After Question 1.2]{d-ramsey}.} show that \emph{dense} linear orders do {\bf not} form a combinatorial Erd\H{o}s-Rado class.  One could work to develop  an infinite version of Ramsey degree, but we defer that for later.

\begin{example}[$\chi$-disjoint linear orders]\label{chior-spr-ex}
 As discussed in the context of Example \ref{bad-color-ex}, the canonical bigness notion for $\cK^{\chi-or}$ says that $\bar{I}$ is $\mu$-big iff every piece has size at least $\mu$.  Proposition \ref{chior-part-prop} (in which the heavy lifting is done by Shelah \cite[Appendix, Theorem 2.7]{sh:c}, see Remark \ref{ehr-rem}) proves\footnote{By inspection, Proposition \ref{chior-part-prop} gives a stronger result, but we use this weakening for easier comparison with other structural partition relations since
$$\beth_{1}\left(\beth_{n(n+1)}(\kappa)^+\right)^+  \leq \beth_{1}\left(\beth_{n(n+1)+1}(\kappa)\right)^+ = \beth_{n(n+1)+2}(\kappa)^+ $$} , for all $n<\omega$ and $\chi \leq \kappa$,
$$\beth_{n(n+1)+2}(\kappa)^+ \xrightarrow{\chi-or}(\kappa^+)^n_\kappa$$
Thus, $\cK^{\chi-or}$ is a combinatorial Erd\H{o}s-Rado class witnessed by $(\kappa^+, n) \mapsto \beth_{n(n+1)+2}(\kappa)^+$ (and so the threshold for limit $\kappa$ are the same as for $\kappa^+$).
\end{example}

\begin{prop}\label{chior-part-prop}
For all $n<\omega$ and $\chi \leq \kappa$,
$$\beth_{1}\left(\beth_{n(n+1)}(\kappa)^+\right)^+ \xrightarrow{\chi-or}(\kappa^+)^n_\kappa$$
\end{prop}

As alluded to above, Shelah essentially proves this proposition \emph{except} that Shelah's result had a lower cardinal on the left-hand side and assumed each piece to be well-ordered (rather than just linearly ordered); thanks again to the referee for catching this mistake.  We raise the cardinal on the left-hand side to make up for this lack of well-ordering.  The author believes that the theorem could probably be proven with Shelah's original bound (and similarly for the trees in Example \ref{ntr-spr-ex}), but since it doesn't affect the final bound for being an Erd\H{o}s-Rado Class, we don't pursue that here.

{\bf Proof:} Let $\bar{I} = (I_i; <_i)_{i < \chi} \in \cK^{\chi-or}$ be $\lambda$-big for $\lambda = \beth_{1}\left(\beth_{n(n+1)}(\kappa)^+\right)^+$ and fix a coloring
$$c:\left[\bigcup_{i <\chi} I_i\right]^n \to \kappa$$
WLOG, each $I_i$ has size exactly $\lambda$.  In order to invoke Shelah's result, the $I_i$'s must be well-ordered, so fix a well-ordering ordering $<_i^*$ of each $I_i$.

First, we apply the binary Erd\H{o}s-Rado Theorem to the pieces: define colorings 
$$d_i:[I_i]^2 \to 2$$
by considering the agreement of the two orderings
$$d_i(x,y) :=\begin{cases}
    0 & x <_i y \iff x <^*_i y\\
    1 & x<_iy \iff y<^*_i x
\end{cases}$$
By Erd\H{o}s-Rado, there are $I_i^* \subset I_i$ of size $\beth_{n(n+1)}(\kappa)^+$ that are homogeneous.  Now, we can define a well-ordering $<^+$ on all of $\bar{I}^*:=\bigcup_{i<\chi} I^*_i$ to be the concatenation of the $<^*_i$'s.  The key property we have constructed (crucially using the homogeneity) is that, for any tuples $\ba, \bb \in \bar{I}^*$, we have
\begin{eqnarray*}\tp_{\chi-or}\left(\ba;\left( I^*_i, <^*_i\right)_{i<\chi}\right) &=& \tp_{\chi-or}\left(\bb;\left( I^*_i, <^*_i\right)_{i<\chi}\right)\\
&\iff& \\
\tp_{\chi-or}\left(\ba;\left( I^*_i, <_i\right)_{i<\chi}\right) &=& \tp_{\chi-or}\left(\bb;\left( I^*_i, <_i\right)_{i<\chi}\right)
\end{eqnarray*}
%$$x<^+y \iff \begin{cases}
%    x \in I_i \text{ and }y\in I_j\text{ and }i <j\\
%    x,y\in I_i \text{ and }d(x,y)=0 \text{ and }x<_iy\\
%        x,y\in I_i \text{ and }d(x,y)=1 \text{ and }y<_ix
%\end{cases}$$

Now, we restrict the coloring
$$c^*:\left[\bar{I}^*\right]^n\to \kappa$$
and apply Shelah's result to the well-ordered $\left(I_i^*, <^+\right)_{i<\chi}$.  This gives $J_i \subset I^*_i$ of size $\kappa^+$ so $\left(J_i, <^+\right)_{i<\chi}$ is type-homogeneous for $c^*$ and, therefore, for $c$.  By the key property above, the $\kappa^+$-big
$$\left(J_i, <_i\right)_{i<\chi} \subset (I_i, <_i)_{i<\chi} $$
is also type-homogeneous for $c$, as desired.

\comment{Now we can apply apply Shelah's result \cite[Appendix, Theorem 2.7]{sh:c} to the structure $\bar{I}* := (I_i; <^*)_{i<\chi}$ with the coloring $c$ to get, for each $i < \chi$, $J_i \subset I_i$ of size $\beth_{n-1}(\kappa)^+$ such that $(J_i, <^*)_{i<\chi}$ is type-homogeneous for the coloring $c$.  However, the $(J_i, <^*)_{i<\chi}$-type of a tuple is not the same as the $(J_i, <_I)_{i<\chi}$-type of a tuple, and we need the type-homoegeneity to be in the context of the latter.

To solve this, for each $i < \chi$, color each $n$-sized set from $J_i$ with the $(J_i, <^*)$-type of the set listed in increasing $<_I$-order; this is essentially a permutation of $n$ elements, so there are at most $n!$ many colors.  Then apply the normal Erd\H{o}s-Rado theorem to get $K_i \subset J_i$ of size $\kappa^+$ that is homogeneous for this permutation coloring.  This homogeneity means that two $n$-tuples from $\bar{K}:=(K_i, <_I)_{i<\chi}$ having the same type implies that they have the same $(K_i, <^*)$ as well.

Thus, $\bar{K}$ is s $\kappa^+$-large type-homogeneous subset of $\bar{I}$.}\hfill \dag\\

\begin{remark}\label{ehr-rem}
Shelah \cite[Historical Remarks, A.2]{sh:c} cites Erd\H{o}s, Hajnal, and Rado \cite{ehr-partition} (without further specification) for the partition relation cited as \cite[Appendix, Theorem 2.7]{sh:c}.  However, the author and referee were unable to locate the actual result there. The closest is \cite[Corollary 17, p. 176]{ehr-partition}, which states 
$$\begin{pmatrix}
    \kappa^+\\ \kappa^+
\end{pmatrix} \to \begin{pmatrix}
    \kappa & \kappa\\
    \kappa & \kappa
\end{pmatrix}^{1,1}_2$$
Unwinding this notation (see \cite[Section 3.3]{ehr-partition}), this means that if $X_0$ and $X_1$ are of size $\kappa^+$ and 
$$c:X_0 \times X_1 \to 2$$
is a coloring, then there is $X'_\ell \subset X_\ell$ of size $\kappa$ and a color $k<2$ such that $c$ is constant $k$-valued on $X_0' \times X_1'$.  This can be combined with the classic Erd\H{o}s-Rado Theorem $\beth_1(\kappa)^+ \to (\kappa^+)^2_2$ to prove the following partition relation restricted to well-ordered elements of the class
$$\beth_1(\kappa)^+ \xrightarrow{2-or} (\kappa)^2_2$$
Given a coloring $c: [I_0\cup I_1]^2 \to 2$, first use Erd\H{o}s-Rado to find $I_\ell' \subset I_\ell$ of size $\kappa^+$ such that the each restriction to $[I_\ell']^2$ is constant.  Then use the polarized partition relation to find $I_\ell'' \subset I_\ell'$ of size $\kappa$ such that the the coloring restricted to $I_0''\times I_1''$ is constant.  Since the three types of $\S^2_{2-or}$ are represented by $[I_0'']^2$, $I_0''\times I_1''$, and $[I_1'']^2$, the coloring is type homogeneous.

If \cite{ehr-partition} investigated cases of the polarized with more colors or greater arity, this would similarly prove the result.  However, they restrict their attention to the $2\times 2$ case above with two colors.  Note that we have pieced together several results that color each type to get a type-homogeneous coloring; this technique can be used when there are finitely many types of the appropriate arity.
\end{remark}

\begin{example}[$\chi$-colored linear orders] \label{clo-spr-ex}
The canonical bigness notion says that $(I, <, P_\beta)_{\beta<\chi}$ is $\kappa$-big iff $\chi \cdot \kappa \leq otp(I)$.  Then $\cK^{\chi-color}$ is a combinatorial Erd\H{o}s-Rado class witnessed by $F(\kappa, n) = \beth_{n-1}(\kappa^\chi)^+$ by Proposition \ref{clo-part-prop}.
\end{example}

\begin{prop}\label{clo-part-prop}
If $\lambda \xrightarrow{or}(\kappa)^n_{\mu_*}$ where $\mu_* = \mu^{(\chi^n)}$, then $\lambda \xrightarrow{\chi-color} (\kappa)^n_\mu$.
\end{prop}

{\bf Proof:}  Given a coloring $c:[\chi \cdot \lambda]^n \to \mu$, we define an auxiliary coloring $d:[\lambda]^n \to {}^{\left(\chi^n\right)}\mu$ given by $d(\gamma_1, \dots, \gamma_n)$ is the function that maps $(i_1, \dots, i_n) \in \chi^n$ to $c(\chi\cdot\gamma_1+i_1, \dots, \chi \cdot\gamma_n+i_n)$.  There is a $\kappa$-sized homogeneous $X \subset \lambda$ by assumption.  Enumerate an initial segment as $\{\gamma_\alpha \mid \alpha < \kappa\}$ and define the set
$$\chi \otimes X := \{\chi \cdot \gamma_{\left(\chi\cdot \beta+i\right)} + i : \chi \cdot \beta + i < \kappa\text{ and }i < \chi\}$$
We claim that this structure (with the induced $\chi$-or structure) is type homogeneous for $c$.

Suppose that $(\chi\cdot \gamma_{(\chi \cdot \beta_1 + i_1)} + i_1, \dots, \chi\cdot \gamma_{(\chi \cdot \beta_n + i_n)} + i_n)$, $(\chi\cdot \gamma_{(\chi \cdot \beta'_1 + i'_1)} + i'_1, \dots, \chi\cdot \gamma_{(\chi \cdot \beta'_n + i'_n)} + i'_n) \in \chi \otimes X$ have the same $\chi$-or type.  First, this means that $P_\beta$ holds of $\chi\cdot \gamma_{(\chi \cdot \beta_\ell + i_\ell)} + i_\ell$ iff it holds of $\chi\cdot \gamma_{(\chi \cdot \beta'_\ell + i'_\ell)} + i'_\ell$; thus $i_\ell = i'_\ell$ for all $\ell \leq n$.  Second, for ordinals\footnote{Importantly, this fails if we allow $\gamma = \gamma'$.  This is the reason for the complicated subscript in the definition of $\chi\otimes X$.} $\gamma \neq \gamma'$ and $\beta, \beta' < \chi$, we have that
$$\chi \cdot \gamma + \beta < \chi \cdot \gamma' +\beta' \text{ iff }\gamma < \gamma' $$
Thus, we have that for each $\ell, k < n$,
\begin{eqnarray*}
\gamma_{(\chi \cdot \beta_\ell + i_\ell)} < \gamma_{(\chi \cdot \beta_k + i_k)} &\iff& \chi \cdot \gamma_{(\chi \cdot \beta_\ell + i_\ell)} + i_\ell < \chi \cdot \gamma_{(\chi \cdot \beta_k + i_k)} + i_k\\
&\iff& \chi \cdot \gamma_{(\chi \cdot \beta'_\ell + i'_\ell)} + i'_\ell < \chi \cdot \gamma_{(\chi \cdot \beta'_k + i'_k)} + i'_k\\
&\iff& \gamma_{(\chi \cdot \beta'_\ell + i'_\ell)} < \gamma_{(\chi \cdot \beta'_k + i'_k)}
\end{eqnarray*}

Thus, $(\gamma_{(\chi\cdot \beta_1+i_1)}, \dots, \gamma_{(\chi\cdot \beta_1+i_1)}), (\gamma_{(\chi\cdot \beta'_1+i'_1}, \dots, \gamma_{(\chi\cdot \beta'_1+i'_1)}) \in [X]^n$ have the same $or$-type and $d$ of them is the same function.  Thus, we can conclude
\begin{eqnarray*}
c\left(\chi\cdot \gamma_{(\chi \cdot \beta_1 + i_1)} + i_1, \dots, \chi\cdot \gamma_{(\chi \cdot \beta_n + i_n)} + i_n\right) &=& d(\gamma_{(\chi\cdot \beta_1+i_1)}, \dots, \gamma_{(\chi\cdot \beta_1+i_1)})\left(i_1, \dots, i_n\right)\\
&=& d(\gamma_{(\chi\cdot \beta_1+i_1)}, \dots, \gamma_{(\chi\cdot \beta_1+i_1)})\left(i'_1, \dots, i'_n\right)\\
&=& d(\gamma_{(\chi\cdot \beta'_1+i'_1)}, \dots, \gamma_{(\chi\cdot \beta'_1+i'_1)})\left(i'_1, \dots, i'_n\right)\\
&=&c\left(\chi\cdot \gamma_{(\chi \cdot \beta'_1 + i'_1)} + i'_1, \dots, \chi\cdot \gamma_{(\chi \cdot \beta'_n + i'_n)} + i'_n \right)
\end{eqnarray*}

\hfill\dag\\

\begin{example}[Trees of height $n<\omega$] \label{ntr-spr-ex}
 The canonical bigness notion for $\cK^{n-tr}$ is that of splitting: $I \in \cK^{n-tr}_\mu$ iff every node of the tree on level $<n$ has $\geq \mu$-many successors.  As with Example \ref{chior-spr-ex}, Shelah proved a version of the partition relation for trees built on ${}^{\leq n} \lambda$ (so assuming some nice well-ordering) with a rich history of improving the bound\footnote{Shelah \cite{sh10} mentions the result, which appears with proof in \cite[Appendix, Theorem 2.6]{sh:c}.  There, the threshold cardinal on the left-hand side is mentioned to be some polynomial $k(n,m)$ in the height $n$ of the tree and size $m$ of the tuples being colored. Alternate proofs appear in \cite[Theorem 6.7]{kks-treerevisited} and \cite[Theorem 5.1]{gs-unsuperstab}, with the later lowering the bound on $k(n,n)$ from $2^n+n+1$ to $n^2$.  Our method to remove the well-ordering assumption raises the threshold, so we're not concerned about finding the optimal value.}.  By applying a similar argument, Proposition \ref{ntr-part-prop} shows
$$\beth_{\max\{n, m\}^2+2}(\kappa)^+ \xrightarrow{n-tr} (\kappa^+)^m_\kappa$$
Thus, $\cK^{n-tr}$ is a combinatorial Erd\H{o}s-Rado class witnessed by $(\kappa^+, m) \mapsto \beth_{\max\{n,m\}^2+2}(\kappa)^+$.
\end{example}

\begin{prop}\label{ntr-part-prop}
For all $n, m < \omega$ and $\kappa$
$$\beth_{1}\left(\beth_{k(n,m)}(\kappa)^+\right)^+ \xrightarrow{n-tr} (\kappa^+)^m_\kappa$$
where $k(n,m) < \omega$ is a function as in  \cite[Appendix, Theorem 2.6]{sh:c}.
\end{prop}

{\bf Proof:} The proof follows the same strategy as Proposition \ref{chior-part-prop}.  First, we take a $\lambda=\beth_{1}\left(\beth_{k(n,m)}(\kappa)^+\right)^+$-big $n$-tree $T$ and a coloring $c$.  WLOG $T$ is exactly $\lambda$-splitting at each non-terminal node.  Then we can provide a global well-ordering so that $T$ is a copy of ${}^{\leq n}\lambda $.  First, we use binary Erd\H{o}s-Rado to pass to a $\beth_{k(n,m)}(\kappa)^+$-big $T_* \subset T$ where $T_*$ is either the well-ordered or the reverse well-ordered on each level.  Second, we use the Shelah result on the well-ordered $T_*$ to get a $\kappa^+$-big $T_{**}$ that is type-homogeneous for $c$ with the well-ordering.  Using the initial coloring, this means that $T_{**}$ is type-homogeneous for the original ordering as well.\hfill\dag\\

\comment{Then the global order is redone (and the tree shrank) so that the tree is a copy of ${}^{\leq n} \beth_{k(n, m)}\left(\beth_{m-1}(\kappa)^+\right)^+$.  Shelah's result \cite[Appendix, Theorem 2.6]{sh:c} is applied to get $T^* \subset T$ that is $\beth_{m-1}(\kappa)^+$-large and type-homogeneous for $c$ for types with the redone order.

In order to make the types according to the two orders agree, we use the same coloring as in Proposition \ref{chior-part-prop}, but work level-by-level, starting at the level after the root.  That is, using Erd\H{o}s-Rado, we make sure that the orderings on $T^*_1$ agree; then there are $\kappa^+$-many nodes on level 1, but each has $\beth_{m-1}(\kappa)^+$-many successors on level 2.  Then the Erd\H{o}s-Rado theorem is applied in turn on the remaining successor sets at each level, leaving with a $\kappa^+$-large, type homogeneous tree. \hfill \dag\\}

\begin{example}[Trees of height $\omega$] \label{wtr-spr-ex}
We do not know if $\cK^{\omega-tr}$ is a combinatorial Erd\H{o}s-Rado class (although we would expect the bigness notion to be splitting).  However, we are still able to show that is an Erd\H{o}s-Rado class (see Corollary \ref{wtree-er-cor}).
\end{example}

\begin{example}[Well-founded trees]\label{wftr-spr-ex}
We say that a well-founded tree is $\lambda$-big iff it contains a copy of $\ds(\lambda)$, which is the well-founded tree that is made up of all decreasing sequences of ordinals less than $\lambda$ ordered by end extension.  Then \cite[Conclusion 2.4]{gs-wfd-tree} shows that, for every $n < \omega$ and $\kappa$,
$$(\beth_{1, n}(\kappa))\xrightarrow{wf-tr}(\kappa)^n_\kappa$$
where $\beth_{1, n}(\lambda)$ is defined by:
\begin{itemize}
	\item $\beth_{1,0}(\lambda) = \lambda$ and
	\item $\beth_{1, k+1}(\lambda) = \beth_{\beth_{1, k}(\lambda)^+}(\lambda)$
\end{itemize}
\end{example}

Note that the bound here is much larger than the other bounds (which are all below $\beth_\omega(\kappa)$).  For instance,
$$\beth_{1,2}(\kappa) = \beth_{\beth_{\kappa^+}(\kappa)^+}(\kappa)$$
This impacts the witness for being an Erd\H{o}s-Rado class (see Remark \ref{er-bound-rem}), but we do not know if the left-hand side here is a tight bound.  Also, well-founded trees are closely related to \emph{scattered} linear orders (those not containing a copy of $\bQ$; see \cite[Observation 4]{gs-wfd-tree} building on work of Hausdorff), so this result forms a counterpoint to the nonexample coming from Sierpinski colorings.

$\cK^{wf-tr}$ is also important as an Erd\H{o}s-Rado Class that is not axiomatiable in $\bL_{\infty, \omega}$ and does not even form an Abstract Elementary Class because it is not closed under unions of chains.  Worse, the class of well-founded trees could not even be the models of such a class (ignoring the substructure relation).  These models could be used to define well-ordering and this is impossible in Abstract Elementary Classes.

\begin{example}[Convexly-ordered equivalence relations]\label{ceq-spr-ex}
The canonical bigness notion says that $(I, <, E) \in \cK^{ceq}$ is $\mu$-big iff there are at least $\mu$-many equivalence classes, each of which is of size at least $\mu$.  Proposition \ref{ceq-part-prop} below shows that $\cK^{ceq}$ is a combinatorial Erd\H{o}s-Rado class.
\end{example}

\begin{prop}\label{ceq-part-prop}
Given infinite $\kappa$ and $n < \omega$, we have
$$\beth_{n(n+3)}(\kappa)^+ \xrightarrow{ceq} \left(\kappa^+\right)^n_\kappa$$
\end{prop}

{\bf Proof:} Let $(I, E, <) \in \cK^{ceq}_{\beth_{n(n+3)}(\kappa)^+}$ and color it with $c:[I]^n \to \kappa$.  We will use two already established facts:
$${\beth_{n(n+3)}(\kappa)^+} \xrightarrow{{\beth_{n-1}(\kappa)^+}-or} ({\beth_{n-1}(\kappa)^+})^n_\kappa$$
$${\beth_{n-1}(\kappa)^+} \xrightarrow{or}(\kappa^+)^n_\kappa$$
First, find $\{i_\alpha \in I \mid \alpha < {\beth_{n-1}(\kappa)^+}\}$ that are $E$-nonequivalent.  Set $I_1 = \bigcup_{\alpha<{\beth_{n-1}(\kappa)^+}} (i_\alpha/E)$ and note that $(I_1, i_\alpha/E, <)_{\alpha<{\beth_{n-1}(\kappa)^+}}\in \cK^{{\beth_{n-1}(\kappa)^+}-or}_{\beth_{n(n+3)}(\kappa)^+}$.  Then $c$ still colors $[I_1]^n$, so use the result to find $I_2 \subset I_1$ and $c^*:\S_{{\beth_{n-1}(\kappa)^+}-or}\to \kappa$ so that $(I_2, i_\alpha/E \cap I_2, <)_{\alpha < {\beth_{n-1}(\kappa)^+}} \in \cK^{{\beth_{n-1}(\kappa)^+}-or}_{\beth_{n-1}(\kappa)^+}$ is type-homogeneous for $c$ with $c^*$.

Now consider the structure $\left(\{i_\alpha \mid \alpha < {\beth_{n-1}(\kappa)^+}\}, <\right) \in \cK^{or}_{\beth_{n-1}(\kappa)^+}$.  We want to give an auxiliary coloring $d:[{\beth_{n-1}(\kappa)^+}]^n \to {}^A \kappa$, where $A  = \{ s\in {}^n (n+1) \mid \sum_{i<n} s(i) = n\}$.  Then 
$$d\left(\{\alpha_1 < \dots < \alpha_n\}\right)$$ is the function that takes $s \in A$ to $c(\{j_1, \dots, j_n\})$ for $j_1, \dots, j_n \in I_2$ such that, for each $k$, $s(k)$-many of the $j_\ell$'s come from the equivalence class of $i_{\alpha_k}$.  Note that this is a well-defined coloring because $I_2$ was $\cK^{\beth_{n-1}(\kappa)^+-or}$-type-homogeneous for $c$.  Then we can find $X \subset {\beth_{n-1}(\kappa)^+}$ of size $\kappa^+$ and $d^*:A \to \kappa$ such that $X$ is $\cK^{or}$-type-homogeneous for $d$ with color $d^*$.

Set $I_* = \{i\in I_2 \mid i E i_\alpha \text{ for some }\alpha \in X\}$, $E_* = E \rest (I_*^2)$, and $<_* = < \rest (I_*^2)$.
 
{\bf Claim:} $(I_*, E_*, <_*) \in \cK^{ceq}_{\kappa^+}$ is type-homogeneous for $c$.
Since $|X|={\kappa^+}$, $I_*$ has ${\kappa^+}$-many equivalence classes.  For each $\alpha \in X$, $i_\alpha/E_* = i_\alpha/E \cap I_2$ and has size at least ${\beth_{n-1}(\kappa)^+} > {\kappa^+}$.  Thus, $(I_*, E_*, <_*)$ is ${\kappa^+}$-big.

For homogeneity, let $j_1 <_* \dots <_* j_n; j_1' <_* \dots <_* j_n' \in I_*$ have the same $\cK^{ceq}$-type.  Then these tuples are each $<$-increasing, from $I_2$, and each element of each tuple is equivalent to an element of $\{i_\alpha \mid \alpha \in X\}$.  Because they have the same $\cK^{ceq}$-type, there are $\alpha_1 < \dots < \alpha_n; \alpha_1' < \dots < \alpha_n'$ from $X$ that contain these witnesses and \emph{a single map} $s \in A$ that maps $\ell<n$ to 
$$|\{k \mid j_k E_* i_{\alpha_\ell}\}| = |\{k \mid j'_k E_* i_{\alpha'_\ell}\}|$$
By the homogeneity of $X$, we have that $d^* = d\left(\{\alpha_1, \dots, \alpha_n\}\right) = d\left(\{\alpha'_1, \dots, \alpha'_n\}\right)$.  Thus,
\begin{eqnarray*}
c\left(\{j_1, \dots, j_n\}\right) &=& d\left(\{\alpha_1, \dots, \alpha_n\}\right)(s)\\
&=& d^*(s)\\
&=&d\left(\{\alpha'_1, \dots, \alpha'_n\}\right)(s)\\
&=& c\left(\{j'_1, \dots, j'_n\}\right)
\end{eqnarray*}
\hfill\dag\\

\begin{example}[$n$-multi-orders]\label{nmlo-spr-ex}
We can use $\cK^{n-mlo}$ to point out that the choice of bigness notion is very important.  If we say $(I, <_1, \dots, <_n)$ is $\mu$-big when $|I|\geq \mu$, then $\cK^{n-mlo}$ is a combinatorial Erd\H{o}s-Rado class simply because $\cK^{or}$ is.  However, this gives us no new information.  A good bigness notion for this class should say something about the independence of the different linear orders.
\end{example}

\begin{example}[Ordered graphs]\label{og-spr-ex}
Ordered graphs start to indicate that set theory begins to enter the picture.  Hajnal and Komj\'{a}th \cite[Theorem 12]{hk-graphs} (with correction at \cite[Theorem 12]{hk-graphs-errata}) show that it is consistent that there is a graph that never appears as a monochromatic subgraph.  In particular, they start with a model of $GCH$, add a single Cohen real, and construct an uncountable bipartite graph $G$ such that every graph $H$ has a coloring of pairs such that there is no type-homogeneous copy of $G$ in $H$.  On the other hand, the next example (which subsumes this one by considering $\cK^{(2, 2)-hg}$) shows that we can consistently get a combinatorial Erd\H{o}s-Rado result.

\end{example}

\begin{example}[Colored hypergraphs]\label{chg-spr-ex}
Shelah \cite[Conclusion 4.2]{sh289} proved that it is consistent that an Erd\H{o}s-Rado Theorem holds for the classes $\cK^{(k, \sigma)-hg}$ with $k < \omega$.  Specifically, he shows that, after an iterated forcing construction, for every well-ordered $N \in \cK^{(k, \sigma)-hg}$, $m < \omega$, and $\kappa$, there is a $M \in \cK^{(k, \sigma)-hg}$ with $\|M\| < \beth_\omega(\|N\|+\sigma+\kappa)$ such that any coloring of $[M]^m$ with $\kappa$-many colors contains a type-homogeneous substructure isomorphic to $N$.  For the right bigness notion and by using the techniques of Propositions \ref{chior-part-prop} and \ref{ntr-part-prop} to remove the well-ordered assumption, this means that
$$\beth_{\omega}(\kappa) \xrightarrow{(k, \sigma)-hg} (\kappa)^n_\kappa$$

\end{example}

\section{Erd\H{o}s-Rado Classes and the Generalized Morley's Omitting Types Theorem}\label{er-sec}

Erd\H{o}s-Rado classes (Definition \ref{er-def}) are those that allow one to build generalized indiscernibles  in nonelementary classes, especially those definable in terms of type omission.  Since these classes are often axiomatized in stronger logics, one could formulate the modeling property of Ramsey classes in terms of these stronger logics (in fact, Shelah \cite{she59} does this, and we compare the notions in Remark \ref{er-ram-remark}).  However, this is not how order indiscernibles are typically built in Abstract Elementary Classes.  Instead, we continue to work with indiscernability in a first-order (and even quantifier-free) context, but strengthen the modeling property so that type omission is preserved.

Note that there are two variants of being an Erd\H{o}s-Rado class here, and a few more in Definition \ref{er-variants-def}.  The cofinal variant is the most common, and gives the sharpest applications.

\begin{defin}\label{er-def}
Let $\cK$ be an ordered abstract class.
\begin{enumerate}
	\item $\cK$ is a \emph{$(\mu, \chi, \bg)$-Erd\H{o}s-Rado class} iff for every language $\tau$ of size $\leq \mu$, every $I \in \cK_\chi^\bg$, every $\tau$-structure $M$, and every injection $f:I \to M$, there is a blueprint $\Phi \in \Upsilon^\cK[\tau]$ such that
	\begin{enumerate}
		\item $\tau(\Phi) = \tau$; and
		\item \label{erc-item} for each $p \in \S^{inc}_\cK$, there are $i_1<\dots< i_n \in I$ realizing $p$ such that
		$$\tp_{\tau}\left(f(i_1), \dots, f(i_n); M\right) = \Phi(p)$$
	\end{enumerate}
	\item $\cK$ is a \emph{cofinally $(\mu, \chi, \bg)$-Erd\H{o}s-Rado class} iff for every language $\tau$ of size $\leq \mu$, if we have, for each cardinal $\alpha < \chi$, a $\tau$-structure $M_\alpha$, an $\alpha$-\big $I_\alpha \in \cK$, and an injection $f_\alpha:I_\alpha \to M_\alpha$, then there is a blueprint $\Phi \in \Upsilon^\cK[\tau]$ such that
	\begin{enumerate}
		\item $\tau(\Phi) = \tau$; and
		\item \label{erc2-item} for each $p \in \S^{inc}_\cK$, there are cofinally many $\alpha < \chi$ such that there are $i_1<\dots < i_n \in I_\alpha$ realizing $p$ such that
		$$\tp_{\tau}\left(f_\alpha(i_1), \dots, f_\alpha(i_n); M_\alpha\right) = \Phi(p)$$
	\end{enumerate}
\end{enumerate}
In either case, writing `$\cK$ is a [cofinally] \bg-Erd\H{o}s-Rado class' means that `there is a function $f:\Card \to \Card$ such that $\cK$ is [cofinally] $(\mu, f(\mu), \bg)$-Erd\H{o}s-Rado for every $\mu$.'  If \big is the standard bigness notion for $\cK$, then we omit it.
\end{defin}

Note that (2) implies (1) by taking the constant sequences $I_\alpha = I$ and $M_\alpha = M$.  Also note that the definition of cofinal Erd\H{o}s-Rado classes is only new when the $\chi$ is limit; otherwise, the definition is equivalent to a normal Erd\H{o}s-Rado class at its predecessor.

We refer to Definition \ref{er-def}.(\ref{erc-item}) or Definition \ref{er-def}.(\ref{erc2-item}) as the \emph{Erd\H{o}s-Rado condition}.  See Remark \ref{er-ram-remark} for a comparison with Ramsey conditions.

Before our main theorem, we need one more auxilliary definition.  This definition will transform the witnessing function for combinatorial Erd\H{o}s-Rado classes into witnessing functions for cofinal Erd\H{o}s-Rado classes.

\begin{defin}\label{daleth-wit-def}
    Fix a function $F:\Card\times\omega\to \Card$.  The \emph{relative $\daleth$ sequence for $F$} is\footnote{$\daleth$ (`daleth') is the fourth letter of the Hebrew alphabet and is chosen to follow in the line of the well-known $\aleph$ and $\beth$ functions and the less well-known $\gimel$ function (see, e.g., \cite[Equation (5.19), p. 56]{jech}} the ordinal-indexed sequence $\daleth_\alpha=\daleth^F_\alpha$ given by 
    \begin{eqnarray*}
    \daleth_0 &=& \aleph_0\\
    \daleth_{\alpha+1} &=& \sup_{n<\omega}F(\daleth_\alpha, n)\\
    \daleth_\delta &=& \sup_{\alpha<\delta} \daleth_\alpha
\end{eqnarray*}
With $F:\Card\times \omega\to \Card$ and an abstract class $\cK$, we define two more functions
\begin{eqnarray*}
    g^{\cK}(\mu) &=& \sup_{n<\omega} \left(2^{\mu\cdot|\S^n_\cK|}\right)\\
    f^{\cK,F}(\mu)&=& \daleth^F_{g(\mu)^+}
\end{eqnarray*}
\end{defin}
Although a formally new cardinal, in practice, this tends to give the same bounds as the normal Morley's Omitting Types Theorem (see Remark \ref{er-bound-rem}).

The following theorem is the main source of Erd\H{o}s-Rado classes.

\begin{theorem}[Generalized Morley's Omitting Types Theorem]\label{gmott-thm}
Let $\cK$ be combinatorially Erd\H{o}s-Rado witnessed by $F$. Then $\K$ is cofinally Erd\H{o}s-Rado witnessed by $f^{\cK,F}$.
\end{theorem}

The statement of Theorem \ref{gmott-thm} makes no mention of types or their omission in the statement, although this notion appears in the name.  This connection is made in Theorem \ref{gmott-ap-thm}, which shows that if a blueprint is generated from (Skolemized) models that all omit some type, then any model built from the blueprint also omits the type.  Morley is classically credited with the original result and gave the basic argument, but an important technical improvement is due to Chang, see \cite[Proof of (E), p. 49]{c-cpt}.\\

{\bf Proof:} Fix $\daleth_\alpha=\daleth^F_\alpha$, $g=g^\cK$, and $f=f^{\cK,F}$ from Definition \ref{daleth-wit-def}.

Suppose that we are given $f_\alpha:I_\alpha\to M_\alpha$ for cardinals $\alpha < f(\mu)$. First, we thin out the sequence in two ways. Although $\cK$ has a unique $0$-type (recall Remark \ref{0-type-remark}), the models $M_\alpha \in \cK^\tau$ might realize distinct $0$-types.  However, the number of possible $0$-types is $\leq 2^{|\tau|} < g(\mu)^+$.  Thus, by passing to cofinal sequence (and using the monotonicity of bigness)\footnote{A more detailed version of this technique is given at the end of the induction step of the construction.} and restricting the models, we instead work with functions $f_\alpha:I_\alpha \to M_\alpha$  for each $\alpha < g(\mu)^+$ such that each $I_\alpha$ is $\daleth_\alpha$-big and all $M_\alpha$ satisfy the same $0$-type $p^*\in \S_\tau^0$.

We will build, for $n<\omega$ and $\alpha < g(\mu)^+$,
\begin{itemize}
    \item $\Phi_n:\S^{inc, n}_\cK\to S^n_\tau$;
    \item $\beta_n(\alpha)<g(\mu)^+$;
    \item $\gamma_{n+1}(\alpha) < g(\mu)^+$;
    \item $I^n_\alpha \in\cK$ that is $\daleth_\alpha$-big;
    \item $h^{n+1}_\alpha:I^{n+1}_\alpha \to I^n_{\gamma_{n+1}(\alpha)}$ that is a $\cK$-morphism; and 
    \item $f^n_\alpha:I^n_\alpha \to M_{\beta_n(\alpha)}$
\end{itemize}

such that 
\begin{enumerate}
    \item $\beta_0(\alpha) = \alpha$; $I^0_\alpha = I_\alpha$; and $f^0_\alpha = f_\alpha$;
    \item\label{phi-cond} for each $\alpha < g(\mu)^+$ and $i_1<\dots<i_n \in I^n_\alpha$, we have that
    $$\Phi_n\left(\tp_\K(i_1, \dots, i_n; I_\alpha^n)\right) = \tp_\tau\left(f_\alpha^n(i_1), \dots, f_\alpha^n(i_n); M_{\beta_n(\alpha)}\right)$$
    \item \label{coh-cond} the $\Phi_n$ are coherent in the following sense: if $p \in \S_\K^{inc,n}$  and $s \subset n$, then
    $$\Phi_n\left(p\right)^s = \Phi_{|s|}\left(p^s\right)$$
    (see Definition \ref{stone-def}.(\ref{type-rest-item}) for this notation);\footnote{Note that this condition follows from the previous one.} and
    \item\label{gmott-cond} given $\alpha < g(\mu)^+$ and $n<\omega$, we have that $\alpha \leq \beta_n(\alpha)$; $\alpha \leq \gamma_n(\alpha)$; and $\beta_{n+1}(\alpha) = \beta_n\left(\gamma_{n+1}(\alpha)\right)$, and the following commutes
    \[\xymatrix{I_\alpha^{n+1} \ar[rr]^{f_\alpha^{n+1}} \ar[dr]_{h_\alpha^{n+1}}& & M_{\beta_{n+1}(\alpha)}\\
    & I_{\gamma_{n+1}(\alpha)}^n \ar[ur]_{f^n_{\gamma_{n+1}(\alpha)}} &}\]

\end{enumerate}
{\bf This is enough:}  Set $\Phi := \bigcup_{n<\omega}\Phi_n$.  Then this is a function with domain $\S^{inc}_\cK$ and range $\S_\tau$.  Moreover, the coherence condition (\ref{coh-cond}) implies that it is proper for $\K$.  Now we wish to show that it has the type reflection required by the Erd\H{o}s-Rado condition, see Definition \ref{er-def}.(\ref{erc2-item}).

Let $p \in \S_{\cK}^{inc,n}$ and $\alpha_0 < f(\mu) = \daleth_{g(\mu)^+}$.  Then there is $\alpha_1<g(\mu)^+$ so $\alpha_0<\daleth_{\alpha_1}$. $I^n_{\alpha_1}$ is $\aleph_0$-\bg, so there is $i_1< \dots< i_n \in I^n_{\alpha_1}$ realizing $p$.  Then, by (\ref{phi-cond}) of the construction
\begin{eqnarray*}
\Phi(p) = \Phi_n\left(\tp_\K(i_1, \dots, i_n; I_{\alpha_1}^n)\right) = \tp_\tau\left(f_{\alpha_1}^n(i_1), \dots, f_{\alpha_1}^n(i_n); M_{\beta_n(\alpha_1)}\right)
\end{eqnarray*}
If $n=0$, then $I^0_{\alpha_1} = I_{\alpha_1}$ and we are done.  If $n>0$, we compose the $h$-embeddings to define
$$h^*:=h^1_{\gamma_1^{-1}(\beta_n(\alpha_1))}\circ\dots\circ h^{n-1}_{\gamma_{n+1}(\alpha_1)}\circ h^n_{\alpha_1}: I^{n+1}_\alpha \to I_{\beta_n(\alpha_1)}$$
This satisfies $f^n_{\alpha_1} = f_{\beta_n(\alpha_1)} \circ h^*$.  Thus, $h^*(i_1), \dots, h^*(i_n) \in I_{\beta_n(\alpha_1)}$ realize $p$ and 
$$\Phi(p) = \tp_\tau\left(f_{\beta_n(\alpha_1)}\left(h^*(i_1)\right), \dots, f_{\beta_n(\alpha_1)}\left(h^*(i_n)\right); M_{\beta_n(\alpha_1)}\right)$$
Since $\daleth_{\beta_n(\alpha_1)}> \alpha_0$, this completes the proof.\\

{\bf Construction:} We work by induction on $n$.\\  

For $n=0$, we must deal with the unique $0$-type $p_0 \in \S^{inc, 0}_\cK$.  Recalling the first step to ensure each $M_\alpha$ has the same $0$-type $p^*$, we can set $\Phi_0:= \{(p_0, p^*)\}$ and use what we are given: $\beta_0(\alpha)=\alpha$; $I_\alpha^0 =I_\alpha$; and $f_\alpha^0=f_\alpha$.\\

For $n+1$, suppose we have completed the construction up to stage $n$.  Fix some $\alpha < g(\mu)^+$.  Then we have $F(\daleth_\alpha, n+1) \leq \daleth_{\alpha+1}$ by definition.  Consider the coloring
$$c^{n+1}_\alpha:[I^n_{\alpha+1}]^{n+1}\to \S^{n+1}_\tau$$
given by, for $i_1<\dots<i_{n+1} \in I^n_{\alpha+1}$,
$$c^{n+1}_\alpha \left(\{i_1,  \dots , i_{n+1}\}\right) = \tp_\tau\left(f^n_{\alpha+1}(i_1), \dots, f^n_{\alpha+1}(i_{n+1}); M_{\beta_n\left(\alpha+1\right)}\right)$$
By monotonicity of bigness, we have   $$\left(\daleth_{\alpha+1}\right) \xrightarrow{\cK} \left(\daleth_\alpha\right)^{n+1}_{2^\mu}$$
Thus there is $\daleth_\alpha$-big $\bar{I}^{n+1}_\alpha \in \cK$; $\bar{h}^{n+1}_\alpha: \bar{I}^{n+1}_\alpha \to I^n_{\alpha+1}$; and $c^{*, n+1}_\alpha:\S^{inc, n+1}_\cK\to \S^{n+1}_\tau$ witnessing the type-homogeneity, that is, such that for all $i_1<\dots<i_{n+1} \in \bar{I}^{n+1}_\alpha$, we have
$$\tp_\tau\left(f_{\alpha+1}^{n}\circ \bar{h}^{n+1}_\alpha(i_1), \dots, f_{\alpha+1}^{n}\circ \bar{h}^{n+1}_\alpha(i_{n+1}); M_{\beta_n\left(\alpha+1\right)}\right) = c^{*, n+1}_\alpha\left(\tp_\cK(i_1, \dots, i_{n+1}; \bar{I}_\alpha^{n+1})\right)$$
For each $\alpha < g(\mu)^+$, we have built a function $c_{\alpha}^{*, n+1}:\S^{inc, n+1}_\cK\to S^{n+1}_\tau$. Since $\cf(g(\mu)^+)=g(\mu)^+$ is greater than the number of these functions, there is $X\subset g(\mu)^+$ of size $g(\mu)^+$ and $c^{*, n+1}:\S^{inc, n+1}_\cK\to \S^{n+1}_\tau$ such that, for all $\alpha \in X$, $c^{*, n+1}_\alpha=c^{*, n+1}$.  Set $\pi:X\cong g(\mu)^+$ to be the collapse of $X$.  Then $\alpha \leq \pi^{-1}(\alpha)$ for all $\alpha \in g(\mu)^+$.  

Set
\begin{enumerate}
    \item $\Phi_{n+1} = c^{*, n+1}$;
    \item $I_\alpha^{n+1} = \bar{I}^{n+1}_{\pi^{-1}(\alpha)}$;
    \item $\gamma_{n+1}(\alpha) = \pi^{-1}(\alpha)+1$;
    \item $\beta_{n+1}(\alpha) =\beta_n\left(\pi^{-1}(\alpha)+1\right)$;
    \item $h_\alpha^{n+1}=\bar{h}^{n+1}_{\pi^{-1}(\alpha)}$;
    \item $f_\alpha^{n+1} = f^n_{\pi^{-1}(\alpha)+)} \circ h_\alpha^{n+1}$
\end{enumerate}
These satisfy the pieces of the construction: by induction, each of the $c^{*,n+1}_\alpha$'s extend $c^{*, n}$ in the sense that the restriction to $n$-types is determined by $c^{*,n}$.  This gives the coherence.  The other properties are routine to verify.\hfill \dag\\

\begin{corollary}\label{cer-er-cor}
Each of the examples of combinatorial Erd\H{o}s-Rado classes in Section \ref{spr-example-ssec} are cofinal Erd\H{o}s-Rado classes.
\end{corollary}

\begin{remark}\label{er-bound-rem}
Whenever $F(\mu, n) \leq \beth_\omega(\mu)$ and $\mu \geq |\S^n_\cK|$ for all $n<\omega$, then this gives the bound $f(\mu) = \beth_{\left(2^\mu\right)^+}$ that often appears in the theory of nonelementary classes.  In the case of well-founded trees, we get the bound $\beth_{1, \left(2^\mu\right)^+}$.  These bounds can be improved by phrasing in terms of the undefinability of well-ordering of certain PC classes (this is done for specific cases in \cite{sh:c, gs-unsuperstab}).
\end{remark}

The following extends the normal notion of PC classes to include classes with a strong substructure relation.  Note that Chang's Presentation Theorem \cite{c-cpt} implies any $\bL_{\infty, \omega}$-axiomatizable class with `elementary according to a fragment' as the strong substructure is what we will call a \emph{PC pair}, and Shelah's Presentation Theorem \cite{sh88} extends this to Abstract Elementary Classes.

\begin{defin} Let $(\bK, \prec_\bK)$ be an abstract class with $\tau = \tau(\bK)$.
\begin{enumerate}
    \item Let $\bK$ be an abstract class with $\tau = \tau(\bK)$.  $\bK$ is a \emph{PC class} iff there is a language $\tau_1 \supset \tau$, a (first-order) $\tau_1$-theory $T_1$, and a collection $\Gamma$ of $\tau_1$-types such that, for any $\tau$-structure $M$, $M \in \bK$ iff there is an expansion $M_1$ of $M$ to $\tau_1$ that models $T_1$ and omits all types in $\Gamma$.
    \item Let $\bK$ be a class of $\tau$-structures and $\prec_\bK$ be a partial order on $\bK$.  $(\bK, \prec_\bK)$ is a \emph{PC pair} iff there is a language $\tau_1 \supset \tau$, a $\tau_1$-theory $T_1$, and a collection $\Gamma$ of $\tau_1$-types such that
    \begin{itemize}
        \item for any $\tau$-structure $M$, $M \in \bK$ iff there is an expansion $M_1$ of $M$ to $\tau_1$ that models $T_1$ and omits all types in $\Gamma$; and
        \item for any $M, N \in \bK$, $M \prec_\bK N$ iff there are expansions $M_1$ of $M$ and $N_1$ of $N$  to $\tau_1$ that models $T_1$ and omits all types in $\Gamma$ such that $M_1 \subset_{\tau_1} N_1$
    \end{itemize}
\end{enumerate}
\end{defin}

\begin{theorem} \label{gmott-ap-thm}
Let $\cK$ be a cofinal Erd\H{o}s-Rado class witnessed by $f$ and let $(\bK, \prec_\bK)$ be a PC pair with $\tau=\tau(\bK)$ and $\tau_1$ the witnessing language.  Suppose that, for every $\alpha<f(|\tau_1|)$, there is $M_\alpha \in \bK$; $\alpha$-big $I_\alpha \in \cK$; and $f_\alpha:I_\alpha \to M_\alpha$.  Then, there is $\Phi\in \Upsilon^{\cK}_{|\tau_1|}[\bK]$ such that, for every $p \in S_{\cK}^{inc}$, there are cofinally many $\alpha < f(|\tau_1|)$ such that there are $i_1< \dots< i_n \in I_\alpha$ realizing $p$ such that
$$\tp_\tau\left(f_\alpha(i_1), \dots, f_\alpha(i_n); M_{\alpha}\right) = \Phi(p)$$
\end{theorem}

A version of Theorem \ref{gmott-ap-thm} also holds for Erd\H{o}s-Rado classes (without the cofinal adjective) when there is a single embedding from a $f(|\tau_1|)$-big member of $\cK$ into $M$.

{\bf Proof:} Let $T_1$ and $\Gamma$ in the language $\tau_1$ witness that $\bK$ is a PC pair.  By a further Skolem expansion, we can assume that $T_1$ is universal and the types of $\Gamma$ are quantifier-free. Let $f_\alpha:I_\alpha \to M_\alpha$ for $\alpha < f(|\tau_1|)$ as in the hypothesis.  Since $\cK$ is $\left(|\tau_1|, f(|\tau_1|)\right)$-cofinally Erd\H{o}s-Rado, we can find a blueprint $\Phi \in \Upsilon^{\cK}[\tau_1]$ satisfying the Erd\H{o}s-Rado condition, Defintion \ref{er-def}.(\ref{erc2-item}).

First, we wish to show that $\Phi$ is proper for $(\cK, \bK)$. For membership in $\bK$, let $I \in \cK$.  First, suppose that $`\forall \bx \phi(\bx)\text{'} \in T_1$ with $\phi(\bx)$ quantifier-free.  If $EM(I, \Phi) \vDash \neg \forall \bx \phi(\bx)$, then there is $\ba \in EM(I, \Phi)$ witnessing this. Since $EM(I, \Phi)$ is generated by $\tau_1$-terms, there are $\tau_1$-terms $\sigma_1,\dots, \sigma_n$ and $i_1, \dots, i_k \in I$ such that $\ba = \sigma_1^{EM(I, \Phi)}(\bi), \dots, \sigma_n^{EM(I, \Phi)}(\bi)$; without loss, these satisfy $i_1 < \dots < i_k$.

Set $p = \tp_\cK(\bi; I)$.  By the Erd\H{o}s-Rado condition, there is some $\alpha < f(\mu)$ and $j_1<\dots<j_k$ such that 
$$\tp_{\tau_1}\left(i_1, \dots, i_k; EM(I, \Phi)\right) = \Phi(p) = \tp_{\tau_1}\left(f_\alpha(j_1), \dots, f_\alpha(j_k); M_\alpha\right)$$
In particular, 
$$M_\alpha \vDash \neg\phi\left(\sigma_1(\bj), \dots, \sigma_n(\bj)\right)$$
But this contradicts that $M_\alpha$ models $T_1$.

The same argument shows that any quantifier-free type which is realized in $EM(I, \Phi)$ is realized in cofinally many $M_\alpha$.  Since all types in $\Gamma$ are quantifier-free, $EM(I, \Phi)$ must omit all of them.  So $EM_\tau(I, \Phi) \in \bK$

For substructure, this follows from the definition for PC pair and the fact that $\Phi$ is proper for $(\cK, \bK^{\tau_1})$. \hfill\dag\\

\begin{remark}
Since a $\Phi$ proper for $(\cK, \bK)$ is a map $\S^{inc}_\cK \to \S_\bK$, the blueprint $\Phi$ also determines Galois types in $\bK$ in the following sense: if $\sigma_1, \dots, \sigma_k$ are $\tau(\Phi)$-terms; $I, J \in \cK$; and $i_1, \dots, i_n \in I$ and $j_1, \dots, j_n\in J$ are tuples such that
$$\tp_{\cK}\left(i_1, \dots, i_n; I\right) = \tp_\cK\left(j_1, \dots, j_n; J\right)$$
then
\begin{eqnarray*}\tp_\bK\left(\sigma_1(i_1, \dots, i_n), \dots, \sigma_k(i_1, \dots, i_n);EM_{\tau(\bK)}(I, \Phi)\right) 
=\tp_\bK\left(\sigma_1(j_1, \dots, j_n), \dots, \sigma_k(j_1, \dots, j_n);EM_{\tau(\bK)}(J, \Phi)\right)
\end{eqnarray*}
This could also be proved by applying the functor $EM_\tau(\cdot, \Phi)$ to the diagram witnessing type equality in $\cK$ to get a diagram proving the type equality in $\bK$.

In particular, $I \subset EM_\tau(I, \Phi)$ is a collection of $\cK$-indiscernibles.

\end{remark}

\begin{remark}\label{er-ram-remark}
We want to highlight the differences between the Erd\H{o}s-Rado condition (Definition \ref{er-def}.(\ref{erc-item})) to the relevant condition in uses of Ramsey classes, such as \cite[Definition 1.15]{she59} or \cite[Definition 2.12]{ghs-collapse}.  We rephrase the Ramsey modeling condition and the Erd\H{o}s-Rado condition to highlight the comparison:
\begin{enumerate}
    \item[{\bf Ramsey:}] for each $p \in \S^{inc}_\cK$ and \underline{for each $\phi(\bx) \in \Phi(p)$}, there is $i_1 < \dots < i_n \in I$ realizing $p$ so
    $$M \vDash \phi\left(f(i_1), \dots, f(i_n); M \right)$$
    \item[{\bf Erd\H{o}s-Rado:}] for each $p \in \S^{inc}_\cK$, there is $i_1 < \dots < i_n \in I$ realizing $p$ so \underline{for each $\phi(\bx) \in \Phi(p)$}
    $$M \vDash \phi\left(f(i_1), \dots, f(i_n); M \right)$$
\end{enumerate}
The witnesses for Ramsey condition depend on the formula under consideration, but the witness for the Erd\H{o}s-Rado condition is uniform for all formulas.

This makes Ramsey classes ill-equipped to handle type omission and nonelementary classes.  This is because, after Skolemization, the generating sequence might not agree on where terms omit the types, so the blueprint is not guaranteed to omit types.  Shelah \cite[Definition 1.15]{she59} addresses this by introducing $\cL$-nice Ramsey classes (for a logic fragment $\cL$) that considers formulas in $\cL$.  However, it is unclear how to get a $\cL$-nice Ramsey class outside of Erd\H{o}s-Rado classes.  He also considers the notion of a strongly Ramsey class, which is similar to our notion.
\end{remark}
\section{Further results}\label{ext-sec}

\subsection{Reversing Generalized Morley's Omitting Types Theorem}

We would like to have a converse to the Generalized Morley's Omitting Types Theorem \ref{gmott-thm} that says that all Erd\H{o}s-Rado classes come from a combinatorial result.  However, this seems unlikely to be true (and we discuss candidates for this in Section \ref{er-notcer-ssec}).  The issue is that the definition of a (cofinally) Erd\H{o}s-Rado class is not as tied to the relevant bigness notion, but the definition of a combinatorial Erd\H{o}s-Rado class is.  In particular, the definition leaves open the possibility that there is only a single witness to the Erd\H{o}s-Rado condition, while combinatorial Erd\H{o}s-Rado classes require a big set of witnesses to the type-homogeneity.  If we strengthen this requirement, then we get a converse.

\begin{defin}\label{er-variants-def}
We say that $\cK$ is \emph{strongly $(\mu, \chi, \bg)$-Erd\H{o}s-Rado} iff for all $I \in \cK_\chi^\bg$ and every injection $f:I \to M$ with $|\tau(M)|\leq \mu$, there is a blueprint $\Phi \in \Upsilon^\cK[\tau]$ with $\tau(\Phi) = \tau(M)$ such that for all $\alpha < \chi$ and $n < \omega$, there is an $\alpha$-big $I^n_\alpha \leq_\cK I$ such that, for every $i_1< \dots< i_n \in I^n_\alpha$, we have
$$\tp_{\tau(M)}\left(f(i_1), \dots, f(i_n); M\right) = \Phi \left(\tp_\K(i_1, \dots, i_n; I)\right)$$
We define the cofinal variant and what it means to omit the $(\mu, \chi, \bg)$-prefix as in Definition \ref{er-def}.
\end{defin}

\begin{theorem}\label{gmott-con-thm}
Let $\cK$ be an ordered abstract class.
\begin{enumerate}
    \item If $\cK$ is combinatorially Erd\H{o}s-Rado witnessed by $F$, then $\cK$ is a strongly, cofinally Erd\H{o}s-Rado witnessed by the function in Theorem \ref{gmott-thm}.
    \item If $\cK$ is strongly, cofinally $(\mu, \chi, \bg)$-Erd\H{o}s-Rado, then $\cK$ is strongly $(\mu, \chi, \bg)$-Erd\H{o}s-Rado.
    \item If $\cK$ is strongly $(\mu, \chi, \bg)$-Erd\H{o}s-Rado, then, for each $n < \omega$ and $\lambda < \chi$,
    $$(\chi) \xrightarrow[\big]{\cK} (\lambda)^n_\mu$$
    \item If $\cK$ is strongly Erd\H{o}s-Rado witnessed by $f$, then $\cK$ is combinatorially Erd\H{o}s-Rado witnessed by $F(\mu^+, n) = f(\mu)$.
\end{enumerate}
\end{theorem}

{\bf Proof:} The proof of the Generalized Morley's Omitting Types Theorem \ref{gmott-thm} proves (1): the $I^n_\alpha$ built in that proof are exactly the ones needed to witness `strong.'  The proof of (2) is straightforward.  We prove (3), which is enough to prove (4). The idea is that a potential coloring is turned into a structure, and the derived blueprint is used to figure out the colors for the large set.

Let $\lambda < \chi$ and $c:[I]^n \to \mu$ be a coloring of $I \in \cK^{\bg}_\chi$.  We build this into a two-sorted structure
$$M = \langle I, \mu; c, \alpha\rangle_{\alpha<\mu}$$
We have an embedding $f:I \to M$ given by the identity.  Then $|\tau(M)| =\mu$, so the strong Erd\H{o}s-Rado property gives us a blueprint $\Phi:\S^{inc}_\cK \to \S_{\tau(M)}$ as in Definition \ref{er-variants-def} with witnessing sets $I^n_\alpha$.\\

{\bf Claim 1:} For every $p \in \S^{inc, n}_\cK$, there is a unique $\alpha_p < \mu$ such that $``c(x_1, \dots, x_n) = \alpha_p\text{"} \in \Phi(p)$.

Take $I^n_\omega \leq_\cK I$ witnessing the strong Erd\H{o}s-Rado property and find $i_1 < \dots < i_n \in I^n_\omega$ realizing $p$; such a tuple exists by the definition of a bigness notion.  Then $\Phi(p) = \tp_{\tau(M)}(i_1, \dots, i_n; M)$.  This has a color, so $\alpha_p = c(i_1, \dots, i_n)$.\hfill$\dag_{\text{Claim 1}}$\\

Set $c^*:\S^{inc, n}_\cK\to \mu$ to be the function that takes $p$ to $\alpha_p$.\\

{\bf Claim 2:} $I^n_\lambda$ is type-homogeneous for $c$ as witnessed by $c^*$.

Straightforward.\hfill $\dag_{\text{Claim 2}}$\\

Since $I^n_\lambda$ is $\lambda$-\bg, this proves the theorem.\hfill $\dag_{\text{Theorem \ref{gmott-con-thm}}}$\\

Note that this is not an exact converse because there is some slippage in the witnessing functions.  However, this doesn't affect the bounds on the Erd\H{o}s-Rado class.

\subsection{A category theoretic interpretation of blueprints}\label{ct-ind-ssec}

This section gives a category theoretic perspective on the results we've proven and indiscernibles in general.  It requires more category theoretic background  than the rest of the paper (such as \cite{mp-accessible} or \cite{ar-presentable}), but can be skipped. This theme will be explored further in \cite{b-cat-ind}.

Makkai and Par\'{e} give the following statement credited to Morley.  For a logic $\cL$, an `$\cL$-elementary category' is (a category equivalent to) one where the objects are models of some fixed $\cL$-theory $T$ and arrows are $\tau$-homomorphisms between models that are elementary for some fragment of $\cL$ containing $T$.

\begin{fact}[{\cite[Theorem 3.4.1]{mp-accessible}}]\label{mor-mp-fact}
$\cK^{or}$ is a ``minimal'' large, $\bL_{\infty, \omega}$-elementary category.  This means that if $\bK$ is a large, $\bL_{\infty, \omega}$-elementary category, then there is a faithful functor $\Phi:\cK^{or} \to \bK$ that preserves directed colimits.
\end{fact}

This is not phrased as Morley (likely) ever wrote it, but this is the classic proof of Morley's Omitting Types Theorem.  The functor $\Phi$ comes from the blueprint that takes $I \in \cK^{or}$ to $EM_\tau(I, \Phi) \in \bK$.  With generalized indiscernibles in hand, we have a generalization.

\begin{theorem}
Erd\H{o}s-Rado classes are below every large, $\bL_{\infty, \omega}$-elementary category (in the sense of Fact \ref{mor-mp-fact}).  In particular, every $\bL_{\infty, \omega}$-axiomatiable Erd\H{o}s-Rado class is minimal amongst the large, $\bL_{\infty, \omega}$-elementary categories.
\end{theorem}

We include a proof to make the translation more clear (and in part because Makkai and Par\'{e} do not give a proof). Recall that $\cK^{wf-tr}$ is an Erd\H{o}s-Rado class that is not $\bL_{\infty, \omega}$-axiomatizable.\\

{\bf Proof:} Let $\cK$ be an Erd\H{o}s-Rado class and $\bK$ be a large, $\bL_{\infty, \omega}$-elementary category.  Fix $f:\Card \to \Card$ witnessing that $\cK$ is Erd\H{o}s-Rado.  By virtue of being large, there is $M \in \bK$ such that $\|M\| \geq f(\mu)$, where $\mu$ is the size of fragment witnessing that $\bK$ is a $\bL_{\infty, \omega}$-elementary category.  Thus by Theorem \ref{gmott-ap-thm} and Chang's Presentation Theorem, there is a blueprint $\Phi \in \Upsilon^\cK[\bK]$.  Define a functor $F: \cK \to \bK$ by, for $I \in \cK$, $F(I) = EM_\tau(I, \Phi)$ and, for $f: I \to J$ in $\cK$, $Ff$ the map that takes $\sigma^{EM(I, \Phi)}(i_1, \dots, i_n)$ for a $\tau(\Phi)$-term $\sigma$ to $\sigma^{EM(J, \Phi)}\left(f(i_1) , \dots, f(i_n)\right)$.

This is clearly faithful.  Moreover, the $EM$ construction commutes with directed colimits, so $F$ preserves them.\hfill \dag\\

This proof works by noting that blueprints can be seen as well-behaved functors.  We can actually specify the properties of these functors to obtain a converse.  The one additional property that we need is that the size of $EM_\tau(I,\Phi)$ is determined by $|I|$ and an additional cardinal parameter representing $|\tau(\Phi)|$.  The following is based on an argument developed with John Baldwin in the case $\cK = \cK^{or}$.

For this, we need the following definition:
\begin{defin}
Fix an abstract class $\cK$.
\begin{enumerate}
    \item $\cK$ is \emph{universal} iff $\leq_{\cK}$ is $\subseteq_\tau$ and $\cK$ is closed under substructure.
    \item $\cK$ is a \emph{universal Erd\H{o}s-Rado class} iff it is universal and an Erd\H{o}s-Rado class.
\end{enumerate}
\end{defin}

Note that all universal classes $\cK$ are $\bL_{\infty, \omega}$-axiomatizable by saying no finite tuple generates a structure not in $\cK$.

\begin{theorem}\label{der-bp-thm}
Suppose $\cK$ is a universal Erd\H{o}s-Rado class and $\bK$ is a large, $\bL_{\infty, \omega}$-elementary category.  Let $F:\cK \to \bK$ be a faithful functor that preserves directed colimits such that there is a cardinal $\mu_F$ so that $\|F(I)\| = |I| + \mu_F$ for every $I \in \cK$.  Then there is a blueprint $\Phi \in \Upsilon^\cK_{\mu_F}[\bK]$ such that the functor induced by $I \in \cK \mapsto EM_\tau(I, \Phi)$ is naturally isomorphic to $F$.
\end{theorem}

{\bf Proof:}  Let $T \subset \bL_{\infty, \omega}(\tau)$ such that $\bK$ is (equivalent to) $\Mod T$. Enumerate the $\cK$-types as $\seq{p_i^n \in \S^n_\cK \mid i < \mu_n}$, and pick some $I_i^n \in \cK$ that is generated by elements $a^{i,n}_1, \dots, a^{i, n}_n$ that realize $p^n_i$.  We expand each $F(I^n_i)$ to a $\tau^*:= \tau(\bK) \cup \{F^n_\alpha: \alpha < \mu_F, n<\omega\}$-structure as in Shelah's Presentation Theorem.  In fact, we only give an explicit description of the $\{F^n_\alpha:\alpha<\mu_F\}$ structure on the $F(I^n_i)$: for each $n < \omega$ and $i < \mu_n$, define these functions so that $\{F^n_\alpha(a^{i, n}_1, \dots, a^{i,n}_n) : \alpha < \mu\}$ enumerates the universe of $F(I_i^n)$.  Then define the remaining functions arbitrarily.  

Since $F$ preserves directed colimits and $\cK$ is generated by the $I_i^n$ under directed colimits, we can lift these expansions to the rest of $F"\cK$.  Taking $I$ large enough, we can define a blueprint $\Phi \in \Upsilon^{\cK}[\tau^*]$.  For all $I$, the $\tau$-reduct of $EM_{\tau^*}(I, \Phi)$ is canonically isomorphic to $F(I)$.  Thus, $\Phi$ is as desired.\hfill\dag\\

Note that this converse requires that models be of a predictable size.  Specializing to linear order, we demand that $\Phi(n)$ be the same for all $n < \omega$.  This is necessary for the formalism we've given where $\tau(\Phi)$ consists of functions that can be applied to any element of $EM(I, \Phi)$.  To state the most general result, we could change this to only apply the functions of $\tau(\Phi)$ to the skeleton $I$.  Then the different sizes of $\Phi(n)$ could be dealt with by having different numbers of functions of different cardinalities.  But this seems like a marginal gain after what would be significant technical pain.  Additionally, the requirement that $\cK$ be universal can be removed.

\comment{Putting these results together, we have the following nice corollary.

\begin{cor}
Erd\H{o}s-Rado classes are precisely the large, finitely accessible categories that are ``minimal'' (in the sense Fact \ref{mor-mp-fact}).
\end{cor}}

We return to this category theoretic perspective in Section \ref{ind-coll-ssec} when discussing indiscernible collapse.  The existence of a minimal, large $\bL_{\infty, \omega_1}$-elementary category (and a version of EM models for those theories) is still open, although \cite[Section 4]{r-minimal} discusses this issue and places some restrictions on it.

\subsection{Generalized Shelah's Omitting Types Theorem}\label{gsott-ssec}

The application of Morley's Omitting Types Theorem to Abstract Elementary Classes is normally done through Shelah's Presentation Theorem, which gives a type omitting characterization of these classes (see Theorem \ref{gmott-ap-thm} for this argument).  Moving beyond this, Shelah has proved an omitting types theorem that strengthens this and specifically applies to to Abstract Elementary Classes in that it references Galois types rather than syntactic types (\cite{ms-compact} and \cite[Lemma 8.7]{sh394} both use some version of this).  The key addition is a reduction in the cardinal threshold for type omission at the cost of less control over what types are omitted.  The main combinatorial tool is still using the Erd\H{o}s-Rado Theorem to build Ehrenfeucht-Mostowski models, so we can similarly prove a version for any Erd\H{o}s-Rado class.

One nonstandard piece of notation is necessary.

\begin{defin}
Suppose $\bK$ is an AEC, and let $N \prec_\bK M$, $p \in \S_\bK(N)$, and $\chi \leq \|N\|$.  We say \emph{$M$ omits $p/E_\chi$} iff, for every $c \in M$, there is some $N_0 \prec_\bK N$ of size $< \chi$ such that $c$ does not realize $p \rest N_0$.
\end{defin}

\begin{theorem}[Generalized Shelah's Omitting Types Theorem]\label{gsott-thm}
Let $\cK$ be an Erd\H{o}s-Rado class and $\bK$ be an Abstract Elementary Class and $|\tau(\cK)| + \LS(\bK)\leq \chi \leq \lambda$ with
\begin{enumerate}
	\item $f(\mu) < \beth_{\LS(\bK)}(\mu)$ where $f$ witnesses that $\cK$ is Erd\H{o}s-Rado (for simplicity);
	\item $N_0 \prec_\bK N_1$ with $\|N_0\| \leq \chi$ and $\|N_1\| = \lambda$;
	\item $\Gamma_0 = \{p^0_i : i < i_0^*\} \subset \S_\bK(N_0)$; and
	\item $\Gamma_1 = \{p^1_i : i < i_1^*\} \subset \S_\bK(N_1)$ with $i_1^* \leq \chi$.
\end{enumerate}
Suppose that, for each $\alpha < \left(2^\chi\right)^+$, there is $M_\alpha \in \bK$ such that
\begin{enumerate}
	\item $f_\alpha^0: I^0_\alpha \to M_\alpha$ for $I^0_\alpha \in \cK_{\beth_\alpha(\lambda)}$;
	\item $N_1 \prec M_\alpha$;
	\item $M_\alpha$ omits $\Gamma_0$; and
	\item $M_\alpha$ omits $p^1_i/E_\chi$ for each $i <i^*_1$.
\end{enumerate}
Then there is $\Phi \in \Upsilon^\cK[\bK]$; increasing, continuous $\{N_q' \in \bK_{\leq \chi} \mid q \in \S_{\cK}^{inc}\}$; and increasing Galois types $p^1_{i, q} \in \S_\bK(N_q')$ for $q \in \S_\cK^{inc}$ and $i < i_1^*$ such that
\begin{enumerate}
	\item $N_0 = N_0' = EM_\tau(\emptyset, \Phi)$;
	\item for each $q \in \S_\cK^{inc}$, there is $\bi^q \in I \in \cK$ realizing $q$ such that $N_q' \prec_\bK EM_\tau(\bi^q, \Phi)$ and $f_q:EM_\tau(\bi^q, \Phi) \to M_{\alpha_q}$ for some $\alpha_q < (2^\chi)^+$ such that $f_q(N_q') \prec_\bK N_1$;
	\item $p^1_{i,q} = f_q^{-1}\left( p_i^1 \rest f_q(N_q')\right)$; and
	\item For every $I \in \cK_\omega$, $EM_\tau(I, \Phi)$ omits every type in $\Gamma_0$ and omits any type that extends $\{p_{i, q}^1 : q \in \S_\cK^{inc}\}$ in the following strong sense: if $a \in EM_\tau(I, \Phi)$ is in the $\tau(\Phi)$-closure of $\bi \in I$, then $a$ doesn't realize $H(p^1_{i, q})$, where $H:EM_\tau(\bi^q, \Phi) \cong EM_\tau(\bi,\Phi)$ is the lifting of $\bi \mapsto \bi^q$.
\end{enumerate}
\end{theorem}

The proof of the above adapts the proof of Shelah's Omitting Types Theorem just as Theorem \ref{gmott-thm} adapts Morley's original; see the notes by the author for a very detailed proof of (the ordinary) Shelah's Omitting Types Theorem \cite{b-sott-notes}.

\subsection{End-approximations} \label{er-notcer-ssec}

Generalized Morley's Omitting Types Theorem \ref{gmott-thm} is the primary source of Erd\H{o}s-Rado classes, but not the only source.  The proof that $\cK^{\omega-tr}$ is an Erd\H{o}s-Rado class uses the fact that elements of it can be approximated by trees of height $n$, and each $\cK^{n-tr}$ is combinatorial Erd\H{o}s-Rado class.  We capture this behavior in an abstract condition in Defintion \ref{endex-def} below, and Proposition \ref{wtr-endex-prop} shows that $\cK^{\omega-tr}$ fits into this framework.

\begin{defin}\label{endex-def}
We say that $\cK$ is \emph{end-approximated by combinatorial Erd\H{o}s-Rado classes} iff there are combinatorial Erd\H{o}s-Rado classes $\{\cK^n \mid n < \omega\}$ such that
\begin{enumerate}
    \item $\tau(\cK^n) \subset \tau(\cK^{n+1})$ and $\tau(\cK) = \cup_{n<\omega} \tau(\cK^n)$;
    \item there are functorial coherent restriction maps 
    $$\cdot \rest n : \cK^{\geq n} \to \cK^n$$
    where $\cK^{\geq n} = \cK \cup\bigcup_{k\geq n} \cK^k$ and 
    \begin{enumerate}
        \item `functorial' means that if $f:I \to J$ is a morphism in some $\cK^m$ for $m > n$ or in $\cK$, then 
        $$f\rest(I\rest n) : (I\rest n) \to (J\rest n)$$
        is a morphism in $\cK^n$ (note that `restriction map' includes that $I\rest n\subset I$)
        \item `coherent' means that for $n > m$ and $I \in \cK^{\geq n}$,
        $$I \rest m = \left( I\rest n\right) \rest m$$
        \item `restriction maps' has the normal meaning (restricting a structure to a smaller language) with the important detail that any sorts in $\tau(\cK)$ that are not $\tau(\cK^n)$ are removed from the structure, so the universe might shrink (see Proposition \ref{wtr-endex-prop}); in particular, $I\rest n \subset I$. 
    \end{enumerate}
    \item \label{endex-def-colimit} various structures involving $\cK$ are the (co)limit of the same structures on the $\cK^n$:
    \begin{enumerate}
        \item if $\seq{I_n \in \cK^n : n < \omega}$ is a sequence so $\left(I_{n+1}\right) \rest n = I_n$ for all $n < \omega$, then the $\tau$-structure 
        $$\bigcup_{n<\omega} I_n$$
        is in $\cK$;
        \item in the above, if each $I_n$ is $\mu$-big, then the union is $\mu$-big;
        \item if $\left\langle p^n = \tp_{\cK^n}\left(\ba^n; I_n\right) \in \S_{\cK^n} \mid k_0 \leq n < \omega \right\rangle$ is a sequence for some $k_0<\omega$ so $\ba^{n+1} \in I_{n+1}\rest n$ and
        $$p^n = p^{n+1}\rest n := \tp_{\cK^{n+1}}\left(\ba^{n+1}; I_{n+1}\rest n\right)$$
        for all $n \geq k_0$, then there is a unique $p = \tp_{\cK}(\ba; I) \in \S_\cK$ such that
        $$p^n = p \rest n := \tp_{\cK}(\ba; I \rest n)$$
    \end{enumerate}
    \item \label{big-rest-cond} the restriction of a $\mu$-big model (with bigness computed in the domain) is $\mu$-big (in the restricted class); and
    \item \label{lift-cond} if $I \in \K^n$ and $J \in \cK$ are both $\mu$-big (in their respective contexts) and there is a $\cK^n$-embedding
    $$f:I \to J\rest n$$
    then there is a lift that consists of $\hat{I} \in \cK$ that is $\mu$-big and a $\cK$-morphism $\hat{f}:\hat{I}\to J$ such that $\hat{I}\rest n = I$ and $\hat{f}\rest I = f$.
    \item \label{type-det-cond} any $p \in \S(\cK)$ is determined at some finite stage $k_p < \omega$, which means that
    \begin{enumerate}
        \item $p$ is the unique extension of $p \rest k_p$ to $\S_\cK$ and, furthermore, for any $n \geq k_p$, $p \rest n$ is the unique extension of $p \rest k_p$ to $\S_{\cK^n}$
        \item if $I\in  \cK$ and $\ba \in I$ realizes $p$, then $\ba \in (I\rest k_p)$
        \item for $p \in \S_\cK^0$, $k_p = 0$
        \item for $s \subset \ell(p)$ and $k \leq k'$, we have
        \begin{eqnarray*}
             k_{p^s} &\leq& k_p\\
            (p \rest k')^s \rest k &=& (p^s) \rest k
        \end{eqnarray*}
    \end{enumerate}
\end{enumerate}
\end{defin}

Note that since the language is finitary, we can always decompose the objects covered in Definition \ref{endex-def}.(\ref{endex-def-colimit}) as the canonical colimits of it's restrictions; that is, if $I \in \cK$, then
$$I = \bigcup_{n<\omega} \left(I\rest n\right)$$

The lifting condition Definition \ref{endex-def}.(\ref{lift-cond}) and the type determination condition Defintion \ref{endex-def}.(\ref{type-det-cond}) are the key properties.  Our initial motivation for this framework was to show that $\cK^{\omega-tr}$ is an Erd\H{o}s-Rado class.  After sending him a draft, Baldwin pointed us to \cite[Section 4]{bs-atomic-stability}, where this is already shown.  We hope this framework can be applied in other situations.

\begin{prop}\label{wtr-endex-prop}
$\K^{\omega-tr}$ is end approximated by $\{\K^{n-tr} \mid n < \omega\}$.
\end{prop}

{\bf Proof:} The proof is straightforward.  The truncation map $\cdot \rest n$ truncates a tree of height $\geq n$ to its $\leq n$ levels.  Any $p \in S^n_{\cK^{\omega-tr}}$ specifies the max height $k_p$ of a realization.

For condition (\ref{lift-cond}), let $I, J, f$ be as there.  Build $\hat{I}$ by specifying $\hat{I} \rest n = I$ and, given maximal $\eta \in I$, the successors of $\eta$ in $\hat{I}$ are an isomorphic copy the successors of $f(\eta)$ in $J$.  Since $J$ is at least $\alpha$-splitting, so is $\hat{I} \in \cK^{n+1}_\alpha$, and the isomorphisms give the lift $\hat{f}:\hat{I} \to J \rest n+1$.

For condition (\ref{type-det-cond}), given $p \in \S(\cK)$, set $k_p$ to the maximum height of the realizations.  This is finite and $p \rest k_p$ determines all of the information in the type.\hfill \dag\\

The other natural candidate is that $\cK^{(\omega,\sigma)-hg}$ might be end approximated by $\{\cK^{(k, \sigma)-hg}\mid k<\omega\}$.  However, the crucial lifting condition fails: we can arrange a saturated/big $H \in \cK^{(3, \sigma)-tr}$ and pick a subgraph $G \subset H$ that is saturated/big for unary types that the $\cK^{(2, \sigma)-hg}$ part can see, but not for the binary relations in $\cK^{(3, \sigma)-hg}$.  The author suspects that more refined methods could still prove $\cK^{(\omega, \sigma)-hg}$ is an Erd\H{o}s-Rado class, but we don't do that here.

\comment{Similarly, we can show the following.

\begin{prop}
$\cK^{(\omega, \sigma)-hg}$ is end approximated by $\{\cK^{(k, \sigma)-hg} \mid k<\omega\}$
\end{prop}}

\begin{theorem} \label{endapp-er-thm}
Let $\K$ be end-approximated by combinatorial Erd\H{o}s-Rado classes $\{\cK^n \mid n < \omega\}$.  Then $\K$ is a cofinal Erd\H{o}s-Rado class.  If $F$ is a witnessing function of all $\cK^n$, then $f^{\cK, F}$ is a witnessing function for $\cK$ (recall Definition \ref{daleth-wit-def}).
\end{theorem}

The goal is to repeat the proof of the Generalized Morley's Omitting Types Theorem \ref{gmott-thm}, except we restrict our set-up to the $n$th approximation $\cK^n$ in stage $n$.  Then when we move to stage $n+1$, we use the lifting condition Definition \ref{endex-def}.(\ref{lift-cond}) to lift the set-up to the next level. Crucially, the `approximate blueprints' $\Phi_n$ have all types of length $\leq n$ in $\cK^n$ and are \emph{not} strictly increasing.  Instead, the $\Phi_n$ are increasing on the types $p$ satisfying $k_p\leq n$. 

At the sage advice of the referee, we repeat the proof here and provide all details so that we can highlight where the different conditions of Definition \ref{endex-def} are used.\\

{\bf Proof:} As in the proof of Theorem \ref{gmott-thm}, set $\daleth_\alpha = \daleth_\alpha^F$, $g=g^\cK$, and $f=f^{\cK, F}$ and we can begin with $f_\alpha:I_\alpha\to M_\alpha$ for $\alpha < g(\mu)^+$ so $I_\alpha$ is $\daleth_\alpha$-big and all $M_\alpha$ realize the same $0$-type, $p^* \in \S^0_\tau$.

We build the following for $n<\omega$ and $\alpha<g(\mu)^+$.  Note that we have tried to match the corresponding names from the proof of Theorem \ref{gmott-thm}; the $g^n_\alpha$ functions are new since the $I^n_\alpha$ are in $\cK^n$, so cannot literally be $I_{\beta_n(\alpha)}$.
\begin{itemize}
    \item $\Phi_n: \S^{inc,\leq n}_{\cK^n}\to \S^{\leq n}_\tau$;
    \item $\beta_n(\alpha) < g(\mu)^+$;
    \item $\gamma_{n+1}(\alpha) < g(\mu)^+$;
    \item $\daleth_\alpha$-big $I^n_\alpha \in \cK^n$;
    \item $\cK^n$-embeddings $h^{n+1}_\alpha:(I^{n+1}_\alpha \rest n) \to I^n_{\gamma_{n+1}(\alpha)}$;
    \item $\cK^n$-embeddings $g^n_\alpha: I^n_\alpha\to (I_{\beta_n(\alpha)}\rest n)$; and
    \item $f^n_\alpha:I^n_\alpha\to M_{\beta_n(\alpha)}$
\end{itemize}
such that
\begin{enumerate}
    \item \label{0-cond'}$\beta_0(\alpha) = \alpha$, $I^0_\alpha = I_\alpha\rest0$, $g^0_\alpha = \id_{I^0_\alpha}$, and $f^0_\alpha = f_\alpha\rest I^0_\alpha$;
    \item \label{phi-cond'} for all $\alpha<g(\mu)^+$, $i_1<\dots<i_m \in I^n_\alpha$, and $m \leq n$, we have
    $$\Phi_n\left(\tp_{\cK^n}\left(i_1, \dots, i_m;I^n_\alpha\right)\right)=\tp_\tau\left(f^n_\alpha(i_1), \dots, f^n_\alpha(i_m); M_{\beta_n(\alpha)}\right)$$
    \item \label{coh-cond'}the $\Phi_n$ are coherent in the following sense:
    \begin{enumerate}
        \item \label{coh-cond'-1}(within $n$) if $p \in \S^{inc, \leq n}_{\cK^n}$ and $s \subset \ell(p)$, then
        $$\Phi_n(p)^s = \Phi_n(p^s)$$
        \item \label{coh-cond'-2}(across $n$) if $p \in \S^{inc,\leq n}_{\cK^n}$ has a unique extension to $p^* \in \S^{inc, \leq n}_{\cK}$ with\footnote{This implies that $p=p^*\rest n$ is the unique extension of $p\rest k_{p^*}=p^*\rest k_{p^*}$.} $k_{p^*} \leq n$, then  
        $$\Phi_{k_{p^*}}(p \rest k_{p^*}) = \Phi_n(p)$$
    \end{enumerate}
    \item \label{diag-cond'} for $\alpha < g(\mu)^+$ and $n<\omega$,
    \begin{enumerate}
        \item $\alpha \leq \beta_n(\alpha)$, $\alpha \leq \gamma_{n+1}(\alpha)$, and $\beta_n(\gamma_{n+1}(\alpha)) = \beta_{n+1}(\alpha)$; and
        \item\label{diag-cond'-2} the following diagram commutes:
% https://q.uiver.app/#q=WzAsNCxbMiw0LCJJXntuKzF9X1xcYWxwaGEgXFxyZXN0IG4iXSxbMiwyLCJJXm5fe1xcZ2FtbWFfe24rMX0oXFxhbHBoYSl9Il0sWzQsMCwiTV97XFxiZXRhX3tuKzF9KFxcYWxwaGEpfSJdLFswLDAsIklfe1xcYmV0YV97bisxfShcXGFscGhhKX0gXFxyZXN0IG4iXSxbMCwxLCJoXntuKzF9X1xcYWxwaGEiXSxbMSwyLCJmXm5fe1xcZ2FtbWFfe24rMX0oXFxhbHBoYSl9Il0sWzEsMywiZ15uX3tcXGdhbW1hX3tuKzF9KFxcYWxwaGEpfSIsMl0sWzMsMiwiZl97XFxiZXRhX3tuKzF9KFxcYWxwaGEpfVxccmVzdCBcXGxlZnQoSV97XFxiZXRhX3tuKzF9KFxcYWxwaGEpfSBcXHJlc3QgblxccmlnaHQpIiwwLHsiY3VydmUiOi00fV0sWzAsMywiZ157bisxfV9cXGFscGhhXFxyZXN0IFxcbGVmdChJXntuKzF9X1xcYWxwaGEgXFxyZXN0IG5cXHJpZ2h0KSIsMCx7ImN1cnZlIjotNH1dLFswLDIsImZee24rMX1fXFxhbHBoYVxccmVzdCBcXGxlZnQoSV57bisxfV9cXGFscGhhIFxccmVzdCBuXFxyaWdodCkiLDIseyJjdXJ2ZSI6NH1dXQ==
\[\begin{tikzcd}
	{I_{\beta_{n+1}(\alpha)} \rest n} &&&& {M_{\beta_{n+1}(\alpha)}} \\
	\\
	&& {I^n_{\gamma_{n+1}(\alpha)}} \\
	\\
	&& {I^{n+1}_\alpha \rest n}
	\arrow["{f_{\beta_{n+1}(\alpha)}\rest \left(I_{\beta_{n+1}(\alpha)} \rest n\right)}", curve={height=-24pt}, from=1-1, to=1-5]
	\arrow["{g^n_{\gamma_{n+1}(\alpha)}}"', from=3-3, to=1-1]
	\arrow["{f^n_{\gamma_{n+1}(\alpha)}}", from=3-3, to=1-5]
	\arrow["{g^{n+1}_\alpha\rest \left(I^{n+1}_\alpha \rest n\right)}", curve={height=-24pt}, from=5-3, to=1-1]
	\arrow["{f^{n+1}_\alpha\rest \left(I^{n+1}_\alpha \rest n\right)}"', curve={height=24pt}, from=5-3, to=1-5]
	\arrow["{h^{n+1}_\alpha}", from=5-3, to=3-3]
\end{tikzcd}\]
    \end{enumerate}
\end{enumerate}

{\bf This is enough:} We want to define our blueprint as the colimit of the $\Phi_n$'s, but face the extra difficulty that they are not actually increasing: $\S^{inc, \leq n}_{\cK^n}$ consists of types in a different class then $\S^{inc, \leq n+1}_{\cK^{n+1}}$.  Moreover, a type in the first set might have multiple extensions to the second.  We solve this problem by using Defintion \ref{endex-def}.(\ref{type-det-cond}), which says that all types in $\S^{inc}_{\cK}$ are determined at some level.

Fix $p \in \S^{inc}_\cK$.  This is determined at some level $k_p <\omega$.  Define
$$\Phi(p) := \Phi_{k_p}(p \rest k_p)$$
First, for any $n \geq k_p$, we claim that 
$$\Phi_{k_p}(p\rest k_p) = \Phi_n(p\rest n)$$
This means that we have chosen $\Phi(p)$ to be the output of the $\Phi_n$'s after it has stabilized.  We prove this by induction on $n \geq k_p$ by using condition (\ref{phi-cond'}) of the construction:  For $n=k_p$, it is immediate. So assume we know $\Phi_n(p\rest n) = \Phi_{k_p}(p\rest k_p)$ and we want to prove this for $n+1$.

To this end, let $i_1<\dots<i_m\in I^{m+1}_\alpha$ for some $\alpha<g(\mu)^+$ realize $p\rest (n+1)$. By Definition \ref{endex-def}.(\ref{type-det-cond}), we have $i_1, \dots, i_m\in \left(I^{n+1}_\alpha\rest n\right)$, so it is in the domain of $h^{n+1}_\alpha$.  Moreover, 
$$h^{n+1}_\alpha(i_1), \dots, h^{n+1}_\alpha(i_m)\in I^n_{\gamma_{n+1}(\alpha)}$$
realizes $p\rest n$.  By the lower right triangle in the diagram of condition (\ref{diag-cond'-2}) and condition (\ref{phi-cond'}), we have
\begin{eqnarray*}
    \Phi_{n+1}\left(p\rest(n+1)\right) &=& \tp_\tau\left(f^{n+1}_\alpha(i_1), \dots, f^{n+1}_\alpha(i_m); M_{\beta_{n+1}(\alpha)}\right)\\
    &=&\tp_\tau\left(f^n_{\gamma_{n+1}(\alpha)}\left(h^{n+1}_\alpha(i_1)\right),\dots,f^n_{\gamma_{n+1}(\alpha)}\left(h^{n+1}_\alpha(i_1)\right);M_{\beta_{n+1}(\alpha)}\right)\\
    &=&\Phi_n\left(\tp_\tau\left(h^{n+1}_\alpha(i_1), \dots, h^{n+1}_\alpha(i_m);I^n_{\gamma_{n+1}(\alpha)}\right)\right)\\
    &=&\Phi_n(p\rest n)=\Phi_{k_p}(p\rest k_p)
\end{eqnarray*}

Now we want to check the coherence condition from Definition \ref{genbp-def}.  Fix $p \in \S^{inc,\leq n}_\cK$ and $s \subset m=\ell(p)$.  By Definition \ref{endex-def}.(\ref{type-det-cond}).(4), $k_{p^s} \leq k_p$ and $p^s \rest k_{p^s} = \left( p\rest k_p \right)^s \rest k_{p^s}$.  Then we have
\begin{eqnarray*}
    \Phi(p)^s &=& \Phi_{k_p}(p \rest k_p)^s\\
    &=& \Phi_{k_p}\left((p\rest k_p)^s \right) \text{ by condition (\ref{coh-cond'-2})}\\
    &=& \Phi_{k_{p^s}}\left((p\rest k_p)^s \rest k_{p^s} \right) \text{ by condition (\ref{coh-cond'-1})}\\
    &=& \Phi_{k_{p^s}}(p^s \rest k_{p^s}) = \Phi(p^s)
\end{eqnarray*}

{\bf Construction:} \underline{$n=0$:} As before, this is determined by condition (\ref{0-cond'}).\\

\underline{$n+1$:}  Suppose we have the construction at stage $n$.  The construction is similar to the construction in Theorem \ref{gmott-thm} in plan, but there are many additional technicalities, so we give the construction in full.

Fix $\alpha < g(\mu)^+$.  We have a $\cK^n$-map $g^n_\alpha:I^n_\alpha \to I_{\beta_n(\alpha)}\rest n$ with $I^n_\alpha$ being $\daleth_\alpha$-big (in $\cK^n$) and $I_{\beta_n(\alpha)}$ being $\daleth_{\beta_n(\alpha)}$-big (in $\cK$).  By the lifting condition Definition \ref{endex-def}.(\ref{lift-cond}) and since $\alpha \leq \beta_n(\alpha)$, there is a lift
$$\hat{g}^{n+1}_\alpha:\hat{I}^{n+1}_\alpha \to I_{\beta_n(\alpha)}\rest (n+1)$$ with $\hat{I}^{n+1}_\alpha \in \cK^{n+1}$ being $\daleth_\alpha$-big so that $\hat{I}^{n+1}_\alpha \rest n = I^n_\alpha$ and $\hat{g}^{n+1}_\alpha\rest (I^n_\alpha) = g^n_\alpha$.  We have that $F^*(\daleth_\alpha, n+1) \leq \daleth_{\alpha+1}$.  Since $F^*$ is the function associated with the `$\leq n$' partition relation (recall Proposition \ref{leqn-prop}), we color
$$c^{n+1}_\alpha:\left[\hat{I}^{n+1}_{\alpha+1}\right]^{\leq n+1} \to \S^{\leq n+1}_\tau$$
by, for $i_1 < \dots <i_m \in \hat{I}^{n+1}_{\alpha+1}$ with $m \leq n+1$, setting
$$c^{n+1}_{\alpha}\left(\{i_1, \dots, i_m\}\right) = \tp_\tau\left(f_{\beta_n(\alpha+1)}\circ\hat{g}^{n+1}_{\alpha+1}(i_1), \dots,f_{\beta_n(\alpha+1)}\circ\hat{g}^{n+1}_{\alpha+1}(i_m);M_{\beta_n(\alpha+1)} \right)$$
By Proposition \ref{leqn-prop}, we can find $\daleth_\alpha$-big $\bar{I}^{n+1}_\alpha \in \cK^{n+1}$ and
\begin{eqnarray*}
    \bar{h}^{n+1}_\alpha&:&\bar{I}^{n+1}_\alpha \to \hat{I}^{n+1}_{\alpha+1}\\
    c^{*, n+1}_\alpha&:& \S^{inc, \leq n+1}_{\cK^{n+1}} \to \S^{\leq n+1}_\tau
\end{eqnarray*}
that witness the type-homogeneity, that is, so for all $i_1<\dots<i_m \in \bar{I}^{n+1}_\alpha$ with $m \leq n+1$, we have
$$c^{n+1}_\alpha\left(\left\{\bar{h}^{n+1}_\alpha(i_1), \dots, \bar{h}^{n+1}_\alpha(i_m)\right\}\right) = c^{*, n+1}_\alpha\left(\tp_{\cK^{n+1}}\left(i_1, \dots, i_m; \bar{I}^{n+1}_\alpha\right)\right)$$

Now we can thin out and collapse these structures as before: each $c^{*, n+1}_\alpha$ is a function $\S^{inc, \leq n+1}_{\cK^{n+1}}\to \S^{\leq n+1}_\tau$ and there are less than $\cf(g(\mu)^+)=g(\mu)^+$ of these, so there is $X \subset g(\mu)^+$ of size $g(\mu)^+$ and $c^{*, n+1}:\S^{inc, \leq n+1}_{\cK^{n+1}}\to\S^{\leq n+1}_\tau$ such that $c^{*, n+1} = c^{*, n+1}_\alpha$ for all $\alpha \in X$.

Set $\pi:X \cong g(\mu)^+$ to be the collapse.  Now define the following:
\begin{itemize}
    \item $\Phi_{n+1} = c^{*, n+1}$;
    \item $\beta_{n+1}(\alpha) = \beta_n(\pi^{-1}(\alpha)+1)$;
    \item $\gamma_{n+1}(\alpha) = \pi^{-1}(\alpha)+1$;
    \item $I^{n+1}_\alpha = \bar{I}^{n+1}_{\pi^{-1}(\alpha)} \in \cK^{n+1}$ is $\daleth_{\pi^{-1}(\alpha)}$-big (and therefore $\daleth_\alpha$-big);
    \item $h^{n+1}_\alpha = \bar{h}^{n+1}_{\pi^{-1}(\alpha)} \rest \left(I^{n+1}_\alpha \rest n\right)$;
    \item $g^{n+1}_\alpha = \hat{g}^{n+1}_{\pi^{-1}(\alpha)+1} \circ \bar{h}^{n+1}_{\pi^{-1}(\alpha)}$; and
    \item $f^{n+1}_\alpha = f_{\beta_n\left(\pi^{-1}(\alpha)+1\right)}\circ \hat{g}^{n+1}_{\pi^{-1}(\alpha)+1} \circ \bar{h}^{n+1}_{\pi^{-1}(\alpha)}$, where we could have written the first composite as `$f_{\beta_n\left(\pi^{-1}(\alpha)+1\right)}\rest\left(I_{\beta_n\left(\pi^{-1}(\alpha)+1\right)}\rest(n+1)\right)$' for more precision and more notation.
\end{itemize}
Now we verify the various conditions of our construction:
\begin{enumerate}
    \item This is only relevant to the base case.
    \item Let $\alpha < g(\mu)^+$ and $i_1<\dots<i_m\in I^{n+1}_\alpha=\bar{I}^{n+1}_{\pi^{-1}(\alpha)}$ with $m \leq n+1$.  Then $\pi^{-1}(\alpha)\in X$, so $\Phi_{n+1}=c^{*, n+1} = c^{*, n+1}_{\pi^{-1}(\alpha)}$ witnesses the type-homogeneity of $\bar{I}^{n+1}_{\pi^{-1}(\alpha)}$ for $c^{n+1}_{\pi^{-1}(\alpha)}$.  Thus, we compute
    \begin{eqnarray*}
        \Phi_{n+1}\left(\tp_{\cK^{n+1}}\left(i_1, \dots, i_m; I^{n+1}_\alpha\right)\right) &=& c^{*, n+1}_{\pi^{-1}(\alpha)}\left( \tp_{\cK^{n+1}}\left(i_1, \dots, i_m; \bar{I}^{n+1}_{\pi^{-1}(\alpha)}\right)\right)\\
        &=& c^{n+1}_{\pi^{-1}(\alpha)} \left( \left\{ \bar{h}^{n+1}_{\pi^{-1}(\alpha)}(i_1), \dots,\bar{h}^{n+1}_{\pi^{-1}(\alpha)}(i_m)\right\}\right)\\
        &=& \tp_\tau\left(f_{\beta_n(\pi^{-1}(\alpha)+1)}\circ\hat{g}^{n+1}_{\pi^{-1}(\alpha)+1}\circ \bar{h}^{n+1}_{\pi^{-1}(\alpha)}(i_1), \dots,\right.\\
        & &\left.f_{\beta_n(\pi^{-1}(\alpha)+1)}\circ\hat{g}^{n+1}_{\pi^{-1}(\alpha)+1}\circ \bar{h}^{n+1}_{\pi^{-1}(\alpha)}(i_m);M_{\beta_n(\pi^{-1}(\alpha)+1)} \right)\\
        &=& \tp_\tau\left(f^{n+1}_\alpha(i_1), \dots, f^{n+1}_{\alpha}(i_m);M_{\beta_{n+1}(\alpha)}\right)
    \end{eqnarray*}
    as desired.
    \item \begin{enumerate}
        \item This follows from condition (\ref{phi-cond'}), but to be explicit: let $p \in \S^{inc, \leq n+1}_{\cK^{n+1}}$ and $s \subset m = \ell(p)$.  Fix some/any $\alpha < g(\mu)^+$ and $i_1<\dots<i_m \in I^{n+1}_\alpha$ such that
    $$p=\tp_{\cK^{n+1}}\left(i_1, \dots, i_m; I^{n+1}_\alpha\right)$$
        If we write $s = \{i_{j_1}, \dots, i_{j_\ell}\}$, then 
        \begin{eqnarray*}
        \Phi_{n+1}(p^s) &=& \tp_\tau\left(f_\alpha^{n+1}(i_{j_1}), \dots, f^{n+1}_\alpha(i_{j_\ell}); M_{\beta_{n+1}(\alpha)}\right)\text{ by condition (\ref{phi-cond'})}\\
        \Phi_{n+1}(p) &=& \tp_\tau\left(f^{n+1}_\alpha(i_1), \dots, f_{\alpha}^{n+1}(i_m) ; M_{\beta_{n+1}(\alpha)}\right) \text{ by condition (\ref{phi-cond'})}\\
        \Phi_{n+1}(p)^s &=& \tp_\tau\left(f_\alpha^{n+1}(i_{j_1}), \dots, f^{n+1}_\alpha(i_{j_\ell}); M_{\beta_{n+1}(\alpha)}\right)\text{ by definition of this operation}
        \end{eqnarray*}
        So $\Phi_{n+1}(p^s) = \Phi_{n+1}(p)^s$ as desired.
        \item Suppose $p \in \S^{inc, \leq n+1}_{\cK^{n+1}}$ has a unique extension to $p^* \in \S_\cK^{inc, \leq n+1}$ with $k_{p^*}\leq n+1$.\\
        If $k_{p^*} = n+1$, then
        $$\Phi_{k_{p^*}}(p\rest k_{p^*}) = \Phi_{n+1}\left(p \rest (n+1)\right) = \Phi_{n+1}(p)$$
        If $k_{p^*} \leq n$, then $p^*$ is also the unique extension of $p\rest n$.  By induction, we have
        $$\Phi_{k_{p^*}}((p\rest n)\rest k_{p^*}) = \Phi_{k_{p^*}}(p\rest k_{p^*}) = \Phi_n(p\rest n)$$
        Let $\alpha<g(\mu)^+$ and $i_1<\dots<i_m \in I^{n+1}_\alpha$ such that
        $$p =\tp_{\cK^{n+1}}(i_1, \dots,i_m; I^{n+1}_\alpha)$$
        Since $k_{p^*}\leq n$, $i_1, \dots,i_m\in I^{n+1}_\alpha\rest n$ by Definition \ref{endex-def}.\ref{endex-def-colimit}.(c).  Thus, $h_\alpha^{n+1}(i_1), \dots, h_\alpha^{n+1}(i_m)\in I^n_{\gamma_{n+1}(\alpha)}$.  Then we have
        \begin{eqnarray*}
            p\rest n &=& \tp_{\cK^n} (i_1, \dots, i_m; I_\alpha^{n+1}\rest n)\\
            &=& \tp_{\cK^n}(h^{n+1}_\alpha(i_1), \dots, h_\alpha^{n+1}(i_m); I^n_{\gamma_{n+1}(\alpha)})
        \end{eqnarray*}
        Thus, by condition (\ref{phi-cond'}) applied at $n$ and $n+1$ (since we have already proved it in this induction), we have
        \begin{eqnarray*}
            \Phi_n(p\rest n) &=&\tp_\tau\left(f^n_{\gamma_{n+1}(\alpha)}\left(h^{n+1}_\alpha(i_1)\right), \dots, f^n_{\gamma_{n+1}(\alpha)}\left(h^{n+1}_\alpha(i_m)\right);M_{\beta_{n}\left(\gamma_{n+1}(\alpha)\right)}\right)\\
            \Phi_{n+1}(p) &=& \tp_\tau \left(f_\alpha^{n+1}(i_1), \dots,f_\alpha^{n+1}(i_m); M_{\beta_{n+1}(\alpha)}\right)
        \end{eqnarray*}
        These types on the right-hand side are the same: $M_{\beta_{n+1}(\alpha)}=M_{\beta_n\left(\gamma_{n+1}(\alpha)\right)}$ by construction and by (the independently proven) condition (\ref{diag-cond'}).(\ref{diag-cond'-2}), 
        $$f^{n+1}_\alpha = f^n_{\gamma_{n+1}(\alpha)}\circ h^{n+1}_\alpha \text{ on }I^{n+1}_\alpha \rest n$$
        Thus we have shown the desired equality
        $$\Phi_{n+1}(p) = \Phi_n(p\rest n) = \Phi_{k_{p^*}}(p\rest k_{p^*})$$
    \end{enumerate}
    \item for $\alpha < g(\mu)^+$ and $n<\omega$,
    \begin{enumerate}
        \item we have
        \begin{eqnarray*}
            \beta_{n+1}(\alpha) &=& \beta_{n}\left(\pi^{-1}(\alpha)+1\right) \geq \pi^{-1}(\alpha)+1>\alpha\\
            \gamma_{n+1}(\alpha) &=& \pi^{-1}(\alpha)+1 > \alpha
        \end{eqnarray*}
        \item to show the diagram commutes, we repeat the desired diagram and make substitutions for the definitions above
        \[\begin{tikzcd}
	{I_{\beta_{n}(\pi^{-1}(\alpha)+1)} \rest n} &&&& {M_{\beta_{n}(\pi^{-1}(\alpha)+1)}} \\
	\\
	&& {I^n_{\pi^{-1}(\alpha)+1}} \\
	\\
	&& {\bar{I}^{n+1}_{\pi^{-1}(\alpha)} \rest n}
	\arrow["{f_{\beta_{n}(\pi^{-1}(\alpha)+1)}\rest \left(I_{\beta_{n}(\pi^{-1}(\alpha)+1)} \rest n\right)}", curve={height=-24pt}, from=1-1, to=1-5]
	\arrow["{g^n_{\pi^{-1}(\alpha)+1}}"', from=3-3, to=1-1]
	\arrow["{f^n_{\pi^{-1}(\alpha)+1}}", from=3-3, to=1-5]
	\arrow["{\left(\hat{g}^{n+1}_{\pi^{-1}(\alpha)+1} \circ \bar{h}^{n+1}_{\pi^{-1}(\alpha)}\right) \rest \left(\bar{I}^{n+1}_{\pi^{-1}(\alpha)} \rest n\right)}", curve={height=-24pt}, from=5-3, to=1-1]
	\arrow["{\left(f_{\beta_n(\pi^{-1}(\alpha)+1)}\circ \hat{g}^{n+1}_{\pi^{-1}(\alpha)+1}\circ \bar{h}^{n+1}_{\pi^{-1}(\alpha)}\right)\rest \bar{I}^{n+1}_{\pi^{-1}(\alpha)}}"', curve={height=24pt}, from=5-3, to=1-5]
	\arrow["{\bar{h}^{n+1}_{\pi^{-1}(\alpha)} \rest \bar{I}^{n+1}_{\pi^{-1}(\alpha)}}", from=5-3, to=3-3]
\end{tikzcd}\]
First note that the upper triangle commutes by the inductive assumption and the outer triangle commutes by definition.  The lower left triangle commutes by the lifting property: 
$$\hat{g}^{n+1}_{\pi^{-1}(\alpha)+1} \rest \left(I^n_{\pi^{-1}(\alpha)+1}\right) = g^{n}_{\pi^{-1}(\alpha)+1}$$
Since the other parts of the diagram commute, the commutation of the lower right triangle follows (although it could be shown directly).
    \end{enumerate}
\end{enumerate}
\hfill\dag\\

Thus, while we have no combinatorial partition result for $\cK^{\omega-tr}$, it is an Erd\H{o}s-Rado class.

\begin{corollary}\label{wtree-er-cor}
$\K^{\omega-tr}$ is a cofinal Erd\H{o}s-Rado class witnessed by $\mu \mapsto \beth_{\left(2^\mu\right)^+}$.
\end{corollary}

{\bf Proof:} By Theorem \ref{endapp-er-thm} applied to Proposition \ref{wtr-endex-prop}.\hfill\dag\\

\section{Applications}\label{app-sec}

\subsection{Unsuperstability in Abstract Elementary Classes} \label{unsuper-ssec}

In countable first-order theories, strict stability can detected by counting types at cardinals $\lambda$ satisfying $\lambda<\lambda^\omega$: 
\begin{center}
$T$ is superstable iff $T$ is stable in some cardinal $\lambda$ with $\lambda<\lambda^\omega$ iff $T$ is stable in all cardinals $\lambda$ with $\lambda\geq2^\omega$.
\end{center}
This is done by building what is called a `Shelah tree' \cite[p. 85]{b-fund-stab}.  This is a way of embedding the $\omega+1$-height tree ${}^{\leq \omega}\lambda$ into a model of $T$ so the types of branches are differentiated over their initial segments.  In the context of nonelementary classes, Baldwin and Shelah \cite[Theorem 3.3]{bs-atomic-stability} generalized this to \emph{atomic} classes by use of $\cK^{\omega-tr}$-indiscernibles.

Here, we generalize this to tame Abstract Elementary Classes with amalgamation.  Note that we break our convention of always using types over the empty set here.  In fact, we will consider types over arbitrary sets.  Given $N \in \bK$, $a \in N$, and $B \subset N$, we write 
$$\tp_\bK(a/B; N)$$
for the equivalence class of triples where the equivalence relation in Definition \ref{stone-def} is required to fix the parameter set $B$ as well. The presence of amalgamation means that we can always witness $\bK$-type equality with $n=1$ in Definition \ref{stone-def}.(1).  Since we lack a monster model, we fix an ambient model $N$ to give meaning to $B$.  Then $\S_{\bK}(B; N)$ is the collection of all Galois types over $B$ as seen as a set in $N$ and we write $\S_\bK(N)$ for $\S_\bK(N;N)$. We say that $\bK$ is Galois stable in $\lambda \geq \LS(\bK)$ iff for every $M \in \bK_\lambda$, we have $|\S_\bK(M)|\leq \lambda$ and Galois unstable in $\lambda$ for the negation of that statement. We omit other basics of Abstract Elementary Classes, but the key definitions can be found in one of \cite{ramibook, baldwinbook, bv-survey}.

\begin{theorem}\label{unsuper-thm}
Let $\bK$ be a $<\kappa$-tame Abstract Elementary Class with amalgamation (allowing for the possibility that $\kappa < \LS(\bK)$).  One of the following holds:
\begin{enumerate}
    \item there is $\chi < \beth_{\left(2^{\kappa+\LS(\bK)}\right)^+}$ such that for all $M \in \bK_{\geq\chi}$, $|\S_\bK(M)|\leq \|M\|^{<\kappa}$; or
    \item $\bK$ is Galois unstable in every $\lambda$ satisfying $\lambda^\omega > \lambda \geq \kappa+\LS(\bK)$.
\end{enumerate}
\end{theorem}

The first case roughly corresponds to Galois superstability, while the second is not superstable.  However, the necessary involvement of $\|M\|^{<\kappa}$ in the type counting (as opposed to just $\|M\|$) makes comparison with other results awkward.  Here we state two corollaries that rephrase the result directly in terms of Galois stability:

\begin{cor}\
\begin{enumerate}
    \item Let $\bK$ be a $<\omega$-tame Abstract Elementary Class with amalgamation.  One of the following holds:
\begin{enumerate}
    \item there is $\chi < \beth_{\left(2^{\kappa+\LS(\bK)}\right)^+}$ such that $\bK$ is Galois stable in every $\lambda > \chi$; or
    \item $\bK$ is Galois unstable in every $\lambda$ satisfying $\lambda^\omega > \lambda \geq \kappa+\LS(\bK)$.
\end{enumerate}
    \item Let $\bK$ be a $<\kappa$-tame Abstract Elementary Class with amalgamation (allowing for the possibility that $\kappa < \LS(\bK)$).  One of the following holds:
\begin{enumerate}
    \item there is $\chi < \beth_{\left(2^{\kappa+\LS(\bK)}\right)^+}$ such that $\bK$ is Galois stable in every $\lambda > \chi$ so $\lambda^{<\kappa}=\lambda$; or
    \item $\bK$ is Galois unstable in every $\lambda$ satisfying $\lambda^\omega > \lambda \geq \kappa+\LS(\bK)$.
\end{enumerate}
\end{enumerate}
\end{cor}

The subscript `$\left(2^{\kappa+\LS(\bK)}\right)^+$' can be replaced by the relevant undefinability of well-ordering number.  This fits into the project summarized in \cite{gv-superstable}: while superstability for arbitrary Abstract Elementary Classes seems poorly behaved (exhibiting what Shelah terms `schizophrenia' \cite[p. 19]{sh:h}), superstability in the context of tame Abstract Elementary Classes with amalgamation is much better behaved.  Vasey \cite{v-stabspec} computes stability spectra of Abstract Elementary Classes.  For tame classes with amalgamation, \cite[Corollary 4.24]{v-stabspec} uses a technical analysis of nonsplitting to show that failure of `$\bK$ is Galois stable on a tail' implies $\chi(\bK) > \omega$ and \cite[Corollary 4.17]{v-stabspec} shows that for `most $\lambda$,' $\cf \lambda < \chi(\bK)$ implies $\bK$ is Galois unstable in $\lambda$; `most $\lambda$' means all sufficiently large, almost $\lambda(\bK)$-closed cardinals.  Theorem \ref{unsuper-thm} offers a tighter bound on when the tail of stability must start and also a better condition on where the instability must happen.

Allowing for $\kappa \leq \LS(\bK)$, especially the case $\kappa = \omega$ is important because of the requirement that $\lambda^{<\kappa}=\lambda$ in (1) of Theorem \ref{unsuper-thm}.  This allows us to recover comparisons to first-order results and \cite{bs-atomic-stability}.

Our proof follows \cite{bs-atomic-stability}, but adapts the argument to Abstract Elementary Classes.  The following notion of type fragments will make our argument smoother.  These are essentially the partial Galois types that allow us to specify extending or not extending small Galois types.  This is motivated by the idea that small Galois types should occasionally be able to stand in for formulas in tame AECs (e.g., \cite[Section 3]{b-tamelc} or Vasey's Galois Morleyization \cite[Definition 3.3]{v-morleyization}).

\begin{hyp}\label{unsuper-hyp}
    In the rest of the section, we assume that $\bK$ is a $<\kappa$-tame Abstract Elementary class with amalgamation.  Write $\kappa^* = \kappa+\LS(\bK)$.
\end{hyp}

\begin{defin}\
\begin{enumerate}

	\item Given $M \in \bK$, $\cP^*_\kappa M : = \{ M_0 \in \bK_{<\kappa} : M_0 \prec M\}$.

	\item A \emph{$<\kappa$-(Galois) type fragment} over $B \subset N \in \bK$ is a collection $\Sigma$ of objects of the form `$p$' or `$\neg p$' with $p \in \S_\bK(A;N)$ for some $A \in \cP_\kappa B$.
	
	\item Some $a \in M$ \emph{realizes} a $<\kappa$-type fragment $\Sigma$ over $M$ iff $a \vDash p$ for all $p \in \Sigma$ and $a \not \vDash p$ for all $\neg p \in \Sigma$.
	
	\item A $<\kappa$-type fragment is \emph{satisfiable} iff some element realizes it.
	
\end{enumerate}
We won't have use for unsatisfiable type fragments, so all type fragments will be assumed to be satisfiable.
\end{defin}

We fix some important notation for this section: given $q\in \S_\bK(A; N)$, we write
\begin{eqnarray*}
q^0 := q \\
q^1 := \neg q
\end{eqnarray*}

In the following, we will want consider the number of types of length $<\kappa$ over the empty set.  We have 
$$\S_{\bK}^{<\kappa}(\emptyset) = \bigcup_{M \in \bK_{\kappa^*}} \S_{\bK}^{<\kappa}(\emptyset; M)$$
and can bound it's size by
$$|\S_{\bK}^{<\kappa}(\emptyset)| \leq \left(2^{\kappa^*}\right)^{<\kappa}$$
\comment{The following is standard:

\begin{prop}
If $\Sigma_n$ is a consistent $<\kappa$-type fragment over $M_n \in \bK$ for $n <\omega$ such that $M_n \prec M_{n+1}$ and $\Sigma_n \leq \Sigma_{n+1}$, then $\bigcup_{n<\omega}\Sigma_n$ is a $<\kappa$-type fragment over $\bigcup_{n<\omega} M_n$.
\end{prop}

\begin{lemma}\label{split-lem-old}
Suppose $N \in \bK$ and $\Gamma \subset \S_\bK(N)$ has size greater than $\|N\|^{<\kappa}$.  Then there is $M_1 \prec N$ of size $<\kappa$ and $r\neq q \in \S_\bK(M_1)$ such that both $q$ and $r$ have more than $\|N\|^{<\kappa}$-many extension to $\Gamma$.
\end{lemma}

{\bf Proof:}  If not, then, for every $M_1 \in \cP^*_\kappa N$, there is a unique $q_{M_1} \in \S_\bK(M_1)$ that has many extensions to $\Gamma$.  Then every $p \in \Gamma$ is of one of two kinds:
\begin{enumerate}
	\item $q \geq q_{M_0}$ for all $M_0 \in \cP^*_\kappa N$; or
	\item there is $M_0 \in \cP^*_\kappa N$ such that $q \not \geq q_{M_0}$.
\end{enumerate}
By tameness, there is at most one type in the first kind.  By counting, there are at most $|\cP^*_\kappa N|\times 2^{<\kappa}\times\|N\|^{<\kappa} = \|N\|^{<\kappa}$ in the second kind.  But $\Gamma$ is too large, so we have a contradiction.\hfill\dag\\

}

The following is similar to an argument of Baldwin-Kueker-VanDieren \cite[Claim 2.2]{bkv-stability} (see \cite[Theorem 11.11]{baldwinbook}), generalizing first-order arguments of Morley.  It will give us more than we need.

\begin{lemma}\label{split-lem}
Let $N \in \bK$ and $A \subset N$ with $\Gamma \subset \S_\bK(A; N)$ of size at least $\left(|A|^{<\kappa}\right)^+$.  Then there is $B \subset A$ of size $< \kappa$ with $q \neq r \in \S_\bK(B; N)$ such that both $q$ and $r$ have at least $(|A|^{<\kappa})^+$-many extensions to $\Gamma$.
\end{lemma}

{\bf Proof:} If not, then for every $B \in \cP_\kappa A$, there is a unique $q_B \in \S_\bK(B; N)$ that has many extensions to $\Gamma$.  Then every $p \in \Gamma$ falls into one of two categories:
\begin{enumerate}
    \item $p \geq q_B$ for all $B \in \cP_\kappa A$; or
    \item there is $B \in \cP_\kappa A$ such that $p \not \geq q_B$.
\end{enumerate}
By tameness, there is at most one type of the first kind.  For each $B$, there are $\leq |A|^{<\kappa}$-many types of the second kind by the choice of $q_B$ and there are $|A|^{<\kappa}$-many such $B$'s.  Thus, there are $\leq |A|^{<\kappa}$-many types of the second kind, which contradicts that there are at least $(|A|^{<\kappa})^+$-many types in $\Gamma$.\hfill \dag\\

\begin{lemma} \label{step1-lem}
Let $\mu > |\S_{\bK}^{<\kappa}(\emptyset)|$.  If $M \in \bK_{\geq 2^\mu}$ and $\Gamma = \seq{p_\alpha \in \S_\bK(M) : \alpha < \left(\|M\|^{<\kappa}\right)^+}$ are distinct,  then there is $\seq{A^i \in \cP_\kappa M : i < \mu}$ and $q(x; Y) \in \S_\bK^{<\kappa}(\emptyset)$ such that one of the following occur:
\begin{enumerate}
	\item for all $j_1 < \mu$, the following set has size $\left(\|M\|^{<\kappa}\right)^+$
	$$\left\{i < \left(\|M\|^{<\kappa}\right)^+ : q(x; A^{j_1}) \not\leq p_i \text{ and }j_0 < j_1 \text{ implies }q(x; A^{j_0}) \leq p_i\right\}$$
	\item for all $j_1 < \mu$, the following set has size $\left(\|M\|^{<\kappa}\right)^+$
	$$\left\{i < \left(\|M\|^{<\kappa}\right)^+ : q(x; A^{j_1}) \leq p_i \text{ and }j_0 < j_1 \text{ implies }q(x; A^{j_0}) \not\leq p_i\right\}$$
\end{enumerate}
\end{lemma}

{\bf Proof:}  We will construct
\begin{enumerate}
    \item a tree $T \subset {}^{\leq\mu} 2$;
    \item types $\{q_\eta(x; X) \in \S^{<\kappa}_\bK(\emptyset) : \eta \in T\}$; 
    \item sets $\{A^\eta \in \cP_\kappa M:\eta \in T\}$; and
    \item type fragments\footnote{Note that this type fragment is actually determined by the other objects in the construction.}
    $$\Sigma_\eta := \{q_{\eta \rest j}(x; A^{\eta\rest j})^{\eta(j)} : j < \ell(\eta)\}$$
    for $\eta \in T$
\end{enumerate}
such that for each $i \leq \mu$:
\begin{enumerate}
    \item every level $T_i$ of $T$ is nonempty (in particular, there is a branch $\eta \in T_\mu$);
	\item \label{frag-con} if $\eta \in T_i$, then the type fragment $\Sigma_{\eta}$ is contained in at least $\left(\|M\|^{<\kappa}\right)^+$-many of the types in $\Gamma$; and
	\item \label{split-cond} every node on level $T_i$ splits.
\end{enumerate}

{\bf This is enough:} Pick $\eta \in T_\mu$.  Since $\mu > |\S_{\bK}^{<\kappa}(\emptyset)|$, there is some $X \in [\mu]^\mu$, $q \in \S^{<\kappa}_\bK(\emptyset)$, and $k \in 2$ such that $q = q_{\eta\rest j}$ and $\eta(j) = k$ for all $j \in X$.  Write $\pi:X \to \mu$ for the Mostowski collapse and set 
$$A^i:=A^{\eta\rest \pi^{-1}(i)}$$
If $k = 0$, we are in the first case; if $k = 1$, we are in the second case. We give the details of the first case, and the second is similar.  Suppose $k=0$ and fix $j_1<\mu$.  We define the node just off the branch $\eta$ at height $\pi^{-1}(j_1)$
$$\nu = \left(\eta \rest \pi^{-1}(j_1)\right){}^\frown\seq{1} \in {}^{\pi^{-1}(j_1)+1} 2$$
We have that $\eta \rest \pi^{-1}(j_1) \in T$, so $\nu \in T$ by condition (\ref{split-cond}).  Thus, the fragment $\Sigma_\nu$ is contained in at least $(\|M\|^{<\kappa})^+$-many types in $\Gamma$ by condition (\ref{frag-con}).  The goal now is to show that the condition in item (1) of the lemma statement are contained in the fragment $\Sigma_\nu$.

If $j_0<j_1$, then $\eta\rest \pi^{-1}(j_0)  = \nu\rest\pi^{-1}(j_0)$, so 
$$q_{\eta\rest\pi^{-1}(j_0)}(x; A_{\eta \rest \pi^{-1}(j_0)})^{\eta(\pi^{-1}(j_0))} = q(x; A_{j_0})^0 = q(x; A_{j_0})\in \Sigma_\nu$$
For $j_1$, we similarly have
$$q_{\eta\rest\pi^{-1}(j_1)}(x; A_{\eta \rest \pi^{-1}(j_1)})^{\nu(\pi^{-1}(j_1))} = q(x; A_{j_1})^1 = \neg q(x; A_{j_1})\in \Sigma_\nu$$

{\bf Construction:} We work by induction on levels $i \leq \mu$.  At each, we will also guarantee that $\Sigma_{\eta^\frown\seq{0}}$ and $\Sigma_{\eta^\frown\seq{1}}$ satisfy condition (\ref{frag-con}) since they are defined at that stage.

For $i = 0$, we apply Lemma \ref{split-lem} to $\Gamma$ to find $A^\emptyset\in \cP_\kappa M$ and $q_\emptyset(x; A^\emptyset)\neq r(x;A^\emptyset) \in \S_\bK(A^\emptyset; M)$ such that both types extend to at least $(\|M\|^{<\kappa})^+$-many types in $\Gamma$.  Then we have that $\Sigma_{\seq{0}} = \{q_\emptyset(x; A^\emptyset)\}$ and $\Sigma_{\seq{1}} = \{ \neg q_\emptyset(x; A^\emptyset)\}$, with the latter satisfying our conditions because every type extending $r$ extends $\Sigma_{\seq{1}}$.

For $i = j+1$, for each $\eta \in T_i$ we follow the same strategy except we apply Lemma \ref{split-lem} to the elements of $\Gamma$ extending $\Sigma_\eta$.  This gives $A^\eta \in \cP_\kappa M$ and $q_\eta(x; A^\eta) \neq r(x;A^\eta) \in \S_\bK(A^\eta)$ that satisfy the requirements of the construction.

At limit stage $i$, we note that every type $p_\alpha$ (or more generally, every $p \in \S_\bK(M)$) extends one of our type fragments $\Sigma_\eta$ for $\eta \in {}^i 2$.  There are $\leq2^i$ many branches at this stage, and at least $\left(\|M\|^{<\kappa}\right)^+>2^\mu\geq 2^i$-many $p_\alpha$'s, so there must be some $\eta \in {}^i 2$ so many of them extend $\Sigma_\eta$; then set $T_i$ to be the collection of all such $\eta$.  Once again, we apply Lemma \ref{split-lem} to the elements of $\Gamma$-extending $\Sigma_\eta$ to define $A^\eta$ and $q_\eta$.\hfill \dag

The following lemma is the key inductive step that allows us to build our tree of types.

\begin{lemma}\label{split2-lem}
Fix $\mu > |\S^{<\kappa}_\bK(\emptyset)|+\LS(\bK)$. Suppose $M \in \bK_{\geq 2^\mu}$ with $|\S_\bK(M)| \geq \left(\|M\|^{<\kappa}\right)^+$, and let $\hat{M}$ be a $\mu^+$-Galois saturated extension of $M$.  

There are increasing $\{M_n \in \bK_{\mu} : n < \omega\}$; types $\{q_\nu \in \S_\bK(M_{\ell(\nu)}) : \nu \in {}^{<\omega} \mu\}$; sets $\{A^\nu \in \cP_\kappa M_{\ell(\nu)}:\nu \in {}^{<\omega} \mu\}$; and elements $\{a_\nu \in \hat{M} : \nu \in {}^{<\omega} \mu\}$ such that
\begin{enumerate}
    \item $M_n \prec_\bK M$ for all $n<\omega$;
	\item each $q_\nu$ has at least $(\|M\|^{<\kappa})^+$-many extensions to $\S_\bK(M)$ and $a_\nu$ realizes $q_\nu$;
	\item if $n < m$ and $\nu \in {}^m \mu$, then 
	$$q_{\nu \rest n} \leq q_\nu$$
	\item for every $\nu \in {}^{<\omega} \mu$ and $i < j < \mu$,
	$$q_{\nu^\frown\seq{i}} \rest A^{\nu^\frown\seq{i}} \neq q_{\nu^\frown\seq{j}}\rest A^{\nu^\frown\seq{i}}$$

\end{enumerate}
\end{lemma}

{\bf Proof:}  We do this by induction.  For the base case $n=0$, we pick $M_0 \prec_\bK M$ and $A^{\seq{}} \in \cP^*_\kappa M_0$ arbitrarily, then use the pigeonhole principle to find $q_{\seq{}} \in \S_\bK(M_0)$ with at least $\left(\|M\|^{<\kappa}\right)^+$-many extensions of $\S_{\bK}(M)$.

Given stage $n$, we know that each $q_\nu$ for $\nu \in {}^n \mu$ has at least $\left(\|M\|^{<\kappa}\right)^+$-many extensions to $\S_\bK(M)$ and $\|M\| \geq 2^\mu$.  So we can apply Lemma \ref{step1-lem} to get $q(x;Y) \in \S_\bK^{<\kappa}(\emptyset)$ and $\{A^i_\nu \in \cP_\kappa M: i < \mu\}$ and $\ell_\nu \in \{0,1\}$ such that
\begin{center}
for all $j_1 < \mu$, the following has size $\geq \left(\|M\|^{<\kappa}\right)^+$:
$$\{p \in \S_\bK(M) : q_\nu \leq p\text{ and }q(x; A_\nu^{j_1})^{1-\ell_\nu}\leq p \text{ and }j_0 < j_1 \text{ implies }q(x; A_\nu^{j_0})^{\ell_\nu} \not \leq p\}$$
\end{center}
Set $A_{\nu^\frown\seq{i}} = A^i_\nu$.

Let $M_{n+1} \prec M$ contain $M_n$ and $\bigcup_{\rho \in {}^{n+1}\mu} A_\rho$ of size $\mu$.  For each $i < \mu$, set
$$\Sigma'_{\nu, i}:= q_\nu \cup\{q(x; A^{i})^{1-\ell_\nu}, q(x; A^{j})^{\ell_\nu} : j < i\}$$
This is a consistent type fragment over $M$ by definition that can be extended to at least $\left(\|M\|^{<\kappa}\right)^+$-many types over $M$. Since there are at most $2^\mu$ extensions of $\Sigma'_{\nu, i}$ to $\S_{\bK}(M_{n+1})$ and $2^\mu < \left(\|M\|^{<\kappa}\right)^+$, the pigeonhole principle says we can extend $\Sigma'_{\nu, i}$ to a type $q_{\nu^\frown\seq{i}} \in \S_\bK(M_{n+1})$ that can be extended to at least $\left(\|M\|^{<\kappa}\right)^+$-many types over $M$.

By the saturation of $\hat{M}$, we can find $a_\nu \in \hat{M}$ realizing $p_\nu$ for each $\nu \in {}^{<\omega}\mu$.\hfill \dag\\

The following is an easy but useful fact before we begin the proof of the main theorem of this section.

\begin{prop}\label{sat-type-prop}
Suppose $\bK$ has amalgamation and $M \in \bK$ is $\mu^+$-Galois saturated.  If $M_0 \prec N_\ell \prec M$ of size $\leq \mu$ and $\ba_\ell \in N_\ell$ for $\ell = 0,1$ such that
$$\tp_\bK(\ba_0/M_0; N_0)=\tp_\bK(\ba_1/M_0; N_1)$$
then there is $N^* \in \bK_\mu$ such that $M_0 \prec N_0\prec N^*$ and $h:N_1\to_{M_0} N^*$ such that $h(\ba_1)=\ba_0$.
\end{prop}

Recall Shelah's Presentation Theorem (\cite{sh88} or \cite[Section 3.1]{bb-hanf} for a longer discussion).  Given an AEC $\bK$, it expands the lanuage by functions indexed by $\LS(\bK)\times\omega$ and gives an infinitary theory $T$ in the expanded langauge that captures both membership in $\bK$ and the strong substructure relation $\prec$.  However, it behaves like a Skolemization and can't necessarily capture every strong substructure relation simultaneously (although you can do so with a given chain).  In the following proof, we have a vast array of strong substructures that we want to code into the language by use of presentation functions.  We resolve this tension by introducing disjoint copies of this expanded language (often with parameters) to capture the different strong substructure relations we need.  We use the phrase `a collection of presentation functions witnessing $M\prec N$' to refer to a collection of functions (possibly with constant parameters) indexed by $\kappa^*\times \omega$ such that
\begin{itemize}
    \item the expansion by these functions satisfies the translation of the presentation theory $T$ to this language; and
    \item the models mentioned are closed under the functions.
\end{itemize}
This will ensure that any model of the translated theory (including generalized EM models, following the proof in Theorem \ref{gmott-ap-thm}) will satisfy the appropriate strong substructure relations.\\

{\bf  Proof of Theorem \ref{unsuper-thm}:} Recall $\kappa^* = \kappa +\LS(\bK)$ and Hypothesis \ref{unsuper-hyp}.  For the proof, suppose that (1) fails.  Then, for every $\alpha < (2^{\kappa_*})^+$, there is $M^\alpha \in \bK_{\geq \beth_{\alpha+1}}$ such that
$$|\S_\bK(M^\alpha)| \geq \left(\|M^\alpha\|^{<\kappa}\right)^+$$
We are going to build a blueprint $\Upsilon^{\omega-tr}[\bK]$ in the language
$$\tau_+:=\tau(\bK) \cup\{f_\gamma, H, F_{k, \beta}, G_{k, \beta}, G_{0, \beta, k}, G_{1, \beta, k}, h_0\mid \gamma < \kappa, k<\omega, \beta<\kappa_*\}$$
such that the $f_\gamma$ and $H$ are unary; the $G_{1, \beta, k}$ are binary; $h_0$ is ternary; the $F_{\beta, k}$ are $k$-ary; the $G_{\beta, k}$ are $(k+1)$-ary; and the $G_{0, \beta, k}$ are $(k+2)$-ary.

Now we build, for each $\alpha < (2^{\kappa_*})^+$, the $\tau_+$-structure $M^\alpha_+$ which will satisfy $M^\alpha \prec \left(M^\alpha_+\rest \tau(\bK)\right)$ (and much more).  Towards this end, fix $\alpha < (2^{\kappa_*})^+$ and set $\mu_\alpha = \beth_\alpha$; WLOG, $\mu_\alpha > |\S^{<\kappa}_{\bK}(\emptyset)|$.  We can use amalgamation to find $\hat{M}^\alpha \succ M^\alpha$ that is $\mu_\alpha^+$-saturated; this will be the restriction of $M^\alpha_+$ to $\tau(\bK)$.  Now we can apply Lemma \ref{split2-lem} to find $\prec$-increasing models $\{M^\alpha_n\in \bK_{\mu_\alpha} : n <\omega\}$; types $\{q^\alpha_\nu \in \S_\bK(M^\alpha_n): n<\omega, \nu \in{}^n \mu_\alpha\}$; sets $\{A^\alpha_\nu \in \cP_\kappa M^\alpha_n : n<\omega, \nu \in {}^n \mu_\alpha\}$; and elements $\{a^\alpha_\nu \in \hat{M}^\alpha : \nu \in {}^{<\omega}\mu_\alpha\}$ such that
\begin{enumerate}
    \item $M^\alpha_n \prec M^\alpha$ for all $n<\omega$;
    \item each $q^\alpha_\nu$ has at least $ (\|M^\alpha\|^{<\kappa})^{+}$-many extensions to $\S_\bK(M^\alpha)$ with $a^\alpha_\nu \vDash q^\alpha_\nu$;
    \item if $\nu < \eta$, then 
    $$q^\alpha_\nu \leq q^\alpha_\eta$$
    \item if $\nu\in {}^{<\omega}\mu$ and $i<j<\mu$, then
    $$q^\alpha_{\nu{}^\frown \{i\}} \rest A^\alpha_{\nu{}^\frown \{i\}} \ne q^\alpha_{\nu{}^\frown\{j\}}\rest A^\alpha_{\nu{}^\frown\{i\}}$$
\end{enumerate}
Now we will define some auxilliary models and functions.  For the most part, we are defining small witnesses for the various phenomena listed above, and then adding parameterized presentation functions to capture those structures in the blueprint.  When we define the expansions, we define them partially to capture the structure we want, and then implicitly expand their domain arbitrarily.

\begin{enumerate}
    \item For each $\nu \in {}^{<\omega}\mu_\alpha$, pick an arbitray element $b^\alpha_\nu \in M^\alpha_{\ell(\nu)}$.
    \item Define $\{f^\alpha_\gamma(b^\alpha_\nu):\gamma < \kappa\}$ to be an enumeration of $A^\alpha_\nu$ (with repetition); $H^\alpha(b^\alpha_\nu) = a^\alpha_\nu$; and $\{F^\alpha_{\beta, k}:\beta<\kappa^*, k<\omega\}$ be a collection of presentation functions such that
    \begin{enumerate}
        \item each $M^\alpha_n$ and $M^\alpha$ are closed under them (which witnesses $M^\alpha_n\prec M^\alpha \prec \hat{M}^\alpha$); and
        \item each $M^\alpha_n$ is the closure under $\{b^\alpha_\nu:\nu\in {}^n \mu_\alpha\}$ under these functions.
    \end{enumerate}
    \item For each $\nu \in {}^{<\omega}\mu_\alpha$, find $\bar{M}^\alpha_\nu \prec \hat{M}^\alpha$ of size $\mu_\alpha$ containing $a^\alpha_\nu$ and $M^\alpha_\ell(\nu)$.  Thus
    \begin{eqnarray*}
    q^\alpha_\nu = \tp_\bK(a^\alpha_\nu/M^\alpha_{\ell(\nu)}; \hat{M}^\alpha) = \tp_\bK(a^\alpha_\nu/M^\alpha_{\ell(\nu)}; \bar{M}^\alpha_\nu)
    \end{eqnarray*}
    Then set $\{G^\alpha_{\beta,k}(\bx, b^\alpha_\nu):\beta<\kappa^*, k<\omega\}$ to be presentation functions witnessing that $\bar{M}^\alpha_\nu \prec \hat{M}^\alpha$.
    \item For each $\nu < \eta \in {}^{<\omega}\mu_\alpha$, we have 
    \begin{eqnarray*}
        q^\alpha_\nu &=& q^\alpha_\eta\rest M^\alpha_{\ell(\nu)}\\
        \tp_\bK(a^\alpha_\nu/M^\alpha_{\ell(\nu)}; \bar{M}^\alpha_\nu) &=&\tp_\bK(a^\alpha_\eta/M^\alpha_{\ell(\nu)}; \bar{M}^\alpha_\eta)
    \end{eqnarray*}
    By Proposition \ref{sat-type-prop}, there is $\bar{M}_{0,(\nu, \eta)}\prec\hat{M}^\alpha$ of size $\mu_\alpha$ containing $a^\alpha_\nu$ and $M^\alpha_\nu$ and $h^\alpha_{0, (\nu, \eta)}: \bar{M}^\alpha_\eta \to\bar{M}^\alpha_{0,(\nu, \eta)}$ that fixes $M^\alpha_{\ell(\nu)}$ and so $h^\alpha_{0, (\nu, \eta)}(a^\alpha_\eta) = a^\alpha_\nu$.  Define 
    $$h^\alpha(x, b^\alpha_\nu, b^\alpha_\eta) := h^\alpha_{0, (\nu, \eta)}(x)$$
    and $\{G^\alpha_{0, \beta, k}(\bx, \eta, \nu):\beta<\kappa^*, k<\omega\}$ be presentation functions witnessing $h^\alpha_{0, (\nu, \eta)}(\bar{M}^\alpha_\eta)\prec\bar{M}^\alpha_{0, (\nu, \eta)} \prec \hat{M}^\alpha$ (note that $\bar{M}^\alpha_\nu\prec\bar{M}^\alpha_{0, (\nu, \eta)}$ will be guaranteed by the presentation functions $G^\alpha_{\beta,k}(\bx, b^\alpha_\nu)$ and coherence).
    \item For each $\nu \neq \eta \in {}^n\omega_\alpha$ with $\rho = \nu\rest(n-1)=\eta\rest(n-1)$ and $\nu(n-1)<\eta(n-1)$, we know that
    $$\tp_\bK(a^\alpha_\nu/A^\alpha_\nu; \hat{M}^\alpha)  \neq \tp_\bK(a^\alpha_\eta/A^\alpha_\nu; \hat{M}^\alpha)$$
    Fix $\bar{M}^\alpha_{1, (\eta, \nu)} \prec \hat{M}^\alpha$ of size $\kappa^*$ that contains $a^\alpha_\eta, a^\alpha_\nu$ and $A^\alpha_\nu$ and is closed under the functions $F^\alpha_{\beta,k}$ for $\beta<\kappa^*$, $k<\omega$.  Define $\{G^\alpha_{1,\beta, k}(b^\alpha_\nu, b^\alpha_\eta):\beta<\kappa^*, k<\omega\}$ to be an enumeration of $\bar{M}^\alpha_{1, (\nu, \eta)}$.
\end{enumerate}

Now we define the expansion $M^\alpha_+$ of $\hat{M}^\alpha$ to the language $\tau_+$: the expansion is the one suggested by the notation 
\begin{eqnarray*}
    f^{M^\alpha_+}_\gamma &:=& f^\alpha_\gamma\\
    H^{M^\alpha_+} &:=& H^\alpha\\
    &\text{etc.}&
\end{eqnarray*}
along with any necessary expansion of these to be full functions (rather than partial).

Now for each $\alpha <(2^{\kappa^*})^+$, we have a $\tau_+$-structure $M^\alpha_+$ along with an embedding 
\begin{eqnarray*}
{}^{<\omega}\beth_\alpha &\to& M^\alpha_+\\
\nu &\mapsto& b^\alpha_\nu
\end{eqnarray*}
Since $\cK^{\omega-tr}$ is a cofinal Erd\H{o}s-Rado Class (Corollary \ref{wtree-er-cor}), we can build a blueprint $\Phi \in \Upsilon^{\omega-tr}_{\kappa^*}[\bK]$ that is modeled off this embedding.  Given $\lambda\geq \kappa^*=\kappa+\LS(\bK)$, set $T^\lambda$ to be the $\cK^{\omega-tr}$-structure built on ${}^{<\omega} \lambda$ in the normal way.  We then read off the various structures we have built:
\begin{enumerate}
    \item $N^\lambda_+:= EM(T^\lambda, \Phi)$
    \item $\hat{N}^\lambda = EM_{\tau(\bK)}(T^\lambda, \Phi) \in \bK$
    \item $N^\lambda_n$ is the substructure of $N^\lambda$ generated by $\{\nu: \nu\in {}^n \lambda\}$ with $N^\lambda_n\prec \hat{N}^\lambda$
    \item $N^\lambda_\omega = \bigcup_{n<\omega} N^\lambda_n \prec \hat{N}^\lambda$
    \item For $\nu \in {}^{<\omega}\lambda$, set 
    $$B^\lambda_\nu:=\{f^{N^\lambda_+}_\gamma(\nu):\gamma<\kappa^*\}$$
\end{enumerate}
First, note that $\|N^\lambda_n\| = \lambda+\kappa^*=\lambda$ since it is generated by $\lambda$-many elements under $\kappa^*$-many functions. Then $\|N^\lambda_\omega\|=\lambda$ as well. Now we wish to construct our $\lambda^\omega$-many distinct types.  Unsurprisingly, they will come from the branches of the tree ${}^{<\omega}\lambda$.  As we do so, remember that the elements of ${}^{<\omega}\lambda$ are also elements of $N^\lambda_+$ and generate that model under the functions of $\Phi$.

For each $\rho \in {}^{<\omega}\lambda$, define
$$p^\lambda_\rho:=\tp_\bK\left(H^{N^\lambda_+}(\rho)/N_n^\lambda; \hat{N}^\lambda\right)$$
We have two key claims, both of which are guaranteed by coding the necessary information in the $\tau_+$-structure.

{\bf Claim 1:} Given $\nu < \eta \in {}^{<\omega} \lambda$, 
$$p^\lambda_\nu = p^\lambda_\eta \rest N^\lambda_{\ell(\nu)}$$
{\bf Proof:}. Define $\tau(\bK)$-structures as given by their universes below
\begin{eqnarray*}
    \bar{N}^\lambda_\nu&:=& \left\{G^{N^\lambda_+}_{\beta, k}\left(\bm, \nu\right):\beta<\kappa^*, k<\omega, \bm\in \left(\left\{H^\lambda(\nu)\right\}\cup N^\lambda_{\ell(\nu)}\right)\right\}\prec N^\lambda\\
    \bar{N}^\lambda_{0, (\nu, \eta)}&:=&\left\{G^{N^\lambda_+}_{0,\beta, k}\left(\bm, \nu, \eta\right):\beta<\kappa^*, k<\omega, \bm\in \left(\left\{H^\lambda(\nu)\right\}\cup N^\lambda_{\ell(\nu)}\right)\right\}\prec N^\lambda\\
    h^\lambda_{0, (\nu, \eta)}(x)&:=&h^{N^\lambda_+}(x, \nu, \eta)
\end{eqnarray*}
By the modeling property, since they are true at each $\alpha < (2^{\kappa^*})^+$, we have
\begin{enumerate}
    \item $N^\lambda_{\ell(\nu)} \prec \bar{N}^\lambda_\nu \prec \bar{N}^\lambda_{0, (\nu, \eta)}$
    \item $N^\lambda_{\ell(\nu)}\prec \bar{N}^\lambda_{\eta}$
    \item $h^\lambda_{0, (\nu, \eta)}: \bar{N}^\lambda_\eta \to \bar{N}^\lambda_{0, (\nu, \eta)}$ is a $\bK$-embedding that fixes $N^\lambda_\ell(\nu)$ and so $h^\lambda_{0, (\nu, \eta)}(H^{N^\lambda_+}(\eta)) = H^{N^\lambda_+}(\nu)$.
\end{enumerate}
This directly witnesses the type equality.

% https://q.uiver.app/#q=WzAsNSxbMSwzLCJOXlxcbGFtYmRhX3tcXGVsbChcXG51KX0iXSxbMCwyLCJcXGJhcntOfV5cXGxhbWJkYV9cXG51Il0sWzIsMiwiXFxiYXJ7Tn1eXFxsYW1iZGFfXFxldGEiXSxbMSwxLCJcXGJhcntOfV5cXGxhbWJkYV97MCwoXFxudSwgXFxldGEpfSJdLFsxLDAsIk5eXFxsYW1iZGFfKyJdLFswLDFdLFswLDJdLFsxLDNdLFsyLDMsIlxcYmFye2h9X3swLChcXG51LCBcXGV0YSl9IiwyXSxbMyw0XV0=
\[\begin{tikzcd}
	& {N^\lambda_+} & \\
	& {\bar{N}^\lambda_{0,(\nu, \eta)}} \\
	{\bar{N}^\lambda_\nu} && {\bar{N}^\lambda_\eta} \\
	& {N^\lambda_{\ell(\nu)}}
	\arrow[from=2-2, to=1-2]
	\arrow[from=3-1, to=2-2]
	\arrow["{\bar{h}_{0,(\nu, \eta)}}"', from=3-3, to=2-2]
	\arrow[from=4-2, to=3-1]
	\arrow[from=4-2, to=3-3]
\end{tikzcd}\]

{\bf Claim 2:} Given $\rho \in {}^{<\omega}\lambda$ and $i < j < \lambda$, write $\nu = \rho{}^\frown\{i\}$ and $\eta = \rho{}^\frown\{j\}$, so
$$p^\lambda_{\nu} \rest B^\lambda_\nu \neq p^\lambda_\eta\rest B^\lambda_\nu$$
{\bf Proof:} This claim is trickier than above since we have to transfer type inequality (which is a statement about lack of structure) rather than type equality (which is a statement about the existence of structure that can included in the language $\tau_+$). We can do this because the structures defining the types are generated by finitely many elements of $T^\lambda$, which in turn is possible since $\bK$ is $<\kappa$-tame with $\kappa\leq|\tau(\Phi)|$; we give the details.

Define the $\tau(\bK)$-structures given by the universe
$$N^\lambda_{1, (\nu, \eta)} =\left\{G^{N^\lambda_+}_{1,\beta, k}(\nu, \eta):\beta<\kappa^*, k<\omega\right\}\prec N^\lambda$$
Suppose for contradiction that 
$$p^\lambda_\nu \rest B^\lambda_\nu = p^\lambda_\eta\rest B^\lambda_\nu$$
This can be rewritten as
\begin{eqnarray*}
\tp_\bK\left(H^{N^\lambda_+}(\nu)/B^\lambda_\nu; N^\lambda_{1, (\nu, \eta)}\right) &=& \tp_\bK\left(H^{N^\lambda_+}(\eta)/B^\lambda_\nu; N^\lambda_{1, (\nu, \eta)}\right)\\
    \tp_\bK\left(H^{N^\lambda_+}(\nu)/\left\{f^{N^\lambda_+}_\gamma(\nu):\gamma<\kappa^*\right\}; \left\{G^{N^\lambda_+}_{1, \beta, k}(\nu, \eta):\beta<\kappa^*, k<\omega\right\}\right)\\&=&\\\tp_\bK\left(H^{N^\lambda_+}(\eta)/\left\{f^{N^\lambda_+}_\gamma(\nu):\gamma<\kappa^*\right\}; \left\{G^{N^\lambda_+}_{1, \beta, k}(\nu, \eta):\beta<\kappa^*, k<\omega\right\}\right)
\end{eqnarray*}
If these types are equal, then there are $N^+ \in \bK_{\kappa^*}$ and $f_\ell:N^\lambda_{1,(\nu, \eta)}\to N^+$ for $\ell = 0,1$ such that $f_0 \rest B^\lambda_\nu = f_1\rest B^\lambda_\nu$ and
$$f_0\left(H^{N^\lambda_+}(\nu)\right)=f_0\left(H^{N^\lambda_+}(\eta)\right)$$
By the modeling property (recall Definition \ref{er-def}.(\ref{erc2-item})), there is $\alpha < (2^{\kappa^*})^+$ and $\nu', \eta' \in {}^{<\omega}\mu_\alpha$ such that
\begin{eqnarray*}
\tp_{\omega-tr}(\nu', \eta'; T^{\mu_\alpha}) &=& \tp_{\omega-tr}(\nu, \eta; T^\lambda)\\
\tp_{\tau_+}(b^\alpha_{\nu'},b^\alpha_{\eta'}; M^\alpha_+) &=& \tp_{\tau_+}(\nu, \eta; N^\lambda_+)
\end{eqnarray*}
In particular, this $\tau_+$-quantifier free type includes the isomorphism type of the model $N^\lambda_{1, (\lambda, \eta)}$.  Thus, the map
$$G^{N^\lambda_+}_{1, \beta, k}(\nu, \eta) \mapsto G^{M^\alpha_+}_{1, \beta, k}(\nu', \eta')$$
is an isomorphism $g:N^\lambda_{1, (\nu, \eta)} \cong M^\alpha_{1,(\nu', \eta')}$ that sends $A^\alpha_{\nu'}$ to $B^\lambda_\nu$ and sends $H^{N^\lambda_+}(\nu), H^{N^\lambda_+}(\eta)$ to $a^\alpha_{\nu'}, a^\alpha_{\eta'}$.  Thus, composing everything, we have
% https://q.uiver.app/#q=WzAsNixbMCwyLCJBXlxcYWxwaGFfe1xcbnUnfSJdLFswLDEsIk1eXFxhbHBoYV97MSwoXFxudScsXFxldGEnKX0iXSxbMSwyLCJNXlxcYWxwaGFfezEsKFxcbnUnLCBcXGV0YScpfSJdLFsxLDAsIk5eXFxsYW1iZGFfezEsIChcXG51LCBcXGV0YSl9Il0sWzIsMSwiTl5cXGxhbWJkYV97MSwgKFxcbnUsXFxldGEpfSJdLFsyLDAsIk5eKyJdLFswLDFdLFswLDJdLFsxLDMsImdeey0xfSJdLFsyLDQsImdeey0xfSIsMl0sWzQsNSwiZl8xIiwyXSxbMyw1LCJmXzAiXV0=
\[\begin{tikzcd}
	& {N^\lambda_{1, (\nu, \eta)}} & {N^+} \\
	{M^\alpha_{1,(\nu',\eta')}} && {N^\lambda_{1, (\nu,\eta)}} \\
	{A^\alpha_{\nu'}} & {M^\alpha_{1,(\nu', \eta')}}
	\arrow["{f_0}", from=1-2, to=1-3]
	\arrow["{g^{-1}}", from=2-1, to=1-2]
	\arrow["{f_1}"', from=2-3, to=1-3]
	\arrow[from=3-1, to=2-1]
	\arrow[from=3-1, to=3-2]
	\arrow["{g^{-1}}"', from=3-2, to=2-3]
\end{tikzcd}\]

so that the maps agree on $A^\alpha_{\nu'}$ and send $a^\alpha_{\nu'}$ to $a^\alpha_{\eta'}$.  This gives
$$\tp_\bK(a^\alpha_{\nu'}/A^\alpha_{\nu'}; M^{\alpha}_{1, (\nu', \eta')})=\tp_\bK(a^\alpha_{\eta'}/A^\alpha_{\eta'};M^\alpha_{1,(\nu', \eta')})$$
However, since $\nu, \eta$ and $\nu', \eta'$ have the same type, we have that $\rho'=\nu'\rest(n-1)=\eta'\rest(n-1)$ and $\nu'(n-1)<\eta'(n-1)$. Thus, by construction, we have
$$\tp_\bK(a^\alpha_{\nu'}/A^\alpha_{\nu'}; \hat{M}^{\alpha})=\tp_\bK(a^\alpha_{\eta'}/A^\alpha_{\eta'};\hat{M}^\alpha)$$
which is a contradiction since $M^\alpha_{1,(\nu', \eta')} \prec \hat{M}^\alpha$.\hfill \dag${}_{\text{Claim 2}}$\\

Putting these together completes the proof.

{\bf Claim 3:} $|\gS_\bK(N^\lambda_\omega)| \geq \lambda^\omega$\\
{\bf Proof:} For each $\eta \in {}^{\omega}\lambda$, the sequence $\{p^\lambda_{\eta\rest n}:n<\omega\}$ is an increasing chain by Claim 1.  By the $\omega$-compactness of AECs with amalgamation (see \cite[Theorem 11.1]{baldwinbook}), there is $p_\eta \in \gS_\bK(N^\lambda_\omega)$ that extends all of them.  By Claim 2, each $p_\eta$ is distinct.\hfill\dag\\

\subsection{Indiscernible Collapse in Nonelementary Classes} \label{ind-coll-ssec}

One of the uses of generalized indiscernibles in first-order is to characterize various dividing lines via indiscernible collapse.  An old result of Shelah \cite[Theorem II.2.13]{sh:c} says that a theory $T$ is stable iff any order indiscernibles in a model of $T$ are in fact set indiscernibles.  Scow \cite[Theorem 5.11]{s-nip-collapse} proved that $T$ is NIP iff any ordered graph indiscernibles in a model of $T$ are in fact order indiscernibles.  In each of these cases, there are abstract (Ramsey) classes $\cK_0$ and $\cK$ with $\cK_0$ a reduct of $\cK$ where some property of $T$ can be detected by whether or not there are $\cK$-indiscernibles that are not $\cK_0$-indiscernibles (after reducting the index). Guingona, Hill, and Scow \cite{ghs-collapse} have formalized this notion of indiscernible collapse and given several more examples.

Following this work, we can give definitions of several dividing lines in Abstract Elementary Classes making use of the fact that the determining classes are Erd\H{o}s-Rado classes in addition to being Ramsey classes.  Unfortunately, at this time, we don't know of any indiscernible collapses characterizing dividing lines that start with an Erd\H{o}s-Rado class (other than order indiscernibles, but this collapse result is already known for Abstract Elementary Classes\footnotei{Citation?}).  \cite[Theorem 3.4]{ghs-collapse} uses $\cK^{n-mlo}$ and \cite[Theorem 4.7]{ghs-collapse} uses a class of trees that doesn't restrict the height (in a similar way that $\cK^{ceq}$ generalizes $\cK^{\chi-or}$), but neither of these are known to be Erd\H{o}s-Rado classes.  \cite[Corollary 5.9]{ghs-collapse} characterizes $NTP_2$ theories via a collapse of $\cK^{ceq}$-indiscernibles, but involves notions of formulas dividing that does not easily generalize to Abstract Elementary Classes.  This leads us to the following question:

\begin{question}
Is $\cK^{og}$ an Erd\H{o}s-Rado class?
\end{question}

Recall that it is consistently not a combinatorial Erd\H{o}s-Rado class by Example \ref{og-spr-ex}.  However, this does not rule out the possibility it is an Erd\H{o}s-Rado class.  A positive answer for this question would give a prospective definition for the notion of NIP for Abstract Elementary Classes.

\begin{defin}
Suppose that $\cK^{og}$ is an Erd\H{o}s-Rado class and let $\bK$ be an Abstract Elementary Class with arbitrarily large models.  We say that $\bK$ is \emph{NIP} iff for every $\Phi \in \Upsilon^{og}[\bK]$, there is $\Psi \in \Upsilon^{or}[\bK]$ such that $\Phi = \Psi \circ U$, where $U \in \Upsilon^{og}[\cK^{or}]$ forgets the graph structure and $\Psi \circ U$ is the composition of these blueprints.
\end{defin}

This has advantages over other prospective definitions\footnotei{(see [\footnote{Ask Andres}] for some discussion)} in that no amalgamation, tameness, etc. assumption is necessary.   Of course, it has the disadvantage that it needs more results to be viable.  This is being explored further in \cite{bs-erc-plus}

\subsubsection{Category theoretic interpretation of indiscernible collapse}

We can build on the category theoretic interpretation of indiscernibles from Section \ref{ct-ind-ssec} to give the same gloss to indiscernible collapse in terms of injectivity conditions (\cite[Section 4.A]{ar-presentable} gives this background).  In general, if $f:A \to B$ is a morphism, then another object $C$ is injective with respect to $f$ iff every $g:A \to C$ can be lifted along $f$ to a $g':B \to C$ so $g = g' \circ f$.  Then Shelah's result \cite[Theorem II.2.13]{sh:c} can be rephrased as follows.

\begin{theorem}
Let $\cK^{set}$ be the abstract class of sets and $U:\cK^{or} \to \cK^{set}$ be the functor forgetting the ordering.  An elementary class $\bK$ is stable iff it is injective with respect to $U$ (in the category of accessible categories whose morphisms are faithful functors preserving directed colimits).
\end{theorem}

{\bf Proof:}  First, assume $\bK$ is stable and let $G:\cK^{or} \to \bK$.  By Theorem \ref{der-bp-thm}, we may assume that $G$ comes from a blueprint $\Phi$ for order-indiscernibles.  By \cite[Theorem II.2.13]{sh:c}, $\Phi$ gives rise to a blueprint $\Psi$ for set indiscernibles.  Then, using Theorem \ref{der-bp-thm} again, $\Psi$ gives rise to the desired $G'$.

Second, suppose that $\bK$ is injective in this sense.  Then any order indiscernibles are set indiscernibles.  By \cite[Theorem II.2.13]{sh:c}, $\bK$ is stable.\hfill \dag\\

Other indiscernible collapses can be phrased similarly.

\subsection{Interpretability Order}

Finding generalized indiscernibles in nonelementary classes can give stronger negative results in the interpretability order even for comparing first-order theories.  The interpretability order is a three-parameter order $\triangleleft^*_{\lambda, \chi, \kappa}$ on complete first-order theories introduced by Shelah \cite{sh500} in the vein of Keisler's order.  It would say that $T_0$ is less complicated than $T_1$ iff every time a first-order theory interprets both $T_0$ and $T_1$, if the interpretation of $T_1$ is saturated, then so is the interpretation of $T_0$.  \cite[Definition 2.10]{sh500} gives the full definition, but we only need a particular instance, $\triangleleft^*_1$, which we weaken to $\triangleleft^{*, \kappa}_1$.  Moreover, we allow the $\chi$-parameter to be arbitrarily large (rather than countable as in most instances in \cite{ms-interp}) to strengthen our results.  Since this application is not central, we omit some of the definitions, but a good exposition (and the results we reference) can be found in Malliaris and Shelah \cite{ms-interp}.

\begin{defin}
Let $T_0$ and $T_1$ be complete first-order theories and let $\mu$  be an infinite cardinal.  
\begin{enumerate}
	\item We say that $T_0 \triangleleft^*_1 T_1$ iff for all large enough, regular $\mu$, there is a first-order theory $T_*$ of size $\leq |T_0|+|T_1|+\aleph_0$ that interprets $T_\ell$ via $\bar{\phi_\ell}$ such that, for every $M_* \vDash T_*$, if the interpretation $M_*^{[\bar{\phi_1}]}$ of $T_1$ is $\mu$-saturated, then $M_*^{[\bar{\phi_0}]}$ is $\mu$-saturated. (Shelah)
	\item We say that $T_0 \triangleleft^{*,\kappa}_1 T_1$ iff for all large enough, regular $\mu$, there is an $\bL_{\kappa, \omega}$-theory $T_*$ of size $\leq |T_0|+|T_1|+\aleph_0$ that interprets $T_\ell$ via $\bar{\phi_\ell}$ such that, for every $M_* \vDash T_*$, if the interpretation $M_*^{[\bar{\phi_1}]}$ of $T_1$ is $\mu$-saturated, then $M_*^{[\bar{\phi_0}]}$ is $\mu$-saturated.
\end{enumerate}
\end{defin}

So $\triangleleft^{*, \kappa}_1$ differs from $\triangleleft^*_1$ in that it allows for infinitary theories to do the interpreting.  In particular, the statement that $\neg(T_0 \triangleleft^{*,\kappa}_1 T_1)$ is a stronger statement than $\neg(T_0 \triangleleft^*_1 T_1)$.  In \cite{ms-interp}, Malliaris and Shelah show several positive and negative instances of the interpretability order.  The negative instances are proved by using various Ramsey classes to build generalized blueprints that saturate $T_1$ without saturating $T_0$.  When these Ramsey classes are in fact Erd\H{o}s-Rado classes, the stronger negative instance can be shown.  In the following statement, $T_{DLO}$ is the theory of dense linear orders and $T_{RG}$ is the theory of the random graph.

\begin{fact}[{\cite[Observation 5.7]{ms-interp}}]
$\neg (T_{DLO} \triangleleft^*_1 T_{RG})$
\end{fact}

\begin{theorem}
For every cardinal $\kappa$, $\neg (T_{DLO} \triangleleft^{*, \kappa}_1 T_{RG})$.
\end{theorem}

{\bf Proof:} We rely heavily on citations from \cite{ms-interp}.  Note their $\cK = \cK_{\lambda}$ is essentially our $\cK^{\lambda-color}$, which is Erd\H{o}s-Rado by Example \ref{clo-spr-ex} and Corollary \ref{cer-er-cor}.  Also, they use $GEM$ to emphasize that the Ehrenfeucht-Mostowski construction uses a generalized blueprint.  We adopt \cite[Hypotheses 5.6 and 5.8]{ms-interp} with our infinitary change, so 
\begin{enumerate}
	\item $\lambda = \lambda^{<\mu} \geq 2^\mu$;
	\item $T_*$ is a skolemized $\bL_{\kappa, \omega}$-theory with $|T_*| \leq \lambda$ that interprets $T_{RG}$ by $R_{RG}$ and interprets $T_{DLO}$ by $<_{DLO}$.
\end{enumerate}
Note that they point out that their results in this area work for uncountable languages as well.

Since $\cK$ is a combinatorial Erd\H{o}s-Rado class, there is $\Phi \in\Upsilon^{\cK}[T_*]$.  By \cite[Corollary 5.11]{ms-interp}, we can find $\Psi$ extending $\Phi$ such that for every separated $I \in \cK$, $EM_{RG}(I, \Psi)$ is $\mu$-saturated.  Note that, since $\Psi$ agrees with $\Phi$ on $\tau(T_*)$, $\Psi$ is still in $\Upsilon^{\cK}[T_*]$.  By \cite[Claim 5.12]{ms-interp}, if $J$ is a separated linear order, then for any $\Phi^* \in \Upsilon^{\cK}[T_*]$, $EM_{DLO}(J, \Phi^*)$ is not $\kappa^+$-saturated.  Thus, by taking $I$ separated with a $(\kappa, \kappa)$-cut, we have $EM_{RG}(I, \Psi)$ is $\mu$-saturated, but $EM_{DLO}(I, \Psi)$ is not $\kappa^+$-saturated, as desired.\hfill\dag\\

\bibliographystyle{amsalpha}
\bibliography{bib}

\end{document}